\pdfoutput=1

\documentclass[titlepage,twoside,openright,10pt]{book}

\usepackage{amsthm,mathrsfs,bbm,amssymb,here,amsmath,cite,textcomp,stmaryrd}
\usepackage[a5paper, top=.8in, bottom=1.3in, left=.5in, right=.5in]{geometry}
\usepackage{fancyhdr}
\usepackage{makeidx}
\usepackage{graphicx}
\usepackage{xr}

\makeindex

\newcommand{\be}{\begin{equation}}
\newcommand{\ee}{\end{equation}}
\newcommand{\bea}{\begin{eqnarray}}
\newcommand{\eea}{\end{eqnarray}}
\newcommand{\ba}{\begin{array}}
\newcommand{\ea}{\end{array}}
\newcommand{\mc}{\mathcal}

\newcommand{\eopsymbol}{$\blacksquare$}
\newcommand{\naturalnumbers}{{\mathbb N}}
\newcommand{\bigtimes}{\raisebox{-.2em}{\huge $\times$}}

\newcommand{\midvspace}{\vspace{.25in}}

\newtheoremstyle{alltheorem}
	{0in}	
	{0in}	
	{}		
	{}		
	{\bf}	
	{}		
	{\newline}	
	{}		

\theoremstyle{alltheorem}
\newtheorem{defi}{Definition}[chapter]
\newtheorem{prop}[defi]{Proposition}
\newtheorem{rema}[defi]{Remark}
\newtheorem{axio}[defi]{Axiom}
\newtheorem{lemm}[defi]{Lemma}
\newtheorem{theo}[defi]{Theorem}
\newtheorem{coro}[defi]{Corollary}
\newtheorem{exam}[defi]{Example}
\newtheorem{lede}[defi]{Lemma and Definition}

\newcommand{\alltheoremwidth}{\textwidth}

\newcommand{\bdefi}{
\begin{minipage}{\alltheoremwidth}
\begin{defi}}

\newcommand{\edefi}{
\hspace{\stretch{1}} \eopsymbol \end{defi} \end{minipage} \vspace{.25in}}

\newcommand{\bexam}{
\begin{minipage}{\alltheoremwidth}
\begin{exam}}

\newcommand{\eexam}{
\hspace{\stretch{1}} \eopsymbol \end{exam} \end{minipage} \vspace{.25in}}

\newcommand{\brema}{
\begin{minipage}{\alltheoremwidth}
\begin{rema}}

\newcommand{\erema}{
\hspace{\stretch{1}} \eopsymbol \end{rema} \end{minipage} \vspace{.25in}}

\newcommand{\baxio}{
\begin{minipage}{\alltheoremwidth}
\begin{axio}}

\newcommand{\eaxio}{
\hspace{\stretch{1}} \eopsymbol \end{axio} \end{minipage} \vspace{.25in}}

\newcommand{\blede}{
\begin{minipage}{\alltheoremwidth}
\begin{lede}}

\newcommand{\elede}{
\end{lede} \end{minipage} \vspace{.0in}}

\newcommand{\blemm}{
\begin{minipage}{\alltheoremwidth}
\begin{lemm}}

\newcommand{\elemm}{
\end{lemm} \end{minipage} \vspace{.0in}}

\newcommand{\btheo}{
\begin{minipage}{\alltheoremwidth}
\begin{theo}}

\newcommand{\etheo}{
\end{theo} \end{minipage} \vspace{.0in}}

\newcommand{\bprop}{
\begin{minipage}{\alltheoremwidth}
\begin{prop}}

\newcommand{\eprop}{
\end{prop} \end{minipage} \vspace{.0in}}

\newcommand{\bcoro}{
\begin{minipage}{\alltheoremwidth}
\begin{coro}}

\newcommand{\ecoro}{
\end{coro} \end{minipage} \vspace{.0in}}

\newcommand{\bproof}{
\vspace{-.05in} \begin{proof}}

\newcommand{\eproof}{
\end{proof} \vspace{.15in}}

\newcommand{\benum}{
\renewcommand{\theenumi}{\roman{enumi}}
\renewcommand{\labelenumi}{(\theenumi)}
\vspace{-.05in} \begin{enumerate}}

\newcommand{\eenum}{\end{enumerate}}

\usepackage[pdftex]{hyperref}

\fancypagestyle{plain}{
\lhead[\thepage]{}
\chead{}
\rhead[]{\thepage}
\lfoot{\ \\ \small \textcopyright \ 2013 \, Felix Nagel --- Set theory and topology, Part I: Sets, relations, numbers}
\cfoot{}
\rfoot{}
}

\begin{document}

\frontmatter


\thispagestyle{empty}

\vspace{\baselineskip}

\begin{center}
\textbf{\huge Set theory and topology\\[.3em]
\large An introduction to the foundations of analysis
\footnote{For remarks and suggestions please contact: stt.info@t-online.de}\\[3em]
\large Part I: \; Sets, relations, numbers}\\
\vspace{7\baselineskip}
{\sc Felix~Nagel}\\
\vspace{5\baselineskip}
\textbf{Abstract}\\
\vspace{1\baselineskip}
\parbox{0.9\textwidth}{We provide a formal introduction into the classic theorems of general topology and its axiomatic foundations in set theory. Starting from ZFC, the exposition in this first part includes relation and order theory as well as a construction of number systems.}
\end{center}
\vspace{\baselineskip}

\pagebreak

\thispagestyle{empty}

{\small The author received his doctoral degree from the University of Heidelberg for his thesis in electroweak gauge theory. He worked for several years as a financial engineer in the financial industry. His fields of interest are probability theory, foundations of analysis, finance, and mathematical physics. He lives in Wales and Lower Saxony.}


\chapter{Preface}
\label{sec-preface}

\setcounter{page}{1}

\thispagestyle{plain}

This series of articles emerged from the author's personal notes on general topology supplemented by an axiomatic construction of number systems.

At the beginning of the text we introduce our axioms of set theory, from which all results are subsequently derived. In this way the theory is developed {\it ab ovo} and we do not refer to the literature in any of the proofs.

Needless to say, the presented theory is fundamental to many fields of mathematics like linear analysis, measure theory, probability theory, and theory of partial differential equations. It establishes the notions of relation, function, sequence, net, filter, convergence, pseudo-metric and metric, continuity, uniform continuity etc. To derive the most important classic theorems of general topology is the main goal of the text. Additionally, number systems are studied because, first, important issues in topology are related to real numbers, for instance pseudo-metrics where reals are required at the point of the basic definitions. Second, many interesting examples involve numbers.

In our exposition we particularly put emphasis on the following:

\benum
\item The two advanced concepts of convergence, viz.\ nets and filters, are treated with almost equal weighting. Most results are presented both in terms of nets and in terms of filters. We use one concept whenever it seems more appropriate than the other.
\item We avoid the definition of functions on the ensemble of all sets. Instead we follow a more conservative approach by first choosing an appropriate set in each case on which the respective analysis is based. Notably, this issue occurs in the Recursion theorem for natural numbers (see Theorem~\ref{theo recursive def} where, with this restriction, also the Replacement schema is not required), in the Induction principle for ordinal numbers (see Theorem~\ref{theo induction ordinals}), and in the Local recursion theorem for ordinal numbers (see Theorem~\ref{theo recursion}).
\item At many places we try to be as general as possible. In particular, when we analyse relations, many definitions and results are stated in terms of pre-orderings, which we only require to be transitive.
\eenum

\thispagestyle{plain}

Finally, we would like to warn the reader that for some notions defined in this work there are many differing definitions and notations in the literature, e.g.\ in the context of relations and orderings. One should always look at the basic definitions before comparing the results.

The text is structured as follows: All definitions occur in the paragraphs explicitly named {\bf Definition}. Important Theorems are named {\bf Theorem}, less important ones~{\bf Lemma}, though a distinction seems more or less arbitrary in many cases. In some cases a Lemma and a Definition occur in the same paragraph in order to avoid repetition. Such a paragraph is called {\bf Lemma and Definition}. Claims that lead to a Theorem and are separately stated and proven are called {\bf Proposition}, those derived from Theorems are called {\bf Corollary}. The proofs of most Theorems, Lemmas, Propositions, and Corollaries are given. Within lengthier proofs, intermediate steps are sometimes indented and put in square brackets $\left[ \ldots \right]$ in order to make the general outline transparent while still explaining every step. Some proofs are left as excercise to the reader. There is no type of paragraph explicitly named as exercise. Paragraphs named {\bf Example} contain specializations of Definitions, Theorems, etc.\ The analysis of the examples is mostly left to the reader without explicit mention. Statements that do not require extensive proofs and are yet relevant on their own are named~{\bf Remark}. 

Note that Definitions, Theorems etc.\ are enumerated per Chapter. Some references refer to Chapters that are contained in subsequent parts of this work~\cite{Nagel}.
\\

\begin{tabular*}{\textwidth}{l @{\extracolsep{\fill}} r}
{\it Wales, May 2013} & {\it Felix Nagel}
\end{tabular*}

\thispagestyle{plain}

\pagebreak

\thispagestyle{plain}


\cleardoublepage

\phantomsection
\thispagestyle{plain}
\addcontentsline{toc}{chapter}{Table of contents}
\tableofcontents
\thispagestyle{plain}

\mainmatter


\part{Sets, relations, numbers}
\label{part sets relations numbers}

\chapter{Axiomatic foundation}
\label{axiomatic foundation}
\setcounter{equation}{0}

\pagebreak

In our exposition, as is the case with every mathematical text, we do not solely use verbal expressions but need a formal mathematical language. To begin with, let us describe which role the formal mathematical language is supposed to play subsequently.

We initially define certain elementary notions of the formal language by means of ordinary language. This is done in Section~\ref{formal language}. First, we define logical symbols, e.g.\ the symbol $\Longrightarrow$, which stands for "implies", and the symbol $=$ meaning "equals". Second, we define the meaning of set variables. Every set variable, e.g.\ the capital letter $X$ of the Latin alphabet, stands for a set. A set has to be interpreted as an abstract mathematical object that has no properties apart from those stated in the theory. Third, we define the symbol~$\in$, by which we express that a set is an element of a set. In all three cases the correct interpretation of the symbols comprises, on the one hand, to understand its correct meaning and, on the other hand, not to associate more with it than this pure abstract meaning. A fourth kind of elementary formal component is used occasionally. We sometimes use Greek letters, e.g.\ $\varphi$, as variables that stand for a certain type of formal mathematical expressions called formulae. Such formula variables belong to the elementary parts of our formal language because a formula variable may not only be an abbreviation for a specific formula in order to abridge the exposition but is also used as placeholder for statements that we make about more than a single formula. The latter is tantamount to an abbreviation if a finite number of formulae is supposed to be substituted but cannot be regarded as a mere abbreviation if an infinite number of formulae is considered.

After having translated these elementary mathematical thoughts into the formal language in Section~\ref{formal language}, we may form formulae out of the symbols. This allows us to express more complicated mathematical statements in our formal language. We follow the rationale that every axiom, definition, and claim in the remainder of this text can in principle be expressed either in the formal language or in the nonformal language and translated in both directions without ambiguity.
In practice, in some cases the formal language and in other cases the nonformal language are preferable with respect to legibility and brevity. Therefore we make use of both languages, even mix both deliberately. For example, in order to obtain a precise understanding of our axioms of set theory, we tend to use only the formal language in this context. In most other cases we partly use the formal language to form mathematical expressions, which are then surrounded by elements of the nonformal language. In the field of mathematical logic such formal languages and their interpretations in the nonformal language are probed. There also the nature of mathematical proofs, i.e.\ the derivation of theorems from assumptions, is analysed. Ways to formalize proofs are proposed so that in principle the derivation of a theorem could be written in a formal language. We do not formalize our proofs in this way but use the fomal language only to the extent described above.

\section{Formal language}
\label{formal language}

The only objects we consider are sets. Statements about sets are written in our formal language as formulae that consist of certain symbols. We distinguish between symbols that have a fixed meaning wherever they occur, variables that stand for sets, and variables for formulae, which may be substituted if a specific formula is meant. The symbols that have a fixed meaning are the symbol~$\in$ and various logical symbols including~$=$\,. As announced in the introduction to this Chapter we now define all these symbols by their meaning in the nonformal language and explain the meaning of variables that stand for sets and those that stand for formulae. We remark that in the remainder of the text further symbols with fixed meaning on a less elementary level are defined as abbreviations.

A variable that denotes a set is called a set variable\index{Variable!set}\index{Set variable}. We may use as set variables any small or capital letters of the Latin or Greek alphabet with or without subscripts, superscripts, or other add-ons.

Let $x$ and $y$ be set variables. We write the fact that $x$ is an element of~$y$ in our formal language as $(x \in y)$\index{Element}\index{Set!element}. In this case we also say that $x$ is a member of~$y$\index{Member}\index{Set!member}. Furthermore we write the fact that $x$ equals~$y$, i.e.\ $x$ and $y$ denote the same set, as $(x = y)$\index{Equality}\index{Set!equality} in our formal language. Similar translations from the nonformal language to the formal language are applied for membership and equality of any two variables different from $x$ and~$y$.

If $x$ and $y$ are set variables, each of the expressions $(x \in y)$ and $(x = y)$ is called an atomic formula\index{Formula!atomic}\index{Atomic}. The same holds for all such expressions containing set variables different from $x$ and~$y$.

Generally, a formula may contain, apart from variables denoting sets, the symbols $\in$ and~$=$, the logical connectives\index{Connective!logical}\index{Logical connective} $\wedge$ (conjunction, and), $\vee$ (disjunction, non-exclusive or), $\neg$ (negation), $\Longrightarrow$ (implication), $\Longleftrightarrow$ (equivalence), and the quantifiers $\forall$ (for every) and $\exists$ (there exists). All these logical symbols are used in their conventional nonformal interpretation indicated after each symbol above. Additionally, the brackets $($ and $)$ are used in order to express in which order a formula has to be read. Some of the symbols are clearly redundant to express our nonformal thoughts. For instance, if we use the symbols $\neg$ and $\vee$, then $\wedge$ is not required. Or, if we use the symbols $\neg$ and  $\exists$, the symbol $\forall$ is redundant. However, it is often convenient to make use of all logical symbols. As the meaning of all these symbols is defined in our nonformal language, it is clear, for each expression that is written down using these symbols, whether such an expression is meaningful or not. For example, the expressions

\benum
\item  \label{meaningful ex 1} $\big( \neg (x \in y) \big) \Longrightarrow (\exists z \; z \in x)$
\item \label{meaningful ex 2} $\exists x \; (x \in x)$
\item \label{meaningful ex 3} $\big( (x = y) \wedge (z \in x) \big) \;\; \Longleftrightarrow \;\; \big( (x \in z) \vee (z = z) \big)$
\item \label{meaningful ex 4} $(x = y) \wedge \neg (x = y)$
\eenum

where $x$, $y$, and $z$ are set variables, are meaningful, though some may be logically false in a specific context, as e.g.~(\ref{meaningful ex 2}) in the theory presented in this text, or even false in any context like~(\ref{meaningful ex 4}).

In contrast, the expressions

\benum
\item \label{ex 1} $\neg x$
\item \label{ex 2} $\forall \forall z$
\item \label{ex 3} $\exists x \; \big( (x \in y) \vee \big)$
\item \label{ex 4} $\exists \; \big( (x = y) \wedge (z \in x) \big)$
\eenum

where $x$, $y$, and $z$ are again set variables, are not meaningful. In~(\ref{ex 1}), a set variable is used in a place where a formula is expected, in~(\ref{ex 2}) two quantifiers immediately follow each other, in~(\ref{ex 3}) the disjunction requires a formula on the right-hand side, and in~(\ref{ex 4}) a variable is expected on the right-hand side of~$\exists$. If an expression is meaningful, then it is called a formula.

A variable by which we denote a formula is called a formula variable\index{Formula variable}\index{Variable!formula}. If it is evident from the context that a letter is not a formula variable or a defined symbol of the theory, then it is understood that the letter denotes a set variable. For instance, we state the Existence Axiom in Section~\ref{axioms of set theory}, $\exists x \; (x = x)$, and do not explicitly say that $x$ is a set variable.

Given a formula $\varphi$, a set variable that occurs in~$\varphi$ is called free\index{Variable!free} in~$\varphi$ if it does not occur directly after a quantifier. We use the convention that if we list set variables in brackets and separated by commas after a formula variable, then the formula variable denotes a formula that contains as free set variables only those listed in the brackets. For example, $\varphi(x,p)$ stands for a formula that has at most $x$ and $p$ as its free set variables.

We introduce one more symbol,~$\notin$. By $x \notin y$ we mean $\neg (x \in y)$, and similarly for any set variables different from $x$ and~$y$.

Furthermore we introduce some variations of our formal notation, which is often very convenient. First, in a formula we may deliberately omit pairs of parentheses whenever the way how to reinsert the parentheses is obvious. Second, we sometimes use the following simplified notation after quantifiers. If $x$ and $X$ are set variables and $\varphi$ is a formula variable, we write $\exists x \!\in\! X \; \varphi$ instead of $\exists x \; (x \!\in\! X \wedge \varphi)$. Similarly, we write $\forall x \!\in\! X \; \varphi$ instead of $\forall x \; ( x \in X \, \Longrightarrow \, \varphi )$. The same conventions apply, of course, to any other choice of set and formula variables. Third, given sets~$x$, $y$, and $X$, the formula $x \in X \wedge y \in X$ is also written as $x, y \in X$, and similarly for more than two set variables. 

In the nonformal language we adopt the convention that instead of saying that a statement holds "for every $x \in X$" we write $(x \in X)$ after the statement; for example we may write "$x \in Y$ ($x \in X$)" instead of "$x \in Y$ for every $x \in X$". Moreover we use the acronym "iff" which means "if and only if" and thus corresponds to the symbol $\Longleftrightarrow$ in the formal language.

Finally, we remark that the usage of formula variables in this text is restricted to a limited number of occasions. First, formula variables are used in the postulation of two Axiom schemas, the Separation schema, Axiom~\ref{axio sep schema}, and the Replacement schema, Axiom~\ref{axio repl schema}, and its immediate consequences Lemma and Definition~\ref{defi set brackets}, Definitions~\ref{defi set brackets short} and~\ref{defi power set formula}, and Lemmas~\ref{lemm sep schema two par} and~\ref{lemm replacement set}. Whenever any of these is used later in the text, the formula variable is substituted by an actual formula. In particular, no further derivation is undertaken where formula variables are used without previous substitution of specific formulae. Second, formula variables are used in Definitions~\ref{defi ev fre sequence} and~\ref{defi ev fre nets} of the notions "eventually" and "frequently". However, whenever these notions are used later, the formula variable of the definition is substituted by a specific formula.

\section{Axioms of set theory}
\label{axioms of set theory}

The axioms of set theory that we postulate in this Section and use throughout the text are widely accepted in the literature~\cite{Bernays,Ebbinghaus,Jech,Suppes}. They are called ZFC\index{ZFC}\index{Zermelo Fraenkel with choice axiom} ("Zermelo Fraenkel with choice axiom"). There are other axioms that have similar implications for mathematical theories, e.g.\ NBG\index{NBG} ("von Neumann Bernays G\"{o}del")\index{von Neumann Bernays G\"{o}del}, see~\cite{Bernays}. Although we discuss certain aspects of~ZFC in this work, the comparison with NBG or other axioms is beyond our scope.

First we postulate an axiom that says that the world of abstract mathematical objects, which are sets and only sets in our theory, contains at least some object.

\midvspace

\baxio[Existence]
\label{axio exist}
\index{Existence axiom}
\index{Axiom!existence}
\[
\exists x \; x = x
\]

\eaxio

Logically, the formula $x = x$ is always true. Thus Axiom~\ref{axio exist} postulates the existence of a set. The existence of sets does not follow from the other axioms presented subsequently. This is briefly discussed at the end of this Section.

Next we specify a condition under which two sets are equal.

\midvspace

\baxio[Extensionality]
\label{axio extension}
\index{Extensionality axiom}
\index{Axiom!extensionality}
\[
\forall x \; \forall y \; \big( \forall z \; (z \in x \Longleftrightarrow z \in y) \Longrightarrow x = y \big)
\]
\eaxio

The interpretation of Axiom~\ref{axio extension} is that two sets are equal if they have the same elements. The converse implication
\[
\forall x \; \forall y \; \big( x = y \Longrightarrow
\forall z \; (z \in x \Longleftrightarrow z \in y) \big)
\]

is logically true in any theory because of the interpretation of the symbol~$=$.

\midvspace

\bdefi
\label{defi subset}
\index{Subset}
\index{Set!subset}
Given two sets $X$ and $Y$, we say that $Y$ {\bf is a subset of}~$X$, written $Y \subset X$ or $X \supset Y$, if the following statement holds:
\[
\forall y\; y \in Y \Longrightarrow y \in X
\]

We also write $Y \not\subset X$ for $\neg \left( Y \subset X \right)$.
\edefi

Notice that Definition~\ref{defi subset} introduces two new symbols $\subset$ and $\supset$ in the formal language, and also specifies a new notion in the nonformal language. Clearly, in the formal language the new symbol is in principle redundant, that is the same expressions can be written down without it. Thus the new symbols are merely abbreviations. Similarly as for~$\in$, we also adopt the short notation $Y, Z \subset X$ for $Y \subset X \wedge Z \subset X$.

Next we postulate the axioms that allow us to specify a set in terms of a given property which is formalized by a formula.

\midvspace

\baxio[Separation schema]
\label{axio sep schema}
\index{Separation schema}
\index{Axiom!separation schema}
Let $\varphi(x,p)$ be a formula. We have:
\[
\forall p \; \forall X \; \exists Y \; \forall x \; \big(x \in Y \Longleftrightarrow x \in X \wedge \varphi(x,p)\big)
\]

\eaxio

In Axiom~\ref{axio sep schema}, for every formula that contains at most $x$ and $p$ as free variables, one axiom is postulated. Therefore not only a single mathematical expression is postulated here, but a method is given how to write down an axiom for every given formula~$\varphi(x,p)$. This is called a schema. The general analysis of this concept is beyond our scope. However, it is clear that if we would like to specify sets in terms of certain properties one would either write down an axiom for each desired property and restrict oneself to a limited number of properties or define a generic method to specify the axioms. In ZFC the latter possibility is chosen. However, note that only in Lemma and Definition~\ref{defi set brackets}, Definitions~\ref{defi set brackets short} and~\ref{defi power set formula}, and Lemma~\ref{lemm sep schema two par} the schema is used with its formula variable. In all other instances when we refer in this text to the Axiom schema or one of the mentioned definitions or results, we substitute an explicit formula for the formula variable. Since in this text this happens only a finite number of times, we could postulate a finite number of axioms instead of the Separation schema, each with an explicit formula substituted. In this sense, Axiom~\ref{axio sep schema} is only an abbreviated notation of a list of a finite number of axioms that do not contain formula variables. Notice however that the restriction to a finite number of axioms generally also constrains the implications that can be concluded {\it from} the statements proven in this text.

The Separation schema has an important consequence, viz.\ there exists a set that has no element.

\midvspace

\blede
\label{lede empty set}
\index{Empty set}
\index{Set!empty}
There is a unique set~$Y$ such that there exists no set~$x$ with $x \in Y$. $Y$~is called the {\bf empty set}, written~$\O$.
\elede

\bproof
We may choose a set~$X$ by the Existence axiom. Let $\varphi(x)$ denote the formula $\neg (x = x)$. This formula is logically false in any theory for every~$x$. There exists a set~$Y$ such that
\[
\forall x \; x \in Y \Longleftrightarrow x \in X \wedge \neg (x = x)
\]

by the Separation schema. Clearly, $Y$ has no element. The uniqueness of~$Y$ follows by the Extensionality axiom.
\eproof

We have postulated the existence of a set by the Existence axiom and concluded in Lemma and Definition~\ref{lede empty set} that the empty set exists. However, we have not proven so far that any other set exists. This is remedied by the Power set axiom to be introduced below in this Section, and, even without Power set axiom, by the Infinity axiom below.

We now introduce several notations that are all well-defined by the Axioms postulated so far, namely the curly bracket notation for sets, the intersection of two sets, the intersection of one set, and the difference of two sets.

\midvspace

\blede
\label{defi set brackets}
\index{Set!brackets}
Let $X$ and $p$ be sets, and $\varphi(x,p)$ a formula. There is a unique set $Y$ such that
\[
\forall x\; x \in Y \Longleftrightarrow x \in X \wedge \varphi(x,p)
\]

We denote $Y$ by $\left\{ x \in X \, : \, \varphi(x,p) \right\}$.
\elede

\bproof
The existence follows by the Separation schema. The uniqueness is a consequence of the Extensionality axiom.
\eproof

In the particular case where the formula and the parameter are such that they define a set without restricting the members to a given set, we may use the following shorter notation.

\midvspace

\bdefi
\label{defi set brackets short}
\index{Set!brackets}
Let $p$ be a set and $\varphi(x,p)$ a formula. If there is a set $Y$ such that
\[
\forall x\; x \in Y \Longleftrightarrow \varphi(x,p),
\]

then $Y$ is denoted by~$\left\{ x \, : \, \varphi(x,p) \right\}$.
\edefi

\bdefi
\label{defi small intersection}
\index{Intersection}
\index{Set!intersection}
Let $X$ and $Y$ be sets. The set $\left\{ z \in X \, : \,  z \in Y \right\}$ is called {\bf intersection of $X$ and $Y$}, written $X \cap Y$.
\edefi

\brema
Let $X$ and $Y$ be sets. We clearly have
\[
\forall z \; z \in X \cap Y \Longleftrightarrow z \in X \wedge z \in Y
\]

Thus in Definition~\ref{defi small intersection} the sets $X$ and~$Y$ may be interchanged without changing the result for their intersection.
\erema

\bdefi
Let $X$ and $Y$ be two sets. $X$ and $Y$ are called {\bf disjoint} if $X \cap Y = \O$. Given a set $Z$, the members of~$Z$ are called {\bf disjoint} if $x \cap y = \O$ for every $x, y \in Z$.
\edefi

\bdefi
\label{defi great intersection}
\index{Intersection}
\index{Set!intersection}
Let $X$ be a set. If $X \neq \O$, then the set $\left\{ y \, : \, \forall x \in X\; y \in x\right\}$ is called {\bf intersection of}~$X$, written $\bigcap X$.
\edefi

In Definition~\ref{defi great intersection} the short notation introduced in Definition~\ref{defi set brackets short} can be used since we may choose $z \in X$ such that $\bigcap X = \left\{ y \in z \, : \, \forall x \in X\; y \in x\right\}$. As we have not proven so far that any other than the empty set exists, "$X \neq \O$" is stated as a condition in Definition~\ref{defi great intersection}, which may in principle never be satisfied. As already mentioned, the Power set axiom as well as the Infinity axiom each guarantee (without the other one) the existence of a large number of sets.

\midvspace

\brema
\label{rema gen intersection}
The intersection of a set as defined in Definition~\ref{defi great intersection} is sometimes generalized in the following way in the literature (see e.g.~\cite{Jech}):
\\

\hspace{0.05\textwidth}
\parbox{0.95\textwidth}
{Let $p$ be a set and $\varphi(X,p)$ a formula. If there exists a set $X$ such that $\varphi(X,p)$ is true, then the set $Y = \left\{ x \, : \, \forall X \; \varphi(X,p) \Longrightarrow x \in X \right\}$ is well-defined. If, in addition, there is a set $Z$ such that
\[
\forall z \; z \in Z \Longleftrightarrow \varphi(z,p),
\]

then we have $Y = \bigcap Z$.}
\\

The last claim shows that Definition~\ref{defi great intersection} is a special case of the first claim. 

Note that this generalization involves a formula variable, which we prefer to avoid. In this text the generalization of the intersection of a set is only used once, viz.\ in the definition of the natural numbers (Definition~\ref{defi naturalnumbers}). Their existence is a consequence of the Separation schema with a concrete formula.
\erema

The following result states that there exists no set that contains all sets.

\midvspace

\blemm
We have
\[
\forall X\; \exists x\; x \notin X
\]
\elemm

\bproof
Let $X$ be a set. Then we have $\left\{ x \in X \, : \, x \notin x \right\} \notin X$. 
\eproof

We now postulate that for two given sets $X$ and $Y$ there is a set that contains all elements of $X$ and~$Y$.

\midvspace

\baxio[Small union]
\label{axio small union}
\index{Small union axiom}
\index{Axiom!small union}
\[
\forall X \; \forall Y \; \exists Z \; \forall z \; (z \in X \vee z \in Y \Longrightarrow z \in Z)
\]

\eaxio

\bdefi
\index{Difference}
\index{Set!difference}
\index{Complement}
\index{Set!complement}
Given two sets $X$ and $Y$, the set $\left\{ x \in X \, : \, x \notin Y \right\}$ is called {\bf difference of} $X$ and~$Y$, written $X \setminus Y$. If $Y \subset X$, the set $X \setminus Y$ is also called {\bf complement of}~$Y$ whenever the set~$X$ is evident from the context. The complement of $Y$ is also denoted by~$Y^c$.
\edefi

\blede
\label{lede small union}
\index{Union}
\index{Set!union}
Let $X$ and $Y$ be two sets. Furthermore, let $Z$ be a set such that $X, Y \subset Z$. The set $\big(X^c \cap Y^c\big)^c$, where the complement is with respect to~$Z$, is called
{\bf union of $X$ and $Y$}, written $X \cup Y$.
\elede

\bproof
The existence of~$Z$ follows by the Small union axiom. To see that the definition of $X \cup Y$ is independent of the choice of~$Z$, let $W$ be another set with $X, Y \subset W$. We clearly have
\[
Z \setminus\! \big((Z \!\setminus\! X) \cap (Z \!\setminus\! Y)\big) = W \setminus\! \big((W \!\setminus\! X) \cap (W \!\setminus\! Y)\big)
\]

\eproof

\brema
Let $X$ and $Y$ be sets. We have
\[
\forall z \; z \in X \cup Y \Longleftrightarrow \varphi(z,X,Y)
\]

where $\varphi(z,X,Y)$ denotes the formula $(z \in X) \vee (z \in Y)$. Since this formula contains $z$ and two parameters as free variables, we cannot use the Separation schema in the above form (i.e.\ Axiom~\ref{axio sep schema}) to define the union of~$X$ and~$Y$.
\erema

In the following two Lemmas we list several important equalities that hold for the unions, intersections, and differences of two or three sets, and for the complement of subsets of a given set.

\midvspace

\blemm
\label{lemm set equalities}
Given three sets $A$, $B$, and $C$, the following equalities hold:

\benum
\item \label{lemm set equalities 1} $A \cup B = B \cup A$,\;\;\; $A \cap B = B \cap A$
\item \label{lemm set equalities 2} $A \cup (B \cup C) = (A \cup B) \cup C$,\;\;\; $A \cap (B \cap C) = (A \cap B) \cap C$
\item \label{lemm set equalities 3} $A \cap (B \cup C) = (A \cap B) \cup (A \cap C)$,\;\;\; $A \cup (B \cap C) = (A \cup B) \cap (A \cup C)$
\item \label{lemm set equalities 4} $A \setminus (A \setminus B) = A \cap B$
\eenum

\elemm

\bproof
Exercise.
\eproof

\blemm[De Morgan]
\label{lemm demorgan formulas}
\index{De Morgan equalities}
Given a set $X$ and $A, B \subset X$, the following equalities hold:
\[
(A \cap B)^c = A^c \cup B^c,\quad (A \cup B)^c = A^c \cap B^c
\]

where the complement is with respect to~$X$ in each case.
\elemm

\bproof
Exercise.
\eproof

\baxio[Great union]
\label{axio great union}
\index{Great union axiom}
\index{Axiom!great union}
\[
\forall X \;\exists Y \; \forall y \; \big( \exists x \; (x \in X \wedge y \in x) \Longrightarrow y \in Y \big)
\]

\eaxio

\blede
\label{lede great union}
\index{Union}
\index{Set!union}
Let $X$ be a set. If $X \neq \O$, then the set $\left\{ y \, : \, \exists x \; (x \in X \wedge y \in x) \right\}$ is called {\bf union of}~$X$, written~$\bigcup X$.
\elede

\bproof
This set is well-defined by the Great union axiom, the Separation schema, and the Extensionality axiom.
\eproof

In some contexts, in particular when considering unions and intersections of sets, a set is called a system, or a system of sets. Similarly a subset of a set is sometimes called a subsystem. This somewhat arbitrary change in nomenclature is motivated by the fact that in such contexts there intuitively seems to be a hierarchy of, first, the system, second, the members of the system, and, third, the elements of the members of the system. Often this intuitive hierarchy is highlighted by three different types of letters used for the variables, namely script letters for the system (e.g.~$\mc A$), capital Latin letters for the members of the system (e.g.~$A$), and small Latin letters for their elements (e.g~$a$). We emphasize, however, that all objects denoted by these variables are nothing but sets, that is also systems are sets, and that two sets are distinct if an only if they have different elements.

We now postulate that the system of all subsets of a given set is contained in a set.

\midvspace

\baxio[Power set]
\label{axio power set}
\index{Power set axiom}
\index{Axiom!power set}
\[
\forall X \;\exists Y \; \forall y \; (y \subset X \Longrightarrow y \in Y)
\]

\eaxio

\blede
\index{Power set}
\index{Set!power}
Given a set $X$, the set $\left\{ y \, : \, y \subset X \right\}$ is called {\bf power set of}~$X$, written~${\mc P}(X)$. We also write ${\mc P}^2 (X)$ for ${\mc P} \left( {\mc P} (X) \right)$.
\elede

\bproof
This set is well-defined by the Power set axiom, the Separation schema, and the Extensionality axiom.
\eproof

The existence of the power set allows the following variation of Definition~\ref{defi set brackets}.

\midvspace

\bdefi
\label{defi power set formula}
Let $X$ and $p$ be sets and $\varphi(x,p)$ a formula. The set $\left\{ x \in {\mc P}(X) \, : \, \varphi(x,p)\right\}$ is also denoted by $\left\{ x \subset X \, : \, \varphi(x,p) \right\}$.
\edefi

As a consequence of the Power set axiom, for every set $X$ there exists a set that contains $X$ and only $X$ as element. We introduce the following notation and nomenclature.

\midvspace

\bdefi
\label{defi singleton}
For every set $X$, the set $\left\{ Y \subset X \, : \, Y = X \right\}$ is denoted by~$\left\{ X \right\}$.
\edefi

\bdefi
\index{Singleton}
\index{Set!singleton}
Let $X$ be a set. $X$ is called a {\bf singleton} if there is a set $x$ such that $X = \left\{ x \right\}$.
\edefi

With the Axioms postulated so far and without the Power set axiom we do not know about the existence of any other set than the empty set. Including the Power set axiom we conclude that also $\left\{ \O \right\}$ and $\left\{ \left\{ \O \right\} \right\}$ are sets. Clearly these three sets are distinct by the Extensionality axiom.

The following Lemma and Definition states that for every two sets $X$ and~$Y$ there is a set whose members are precisely $X$ and~$Y$.

\midvspace

\blede
\label{lede pair set}
Given two sets $X$ and $Y$, there is a set $Z$ such that
\[
\forall z\; z \in Z \Longleftrightarrow z = X \vee z = Y
\]

We also denote $Z$ by $\left\{ X, Y \right\}$.
\elede

\bproof
Notice that $\left\{ X \right\}$, $\left\{ Y \right\}$ are sets by the Power set axiom. Let $Z = \left\{ X \right\} \cup \left\{ Y \right\}$.
\eproof

This shows that e.g.\ also $\left\{ \O, \left\{ \O \right\} \right\}$ is a set.

\midvspace

\brema
Given a set $X$, we have $\left\{ X, X \right\} = \left\{ X \right\}$.
\erema

\brema
Let $X$ and $Y$ be two sets. We have

\benum
\item $X \cap Y = \, \bigcap \left\{ X, Y \right\}$
\item $X \cup Y = \, \bigcup \left\{ X, Y \right\}$
\eenum

Thus Definition~\ref{defi small intersection} and Lemma and Definition~\ref{lede small union} can be considered as special cases of Definition~\ref{defi great intersection} and Lemma and Definition~\ref{lede great union}, respectively.
\erema

For convenience we also introduce a notation for a set with three members.

\midvspace

\blede
\label{lede triple set}
Given sets $X$, $Y$, and $Z$, there is a set $U$ such that
\[
\forall u\; u \in U \Longleftrightarrow \left( u = X \right) \vee \left( u = Y \right) \vee \left( u = Z \right)
\]

We also denote $U$ by $\left\{ X, Y, Z \right\}$.
\elede

\bproof
We may define $U = \left\{ X, Y \right\} \cup \left\{ Z \right\}$.
\eproof

It is obvious that the order in which we write the sets $X$ and $Y$ in Lemma and Definition~\ref{lede pair set}, or the order in which we write $X$, $Y$, and~$Z$ in Lemma and Definition~\ref{lede triple set} does not play a role. In many contexts we need a system that identifies two (distinct or equal) sets and also specifies their order. This is achieved by the following concept.

\midvspace

\blede
\index{Ordered pair}
\index{Coordinates}
Let $X$ be a set. $X$ is called {\bf ordered pair}, or short {\bf pair}, if there are sets $x$ and $y$ such that $X = \left\{ \left\{ x, y \right\}, \left\{ x \right\} \right\}$. In this case we have $\bigcup \bigcap X = x$ and $\bigcup \big( \bigcup X \!\setminus\! \left\{ x \right\} \!\big) = y$. Moreover, in this case $x$ and $y$ are called {\bf left} and {\bf right coordinates} of~$X$, respectively. Further, if $X$ is an ordered pair, its unique left and right coordinates are denoted by $X_l$ and $X_r$, respectively, and $X$ is denoted by~$(X_l,X_r)$.
\elede

\bproof
Exercise.
\eproof

\brema
Let $x$ and $y$ be two sets. Then $x \neq y$ implies $(x,y) \neq (y,x)$.
\erema

The concept of ordered pair allows the extension of the Separation schema, Axiom~\ref{axio sep schema}, to more than one parameter. The following is the result for two parameters.

\midvspace

\blemm[Separation schema with two parameters]
\label{lemm sep schema two par}
\index{Separation schema}
Let $\varphi(x,p,q)$ be a formula. We have:
\[
\forall p \; \forall q \; \forall X \; \exists Y \; \forall x \; \big(x \in Y \Longleftrightarrow x \in X \wedge \varphi(x,p,q)\big)
\]

\elemm

\bproof
Let $p$, $q$, and $X$ be sets. We define $r = (p,q)$ and $Y = \left\{ x \in X \, : \, \psi(x,r) \right\}$ where $\psi(x,r)$ is the formula $\exists u \; \exists v \; r = (u,v) \wedge \varphi(x,u,v)$. Then $Y$ is the required set.
\eproof

The Separation schema is used in the following two Lemmas. Before stating these Lemmas we would like to relax the rules for the usage of the set brackets $\left\{ \ldots \right\}$ that are defined in Definition~\ref{defi set brackets}. Remember that such a modified notation is already defined in the case where the formula specifies a set (cf.~Definition~\ref{defi set brackets short}) and in the case of a subset of the power set (cf.~Definition~\ref{defi power set formula}). We now agree that we may use a comma instead of $\wedge$ between two or more formulae on the right hand side of the colon; e.g.\ given two sets $X$ and $Y$, we may write $\left\{ x \, : \, x \in X,\; x \in Y \right\}$ instead of $\left\{ x \, : \, x \in X \wedge x \in Y \right\}$. Moreover we agree that we may use all defined symbols on the left hand side of the colon, thereby eliminating one or more $\exists$ and one equality on the right hand side; e.g.\ given two sets $X$ and $Y$, we may write $\left\{ x \cap y \, : \, x \in X,\; x \in Y \right\}$ instead of $\left\{ z \, : \, \exists x \in X \; \exists x \in Y \; z = x \cap y \right\}$, and $\left\{ Y \!\setminus x \, : \, x \in X \right\}$ instead of $\left\{ z \, : \, \exists x \in X \; z = Y \!\setminus x \right\}$. This is precisely the same kind of notation as in Definition~\ref{defi power set formula}.

\midvspace

\blemm
\label{lemm set equalities great}
Let $\mc A$ and $\mc B$ be two systems of sets. If ${\mc A} \neq \O$ and ${\mc B} \neq \O$, the following equalities hold:
\benum
\item \label{lemm set equalities great 1} $\big( \bigcup {\mc A} \, \big) \cap \big( \bigcup {\mc B} \, \big) \, = \, \bigcup \left\{ A \cap B \, : \, A \in {\mc A},\, B \in {\mc B} \right\}$
\item \label{lemm set equalities great 2} $\big( \bigcap {\mc A} \, \big) \cup \big( \bigcap {\mc B} \, \big) \, = \, \bigcap \left\{ A \cup B \, : \, A \in {\mc A},\, B \in {\mc B} \right\}$
\eenum

In particular, the right hand sides of~(\ref{lemm set equalities great 1}) and~(\ref{lemm set equalities great 2}) are well-defined.
\elemm

\bproof
Exercise.
\eproof

\blemm[De Morgan]
\label{lemm demorgan formulas great}
\index{De Morgan equalities}
Given a set $X$ and a system ${\mc A} \subset {\mc P}(X)$ with ${\mc A} \neq \O$, the following equalities hold:
\benum
\item \label{lemm demorgan formulas great 1} $\big( \bigcap {\mc A} \, \big)^c = \, \bigcup \left\{ A^c \, : \, A \in {\mc A} \right\}$
\item \label{lemm demorgan formulas great 2} $\big( \bigcup {\mc A} \, \big)^c = \, \bigcap \left\{ A^c \, : \, A \in {\mc A} \right\}$
\eenum

where the complements refer to the set~$X$. In particular, the right hand sides of~(\ref{lemm demorgan formulas great 1}) and~(\ref{lemm demorgan formulas great 2}) are well-defined.
\elemm

\bproof
Exercise.
\eproof

\brema
Notice that the equalities (\ref{lemm set equalities 3}) in Lemma~\ref{lemm set equalities} are special cases of those in Lemma~\ref{lemm set equalities great}, and that the equalities in Lemma~\ref{lemm demorgan formulas} are special cases of those in Lemma~\ref{lemm demorgan formulas great}.
\erema

We now introduce the set of all ordered pairs such that the left coordinate is in~$X$ and the right coordinate is in~$Y$ where $X$ and $Y$ are two sets. The following Definition uses again the Separation schema with two parameters.

\midvspace

\bdefi
\label{defi cartesian prod pair}
\index{Cartesian product}
Let $X$ and $Y$ be two sets. The set
\[
\left\{ z \in {\mc P}^2 (X \cup Y) \, : \, \exists x \in X\; \exists y \in Y\; z = \left\{ \left\{ x, y \right\}, \left\{ x \right\} \right\} \right\}
\]

is called {\bf Cartesian product of} $X$ and~$Y$, and denoted by~$X \!\times Y$.
\edefi

We now state some properties of the Cartesian product of two sets.

\midvspace

\blemm
Let $U$, $V$, $X$, and $Y$ be sets. Then we have

\benum
\item $\O \times V = U \times \O = \O$
\item $(U \!\times V) \cap (X \!\times Y) = (U \cap X) \times (V \cap Y)$
\item $X \!\times\! (V \cup Y) = (X \!\times V) \cup (X \!\times Y)$
\item $X \!\times\! (Y \!\setminus V) = (X \!\times Y) \setminus (X \!\times V)$
\eenum

\elemm

\bproof
Exercise.
\eproof

\blemm
Let $U$, $V$, $X$, and $Y$ be sets with $X \subset U$ and $Y \subset V$. Then we have
\[
(X \!\times Y)^c = (X^c \!\times Y^c) \cup (X^c \!\times Y) \cup (X \!\times Y^c)
\]

where the first complement refers to $U \!\times V$, the complement of $X$ refers to~$U$, and the complement of $Y$ refers to~$V$.
\elemm

\bproof
Exercise.
\eproof

Similarly to ordered pairs we introduce the notion of ordered triple consisting of three sets in a specific order.

\midvspace

\bdefi
\index{Ordered triple}
Let $X$ be a set. $X$ is called {\bf ordered triple} if there are sets $x$, $y$, and $z$ such that $X = \big((x,y),z\big)$. In this case $X$ is also denoted by $(x,y,z)$.
\edefi

Notice that the Separation schema may be extended to three parameters by using ordered triples. Similarly, we may clearly define ordered quadruples etc.\ and define corresponding Separation schemas.

\midvspace

\baxio[Infinity]
\label{axio infinity}
\index{Infinity axiom}
\index{Axiom!infinity}
\[
\exists X \; \Big( \O \in X \wedge \forall x \; \big( x \in X \Longrightarrow x \cup \left\{ x \right\} \in X \big) \Big)
\]

\eaxio

\bdefi
\index{Inductive set}
\index{Set!inductive}
A set $X$ is called {\bf inductive} if it has the following properties:
\benum
\item $\O \in X$
\item $\forall x \in X \quad x \cup \left\{ x \right\} \in X$
\eenum

\edefi

Obviously, the Infinity axiom says that there exists an inductive set. 

\midvspace

\bdefi
\label{defi naturalnumbers}
\index{Natural numbers}
\index{Numbers!natural}
Let $X$ be an inductive set. The members of the set
\[
\left\{ n \in X \, : \, \forall Y\, \mbox{$Y$ is inductive} \Longrightarrow n \in Y \right\}
\]

are called {\bf natural numbers}. The set of natural numbers is denoted by~$\naturalnumbers$.
\edefi

Notice that Definition~\ref{defi naturalnumbers} does not depend on the choice of the inductive set~$X$. This is an example of the concept discussed in Remark~\ref{rema gen intersection}, which can be understood as a generalized form of intersection.

\midvspace

\brema
The set~$\naturalnumbers$ is inductive.
\erema

The following Axiom is part of~ZFC, essentially in order to obtain the statements in the following Lemma.

\midvspace

\baxio[Regularity]
\label{axio regularity}
\index{Regularity axiom}
\index{Axiom!regularity}
\[
\forall X \; X \neq \O \Longrightarrow \exists x\; \left( x \in X \right) \wedge \left( X \cap x = \O \right)
\]

\eaxio

\blemm
\label{lemm element min}
The following statements hold:
\benum
\item \label{lemm element min 1} $\neg \exists X\; X \in X$
\item \label{lemm element min 2} $\neg \exists X\; \exists Y\; \left( X \in Y \right) \wedge \left( Y \in X \right)$
\item \label{lemm element min 3} $\neg \exists X\; \exists Y\; \exists Z\; \left( X \in Y \right) \wedge \left( Y \in Z \right) \wedge \left( Z \in X \right)$
\eenum

\elemm

\bproof
To see (\ref{lemm element min 1}), let $X$ be a set. Notice that $\left\{ X \right\}$ is non-empty and thus $\left\{ X \right\} \cap X = \O$ by the Regularity axiom.

To see (\ref{lemm element min 2}), we assume that $X$ and $Y$ are two sets such that $X \in Y$ and $Y \in X$. Since the set $Z = \left\{ X, Y \right\}$ is non-empty there is $z \in Z$ such that $Z \cap z = \O$ by the Regularity axiom. However, $z = X$ implies $Z \cap X = Y$, and $z = Y$ implies $Z \cap Y = X$, which is a contradiction.

Finally, to show (\ref{lemm element min 3}), we assume that there are sets $X$, $Y$, and $Z$ such that $X \in Y$, $Y \in Z$, and $Z \in X$. We define $W = \left\{ X, Y, Z \right\}$. Application of the Regularity axiom again leads to a contradiction.
\eproof

\baxio[Replacement schema]
\label{axio repl schema}
\index{Replacement schema}
\index{Axiom!replacement schema}
Let $\varphi(x,y,p)$ be a formula.

\begin{center}
\begin{tabular}{l}
$\forall p \; \Big(\big(\forall x \; \forall y \; \forall z \; (\varphi(x,y,p) \wedge \varphi(x,z,p) \Longrightarrow y = z)\big)$\\[.5em]
$\quad \quad \quad \Longrightarrow \forall X \; \exists Y \; \forall y \; \big(\exists x \; (x \in X \wedge \varphi(x,y,p)) \Longrightarrow y \in Y \big)\Big)$
\end{tabular}
\end{center}

\eaxio

Notice that Axiom~\ref{axio repl schema} is not a single axiom but a schema of axioms. This concept is discussed above in the context of the Separation schema, Axiom~\ref{axio sep schema}. The Replacement schema is applied below to derive Theorem~\ref{theo ordinal injection} where a concrete formula is substituted.

The premise in Axiom~\ref{axio repl schema} says that for every $x$ there is at most one $y$ such that $\varphi(x,y,p)$ is satisfied. The conclusion states that if the sets $x$ are taken out of a given set~$X$, there exists a set $Y$ that has those sets~$y$ among its members. Clearly one may also define a set $Y$ that has precisely those sets~$y$ as its members (and no others) by the Separation schema. This is the result of the following Lemma.

We remark that, as in the context of the Separation schema (cf.\ Lemma~\ref{lemm sep schema two par}), the Replacement schema can be extended to two or more parameters. However, this is not required in the present work.

\midvspace

\blemm
\label{lemm replacement set}
Let $X$ and $p$ be sets, and $\varphi(x,y,p)$ a formula such that
\[
\forall x \; \forall y \; \forall z \; \varphi(x,y,p) \wedge \varphi(x,z,p) \Longrightarrow y = z
\]

holds. Then the set $\left\{ y \, : \, \exists x \in X \; \varphi(x,y,p) \right\}$ is well-defined.
\elemm

\bproof
The existence follows by the Replacement schema and the Separation schema. The uniqueness is a consequence of the Extensionality axiom.
\eproof

\baxio[Choice]
\label{axio choice}
\index{Choice axiom}
\index{Axiom!choice}
\begin{eqnarray*}
&&\forall X \; \exists Z \; \forall Y \in {\mc P}(X) \!\setminus\! \left\{ \O \right\} \; \Big( \big(\exists y \in Y \;\; (Y,y) \in Z \big)\\
&&\quad \quad \quad \; \wedge \; \big( \forall u, v \in X \;\; (Y,u), (Y,v) \in Z \Longrightarrow u = v \big) \Big)
\end{eqnarray*}

\eaxio

\blede
\label{lede choice function}
\index{Choice function}
Let $X$ be a set and $Y \subset {\mc P}(X)$ with $\O \notin Y$. There exists a set $Z \subset Y \!\times X$ such that the following statements hold:

\benum
\item \label{lede choice function 1} $\forall z \in Z \quad z_r \in z_l$
\item \label{lede choice function 2} $\forall y \in Y \quad \exists z \in Z \quad z_l = y$
\item \label{lede choice function 3} $\forall w, z \in Z \quad \left( w_l = z_l \, \Longrightarrow \, w_r = z_r \right)$
\eenum

$Z$ is called a {\bf choice function}.
\elede

\bproof
Given the stated conditions, there exists, by Axiom~\ref{axio choice}, a set $W$ such that $Z = W \cap \left(Y \!\times X\right)$ satisfies (\ref{lede choice function 1}) to~(\ref{lede choice function 3}).
\eproof

Given the other axioms above, the Choice axiom can be stated in several equivalent forms. Two versions are presented in Section~\ref{choice}, viz.\ the Well-ordering principle and Zorn's Lemma.

We finally introduce a notation that is convenient when dealing with topologies, topological bases and other concepts to be introduced below.

\midvspace

\bdefi
\label{defi set system point}
Let $X$ be a set, ${\mc A} \subset {\mc P}(X)$, and $x \in X$. We define ${\mc A}(x) = \left\{ A \in {\mc A} \, : \, x \in A \right\}$.
\edefi

Notice that the existence of sets, or even of a single set, is not guaranteed if the Existence axiom is not postulated. This is because clearly none of the other axioms postulates the existence of a set---without any other set already existing---apart from the Infinity axiom~\ref{axio infinity} below. There however the definition of the empty set is used, which in turn is defined in Lemma and Definition~\ref{lede empty set} by usage of the Separation schema. The Separation schema always refers to an existing set. It is possible to modify the Infinity axiom such that it also postulates the existence of a set, see e.g.~\cite{Jech}.


\chapter{Relations}
\label{relations}
\setcounter{equation}{0}

\pagebreak

\section{Relations and orderings}

In this Section we introduce the concept of relation, which is fundamental in the remainder of the text. Many important special cases are analysed, in particular orderings. Also functions, that are introduced in the next Section, are relations.

\midvspace

\bdefi
\label{defi rel inv prod}
\index{Relation}
\index{Inverse}
\index{Relation!inverse}
\index{Product}
\index{Relation!product}
\index{Relational space}
\index{Space!relational}
\index{Diagonal}
Given two sets $X$ and $Y$, a subset $U \subset X \!\times Y$ is called a {\bf relation on} $X \!\times Y$. The {\bf inverse of}~$U$ is a relation on $Y \!\times X$ and defined as
\[
U^{-1} = \big\{ (y,x) \in Y \!\times X \, : \, (x,y) \in U \big\}
\]

Given another set $Z$ and a relation $V \subset Y \!\times Z$, the {\bf product of~$V$ and~$U$} is defined as
\[
V U = \big\{ (x,z) \in X \!\times Z \, : \, \exists y \in Y \; (x,y) \in U,\, (y,z) \in V \big\}
\]

A relation $R$ on $X \!\times X$ is also called {\bf relation on}~$X$. In this case the pair $(X,R)$ is called {\bf relational space}. Furthermore, the set $\Delta = \left\{ (x,x) \, : \, x \in X \right\}$ is called {\bf diagonal}.
\edefi

Notice the order of $U$ and $V$ in the definition of the product, which may be counterintuitive.

\midvspace

\bdefi
\label{def relation evaluation}
\index{Range}
\index{Relation!range}
\index{Full range}
\index{Range!full}
\index{Domain}
\index{Relation!domain}
\index{Full domain}
\index{Domain!full}
\index{Field}
\index{Relation!field}
\index{Full field}
\index{Field!full}
Given two sets $X$ and $Y$, a relation $R \subset X \!\times Y$, and a set $A \subset X$ we introduce the following notation:
\benum
\item $R \left[ A \right] = \big\{ y \in X \, : \, \exists x \in A \; (x,y) \in R \big\}$
\item $R \left\{ x \right\} = R \left[ \left\{ x \right\} \right]$,\; that is \,\! $R \left\{ x \right\} = \left\{ y \in X \, : \, (x,y) \in R \right\}$
\item $R \langle A \rangle = \big\{ y \in X \, : \, \forall x \in A \; (x,y) \in R \big\}$
\eenum

The set $R \left[ X \right]$ is called the {\bf range of}~$R$, written~$\mathrm{ran}(R)$ or~$\mathrm{ran} \, R$. The set~$R^{-1} \left[ Y \right]$ is called the {\bf domain of}~$R$, written~$\mathrm{dom}(R)$ or~$\mathrm{dom} \, R$. We say that $R$ {\bf has full range} if $\mathrm{ran} \, R = Y$, and {\bf full domain} if $\mathrm{dom} \, R = X$.

Given a relation $S$ on~$X$, the set $\left(S \cup S^{-1}\right) \left[ X \right] = (\mathrm{dom} \, S) \cup (\mathrm{ran} \, S)$ is called the {\bf field of}~$S$, written $\mathrm{field}(S)$ or~$\mathrm{field} \, S$. We say that $S$ {\bf has full field} if \,\! $\mathrm{field} \, S = X$.
\edefi

Clearly, a relation $S$ on $X$ that has full domain or full range also has full field.

Notice that Definition~\ref{defi rel inv prod} of a relation $R$ on a Cartesian product $X \!\times Y$ specifies the sets $X$ and $Y$ from which the product is formed although $X$ may not be the domain and $Y$ may not be the range of~$R$. This is important in the context of functional relations to be defined in Section~\ref{functions}.

\midvspace

\bdefi
Let $\mc R$ be a system of relations on $X \!\times Y$. We define
\benum
\item ${\mc R} \left[ A \right] = \left\{ R \left[ A \right] \, : \, R \in {\mc R} \right\}$
\item ${\mc R} \left\{ x \right\} = {\mc R} \left[ \left\{ x \right\} \right]$,\; that is \,\! ${\mc R} \left\{ x \right\} = \left\{ R \left\{ x \right\} \, : \, R \in {\mc R} \right\}$
\eenum

The sets
\[
{\textstyle \bigcup} \left( {\mc R} \left[ X \right] \right) \, = \, {\textstyle \bigcup} \left\{ R \left[ X \right] \, : \, R \in {\mc R} \right\},
\quad \quad
{\textstyle \bigcup} \left\{ R^{-1} \left[ Y \right] \, : \, R \in {\mc R} \right\}
\]

are called {\bf range} and {\bf domain of}~$\mc R$, respectively. 

If $X = Y$, then the set
\[
{\textstyle \bigcup} \left\{ R \left[ X \right] \, : \, R, R^{-1} \in {\mc R} \right\}
\]

is called {\bf field of}~$\mc R$. In this case $\mc R$~is said to {\bf have full field} if its field is~$X$.
\edefi

\brema
Given two sets $X$ and $Y$, $A \subset X$, and a relation $R \subset X \!\times Y$, the following statements hold:
\benum
\item $R \left\{ x \right\} = R \langle \left\{ x \right\} \rangle$ \,\! for every $x \in X$
\item $R \left[ \O \right] = \O$
\item $R \langle \O \rangle = Y$
\item $R \left[ A \right] = \bigcup \left\{ R \left\{ x \right\} \, : \, x \in A \right\}$
\item $R \langle A \rangle = \bigcap \left\{ R \left\{ x \right\} \, : \, x \in A \right\}$
\eenum

\erema

\bdefi
Let $X$ and $Y$ be sets, $R \subset X \!\times Y$ a relation, and ${\mc A} \subset {\mc P}(X)$. We define
\[
R \, \llbracket {\mc A} \, \rrbracket = \left\{ R \left[ A \right] \, : \, A \in {\mc A} \right\}
\]

\edefi

We now list a few consequences of the above definitions.

\midvspace

\blemm
\label{lemm prop rel}
Given sets $X, Y$ and $Z$, and relations $U, U' \subset X \!\times Y$, and $V, V' \subset Y \!\times Z$ where $U' \subset U$ and $V' \subset V$, and $W \subset Z \times S$ the following statements hold:
\benum
\item $(VU)^{-1} = U^{-1}V^{-1}$
\item \label{lemma relation assoc} $(WV)U = W(VU)$
\item $U'^{-1} \subset U^{-1}$
\item $V'U \subset VU$
\item $VU' \subset VU$
\eenum

\elemm

\bproof
Exercise.
\eproof

By Lemma~\ref{lemm prop rel}~(\ref{lemma relation assoc}) we may drop the brackets in the case of multiple products of relations without generating ambiguities.

\midvspace

\blemm
Given sets $X$, $Y$, and $Z$, relations $U \subset X \!\times Y$ and $V \subset Y \!\times Z$, and a set $A \subset X$, we have $(VU) \left[ A \right] = V \left[ U \left[ A \right] \right]$.
\elemm

\bproof
Exercise.
\eproof

\bdefi
\label{defi relation restriction}
\index{Restriction}
\index{Relation!restriction}
Let $(X,R)$ be a relational space. Then the relation $R \, | \, A = R  \cap \left( A \!\times\! A \right)$ on~$A$ is called the {\bf restriction of} $R$ {\bf to}~$A$.
\edefi

The following properties are important to characterize different types of relational spaces.

\midvspace

\bdefi
\label{properties relation}
\index{Reflexive}
\index{Relation!reflexive}
\index{Antireflexive}
\index{Relation!antireflexive}
\index{Symmetric}
\index{Relation!symmetric}
\index{Antisymmetric}
\index{Relation!antisymmetric}
\index{Transitive}
\index{Relation!transitive}
\index{Connective!relation}
\index{Relation!connective}
\index{Directed}
\index{Relation!directed}
Let $(X,R)$ be a relational space. Then $R$ is called
\benum
\item \label{properties relation reflexive} {\bf reflexive} \, if \, $\Delta \subset R$
\item \label{properties relation antireflexive} {\bf antireflexive} \, if \, $\Delta \cap R = \O$
\item \label{properties relation symmetric}{\bf symmetric} \, if \, $R^{-1} = R$
\item \label{properties relation antisymmetric}{\bf antisymmetric} \, if \, $R \cap R^{-1} \subset \Delta$
\item \label{properties relation transitive} {\bf transitive} \, if \, $R^2 \subset R$
\item \label{properties relation connective} {\bf connective} \, if \, $R \cup R^{-1} \cup \Delta = X \!\times\! X$
\item \label{properties relation directive} {\bf directive} \, if \, $X \!\times\! X = R^{-1} R$
\eenum

\edefi

We remark that antisymmetry of a relation is defined in different ways in the literature, see for example \cite{Gaal}, p.~6, where the definition is $R \cap R^{-1} = \O$, or \cite{Querenburg}, p.~4, where the definition is $R \cap R^{-1} = \Delta$. The expressions "connective" and "directive" are not standard terms in the literature, see however \cite{Ebbinghaus},~p.~58.

\midvspace

\brema
\label{rema rel space prop}
Let $(X,R)$ be a relational space and $A \subset X$. The following statements hold:
\benum
\item $R$ is connective iff for every $x,y \in X$ we have $(x,y) \in R$ or $(y,x) \in R$ or $x = y$.
\item $R$ is directive iff for every $x,y \in X$ there is $z \in X$ such that $(x,z), (y,z) \in R$.
\item $R^{-1}$ is reflexive, antireflexive, symmetric, antisymmetric, transitive, or connective if $R$ has the respective property.
\item $R \, | \, A$ is reflexive, antireflexive, symmetric, antisymmetric, transitive, or connective if $R$ has the respective property.
\eenum

\erema

\blemm
\label{lemma relations union}
Given a set $X$ and relations $U, V$ on $X$ where $V$ is symmetric, we have
\[
VUV = \, \bigcup \big\{ \left(V \! \left\{ x \right\}\right) \times \left(V \! \left\{ y \right\}\right) \, : \, (x,y) \in U \big\}
\]

\elemm

\bproof
If $(u,v) \in VUV$, then there exists $(x,y) \in U$ such that $(u,x), (y,v) \in V$. Therefore we have $(u,v) \in \left(V \! \left\{ x \right\}\right) \times \left(V \! \left\{ y \right\}\right)$. The converse is shown in a similar way.
\eproof

\bdefi
Let $(X,R)$ be a related space and $A \subset X$. $A$~is called a {\bf chain} if $R \, | \, A$ is connective. For definiteness, we also define $\O$ to be a chain.
\edefi

\bdefi
\index{Partition}
Let $X$ be a set and ${\mc A} \subset {\mc P}(X)$. $\mc A$ is called a {\bf partition of}~$X$ if $\bigcup {\mc A} = X$ and $A \cap B = \O$ for every $A, B \in {\mc A}$.
\edefi

\bdefi
\label{defi special relations}
\index{Equivalence relation}
\index{Relation!equivalence}
\index{Equivalence class}
Let $(X,R)$ be a relational space. $R$~is called {\bf equivalence relation} if it is reflexive, symmetric, and transitive. Given a point $x \in X$, the set~$R \left\{ x \right\}$ is called {\bf equivalence class of}~$x$, written~$\left[ x \right]$.
\edefi

\brema
\label{rema equi rel}
Let $X$ be a set, $R$ an equivalence relation on~$X$, and $x, y \in X$. The following equivalences hold:
\[
(x,y) \in R \quad \Longleftrightarrow \quad x \in \left[ y \right] \quad \Longleftrightarrow \quad y \in \left[ x \right] \quad \Longleftrightarrow \quad \left[ x \right] = \left[ y \right]
\]

\erema

\blemm
Given a set $X$ and an equivalence relation~$R$ on~$X$, the system of all equivalence classes is a partition of~$X$, denoted by~$X / R$.
\elemm

\bproof
We clearly have $\bigcup X / R = X$. Now we assume that $u, x, y \in X$ such that $u \in \left[ x \right] \cap \left[ y \right]$. It follows that $\left[ x \right] = \left[ u \right] = \left[ y \right]$ by Remark~\ref{rema equi rel}.
\eproof

\bdefi
\label{defi ordering}
\index{Pre-ordering}
\index{Pre-ordered space}
\index{Space!pre-ordered}
\index{Ordering}
\index{Ordered space}
\index{Space!ordered}
\index{Ordering in the sense of~$<$}
\index{Space ordered in the sense of~$<$}
\index{Ordering in the sense of~$\leq$}
\index{Space ordered in the sense of~$\leq$}
\index{Direction}
\index{Directed space}
\index{Space!directed}
Let $(X,R)$ be a relational space.
\benum
\item \label{defi pre-ordering} $R$ is called a {\bf pre-ordering on}~$X$ if it is transitive. In this case we also write $\prec$ for~$R$, and $x \prec y$ for $(x,y) \in R$. The pair $(X,\prec)$ is called {\bf pre-ordered space}.
\item \label{defi ordering general} $R$ is called an {\bf ordering on}~$X$ if it is transitive and antisymmetric. The pair $(X,R)$ is called {\bf ordered space}.
\item \label{defi ordering antireflexive} $R$ is called an {\bf ordering in the sense of} "$<$" {\bf on}~$X$ if it is antireflexive and transitive. In this case we also write $<$ for~$R$, and $x < y$ or $y > x$ for $(x,y) \in R$. The pair $(X,<)$ is called {\bf space ordered in the sense of}~"$<$".
\item \label{defi ordering reflexive} $R$ is called an {\bf ordering in the sense of} "$\leq$" {\bf on}~$X$ if it is reflexive, antisymmetric, and transitive. In this case we also write $\leq$ for~$R$, and $x \leq y$ or $y \geq x$ for $(x,y) \in R$. The pair $(X,\leq)$ is called {\bf space ordered in the sense of}~"$\leq$".
\item \label{defi direction} $R$ is called a {\bf direction on}~$X$ if it is transitive, reflexive and directive. In this case we also write $\leq$ for~$R$, and $x \leq y$ for $(x,y) \in R$. The pair $(X,\leq)$ is called {\bf directed space}.
\eenum

\edefi

It follows by Definition~\ref{defi ordering} that there are no $x, y \in X$ such that both $x < y$ and $y < x$. Hence an ordering in the sense of "$<$" is antisymmetric. Therefore orderings in the sense of "$<$" and orderings in the sense of "$\leq$" are both orderings.

We remark that subsequently the symbol $\prec$ is used only for pre-orderings, $<$ and $>$ only for orderings in the sense of~"$<$", $\leq$ only for orderings in the sense of~"$\leq$" and for directions, and $\geq$ only for orderings in the sense of~"$\leq$". Regarding the symbol $\leq$ it is clarified in each case which kind of relation is considered. Of course, any of such relations may have additional properties and are still denoted by the same symbol. For example, an ordering may be denoted by~$\prec$ because it is a pre-ordering. Even an ordering in the sense of "$\leq$" or one in the sense of "$<$" may be denoted by~$\prec$ in certain cases. Thus by using the symbols we implicitly imply that the relation satisfies certain properties but we do not exclude that it satisfies more.

We deliberately use the aggregated notation $x \prec y \prec z$ instead of "$x \prec y$ and $y \prec z$", and similarly for the other two symbols.

\midvspace

\brema
\label{rema inv rel prop rel}
Let $(X,R)$ be a relational space and $A \subset X$. $R^{-1}$ and $R \, | \, A$ are pre-orderings, orderings, orderings in the sense of~"$<$", or orderings in the sense of~"$\leq$", if $R$ has the respective property.
\erema

\blemm
\label{lemm both orderings}
Let $(X,R)$ be an ordered space. Then $S = R \cup \Delta$ is an ordering in the sense of~$\leq$, and $T = R \setminus \Delta$ is an ordering in the sense of~$<$.
\elemm

\bproof
Exercise.
\eproof

\bdefi
\index{Successor}
\index{Predecessor}
Let $(X,\prec)$ be a pre-ordered space and $x \in X$. A point $y \in X$ is called {\bf successor of}~$x$ if $x \prec y$, $x \neq y$, and if there is no $z \in X \setminus\! \left\{ x, y \right\}$ such that $x \prec z \prec y$. A point $y \in X$ is called {\bf predecessor of}~$x$ if $y \prec x$, $x \neq y$, and if there is no $z \in X \setminus\! \left\{ x, y \right\}$ such that $y \prec z \prec x$.
\edefi

It is clear that generally a successor or a predecessor of a point $x \in X$ need not exist and if it exists, need not be unique. Obviously, if $y$ is a successor of~$x$, then $x$ is a predecessor of~$y$.

\midvspace

\bexam
Let $X = \left\{ a, b, c \right\}$ and $R = \left\{ (a,b), (b,c), (a,c), (b,b) \right\}$. Then $R$ is an ordering on~$X$, $R \setminus\! \left\{ (b,b) \right\}$ is an ordering in the sense of~"$<$", and $R \cup \left\{ (a,a), (c,c) \right\}$ is an ordering in the sense of~"$\leq$".
\eexam

Given a relational space, one can construct a pre-ordered space such that the set remains the same and the pre-ordering contains the original relation as subset. However this requires a recursive definition and is therefore postponed until Section~\ref{natural numbers}.

Given a pre-ordered space, one can construct an ordered space by an antisymmetrization procedure as follows.

\midvspace

\blemm
\label{lemm pre-ordering ordering}
Let $(X,R)$ be a pre-ordered space and $Q$ a relation on~$X$ defined by
\[
(x,y) \in Q \quad \Longleftrightarrow \quad (x,y), (y,x) \in R \;\; \wedge \;\; x = y
\]

Then $Q$ is an equivalence relation. Let $S \subset (X / Q) \times (X / Q)$ be the relation defined by
\[
(s,t) \in S \quad \Longleftrightarrow \quad \exists x \in s, \, y \in t \;\; (x,y) \in R
\]

Then $(X / Q, S)$ is an ordered space. If $R$ is reflexive, then $S$ is reflexive.
\elemm

\bproof
Exercise.
\eproof

Lemma~\ref{lemm pre-ordering ordering} is used in Theorem~\ref{theo Zorn}. The analysis of the equivalence relation~$Q$ also enhances our understanding how a general pre-ordered space, i.e.\ a set with a relation satisfying transitivity, looks like. It shows that there may be disjoint groups each consisting of several elements of~$X$ and "isolated elements" in the following sense. For every pair of distinct elements $x$ and $y$ within the same group we have $(x,y), (y,x) \in R$ and thus also $(x,x), (y,y) \in R$ by transitivity. For an "isolated element"~$x$ we may have $(x,x) \in R$ or $(x,x) \notin R$. The relation $S$ on $X / Q$ always leads to $(s,s) \in S$ if $s$ corresponds to a group of elements, and it may lead to $(s,s) \in S$ or to $(s,s) \notin S$ for "isolated elements" depending on which statement holds for the original elements of~$X$. Therefore the ordering $S$ need not be in the sense of~"$\leq$" nor in the sense of~"$<$". However, we have already examined another method how to construct such orderings from arbitrary orderings in Lemma~\ref{lemm both orderings}.

\midvspace

\blede
\label{lede subset ordering}
Let $X$ be a set and $\mc A \subset {\mc P}(X)$. Let the relation $R \subset {\mc A} \times {\mc A}$ be defined by $(A,B) \in R$ if $A \subset B$. We also write $({\mc A},\subset)$ for $({\mc A},R)$ and $\left( {\mc A}, \supset \right)$ for $\left( {\mc A}, R^{-1} \right)$. Each of the pairs $({\mc A},\subset)$ and $({\mc A},\supset)$ is an ordered space in the sense of "$\leq$", and if ${\mc A} = {\mc P}(X)$, a directed space.
\elede

\bproof
Exercise.
\eproof

\bdefi
\index{Total ordering}
\index{Ordering!total}
\index{Totally ordered space}
\index{Space!totally ordered}
Given a set $X$, a connective ordering $\prec$ on~$X$ is called {\bf total ordering}. In this case $(X,\prec)$ is called {\bf totally ordered space}.
\edefi

\bexam
Let $X$ be a set. If $X$ has more than one member, then $\left( {\mc P}(X), \subset \right)$ is not a totally ordered space.
\eexam

\bdefi
\label{defi interval}
\index{Proper interval}
\index{Interval!proper}
\index{Lower segment}
\index{Segment!lower}
\index{Upper segment}
\index{Segment!upper}
\index{Improper interval}
\index{Interval!improper}
\index{Interval}
Let $(X,\prec)$ be a pre-ordered space. For every $x,y \in X$ with $x \prec y$ the set $\, \left] x, y \right[ \, = \left\{ z \in X \, : \, x \prec z \prec y \right\}$ is called {\bf proper interval}. Moreover, for every $x \in X$, the set $\, \left] -\infty, x \right[ \,  =  \left\{ z \in X \, : \, z \prec x \right\}$ is called the {\bf lower segment of}~$x$, and the set $\, \left] x, \infty \right[ \, =  \left\{ z \in X \, : \, x \prec z \right\}$ is called the {\bf upper segment of}~$x$. A lower or upper segment is also called an {\bf improper interval}. A proper or improper interval is also called an {\bf interval}.
\edefi

Clearly, if $x \prec y$, then $\left] x, y \right[ \, = \, \left] -\infty, y \right[ \, \cap \, \left] x, \infty \right[ \,$. We remark that $\infty$ and $-\infty$ are merely used as symbols here. In particular, they do not generally refer to any of the number systems to be introduced below in this Chapter, neither does their usage imply that there is an infinite number of elements---for a Definition of "infinite" see Section~\ref{cardinality} below---in an improper interval.

\midvspace

\bdefi
Let $X$ be a set and ${\mc R} = \left\{ R_i \, : \, i \in I \right\}$ a system of pre-orderings on~$X$. Intervals with respect to a pre-ordering $R \in {\mc R}$ are denoted by subscript~$R$, i.e.\ $\,\left] -\infty, x \right[_{\,R}\,$ and $\,\left] x, \infty \right[_{\,R}\,$ where $x \in X$, and $\,\left] x, y \right[_{\,R}\,$ where $x, y \in X$, $(x,y) \in R$. Alternatively intervals with respect to $R_i$ for some $i \in I$ are denoted by index~$i$, i.e.\ $\,\left] -\infty, x \right[_{\,i}\,$, etc.\
\edefi

\brema
\label{rema order union}
Let $(X,\prec)$ be a pre-ordered space and $x \in X$.
\benum
\item If $\prec$ has full range, then $\, \left] -\infty, x \right[ \, = \, {\displaystyle \bigcup} \, \big\{ \, \left] y, x \right[ \; : \, y \in X, \; y \prec x \big\}$
\item If $\prec$ has full domain, then $\, \left] x, \infty \right[ \, = \, {\displaystyle \bigcup} \, \big\{ \, \left] x, y \right[ \; : \, y \in X, \; x \prec y \big\}$
\eenum

\erema

\bdefi
\label{defi order dense}
\index{Dense}
\index{Order dense}
Let $(X,\prec)$ be a pre-ordered space. A subset $Y \subset X$ is called {\bf $\prec$-dense in}~$X$ or {\bf order dense in}~$X$ if for every $x, y \in X$ with $x \prec y$ there exists $z \in Y$ such that $x \prec z \prec y$. $X$ is called {\bf $\prec$-dense} or {\bf order dense} if it is order dense in itself.
\edefi

\bdefi
\label{defi gen order dense}
\index{R-dense}
\index{Dense}
Let $X$ be a set and $\mc R$ a system of pre-orderings on~$X$. A subset $Y \subset X$ is called {\bf $\mc R$-dense in}~$X$ if for every $R \in {\mc R}$ and $x, y \in X$ with $(x,y) \in R$ there exists $z \in Y$ such that $(x,z), (z,y) \in R$. $X$~is called {\bf $\mc R$-dense} if it is $\mc R$-dense in itself.
\edefi

\brema
\label{rema ordered dense}
Let $(X,\prec)$ be a pre-ordered space and $Y \subset X$ order dense. For every $x, y \in X$ the following equalities hold:
\begin{eqnarray*}
\left] -\infty, y \right[  \!\! & = & \!\! \bigcup \big\{ \, \left] -\infty, z \right[ \; : \, z \in Y, \; z \prec y \big\}\\
\left] x, \infty \right[ \!\! & = & \!\! \bigcup \big\{ \, \left] z, \infty \right[ \; : \, z \in Y, \; x \prec z \big\}\\
\left] x, y \right[ \!\! & = & \!\! \bigcup \big\{ \, \left] u, v \right[ \; : \, u, v \in Y, \; x \prec u \prec v \prec y \big\}
\end{eqnarray*}

\erema

\bdefi
\label{defi weak minimum}
\index{Weak minimum}
\index{Minimum!weak}
\index{Weak maximum}
\index{Maximum!weak}
Let $(X,R)$ be a relational space. A member $x \in X$ is called a {\bf weak minimum of}~$X$ if $(y,x) \in R$ implies $(x,y) \in R$. Moreover a member $x \in X$ is called a {\bf weak maximum of}~$X$ if $(x,y) \in R$ implies $(y,x) \in R$. Further let $A \subset X$. Then $x \in A$ is called a {\bf weak minimum} ({\bf weak maximum}) {\bf of} $A$ if it is a weak minimum (weak maximum) of~$A$ with respect to the restriction~$R \, | \, A$.
\edefi

\midvspace

\bdefi
\label{defi minimum}
\index{Minimum}
\index{Least member}
\index{Maximum}
\index{Greatest member}
Let $(X,R)$ be a relational space. A member $x \in X$ is called a {\bf minimum} or {\bf least element of}~$X$ if $(x,y) \in R$ for every $y \in X \!\setminus\! \left\{ x \right\}$. Moreover a member $x \in X$ is called a {\bf maximum} or {\bf greatest element of}~$X$ if $(y,x) \in R$ for every $y \in X \!\setminus\! \left\{ x \right\}$. Further let $A \subset X$. Then $x \in A$ is called a {\bf minimum} ({\bf maximum}) {\bf of}~$A$ if it is a minimum (maximum) of~$A$ with respect to the restriction~$R \, | \, A$. If $A$ has a unique minimum (maximum), then it is denoted by $\mathrm{min} \, A$~($\mathrm{max} \, A$). 
\edefi

Notice that the singleton~$\left\{ x \right\}$, where $x$ is a set, trivially has $x$ as its minimum and maximum. Although Definitions \ref{defi weak minimum} and~\ref{defi minimum} are valid for any relation $R$ on~$X$, they are mainly relevant in the case where $R$ is a pre-ordering.

The following result shows that the notions defined in Definition~\ref{defi minimum} are invariant under a change from the original relation to the relations defined in Lemma~\ref{lemm both orderings}.

\midvspace

\blemm
\label{lemm max min inv}
Let $(X,R)$ be a relational space, $T \in \left\{ R \cup \Delta,\, R \!\setminus\! \Delta \right\}$, and $x \in X$. $x$~is a minimum (maximum) of~$X$ with respect to~$T$ iff it is a minimum (maximum) of~$X$ with respect to~$R$.
\elemm

\bproof
Exercise.
\eproof

\brema
Let $(X,R)$ be a relational space. If $x \in X$ is a minimum (maximum) of~$X$, then $x$ is also a weak minimum (weak maximum) of~$X$.
\erema

\brema
Let $(X,<)$ be an ordered space and $x$ a weak minimum of~$X$. Then there is no $y \in X$ with $y < x$.
\erema

\brema
\label{rema max unique}
Let $(X,\prec)$ be an ordered space. If $X$ has a minimum (maximum), then this minimum (maximum) is unique.
\erema

\brema
Let $(X,\prec)$ be a totally ordered space. If $X$ has a weak minimum (weak maximum), then this weak minimum (weak maximum) is the minimum (maximum) of~$X$.
\erema

\bdefi
\index{Minimum property}
Let $(X,R)$ be a relational space. We say that $R$ {\bf has the minimum property} if every $A \subset X$ with $A \neq \O$ has a minimum.
\edefi

\brema
\label{prop minimum property connective}
Let $(X,R)$ be a relational space. If $R$ has the minimum property, then $R$ is connective.
\erema

\bdefi
\label{defi well-ordering}
\index{Well-ordering}
\index{Well-ordered space}
\index{Space!well-ordered}
Given a set $X$, an ordering $R$ on $X$ that has the minimum property is called {\bf well-ordering}. In this case we say that $R$ {\bf well-orders} $X$, and $(X,R)$ is called {\bf well-ordered space}.
\edefi

Notice that according to Definition \ref{defi well-ordering} a well-ordering may be an ordering in the sense of "$<$" or "$\leq$" or neither. In the literature "well-ordering" is often used only in the sense of "$<$" (see e.g.\ \cite{Kelley} or~\cite{Ebbinghaus}) or only in the sense of "$\leq$" (see e.g.~\cite{Querenburg}).

\midvspace

\blemm
Every well-ordered space is totally ordered.
\elemm

\bproof
This follows from Remark~\ref{prop minimum property connective}.
\eproof

\blemm
\label{lemm both well-orderings}
Let $(X,R)$ be a well-ordered space. Then $S = R \cup \Delta$ is a well-ordering in the sense of~"$\leq$", and $T = R \setminus \Delta$ is a well-ordering in the sense of~"$<$". 
\elemm

\bproof
$S$ and $T$ are clearly well-orderings. The claim follows by Lemma \ref{lemm both orderings}.
\eproof

\blemm
Let $(X,R)$ be a relational space. If $R$ is antisymmetric and has the minimum property, then it is a well-ordering.
\elemm

\bproof
$R$ is connective by Lemma \ref{prop minimum property connective}. Now assume that $R$ is not transitive. Let $x, y, z \in X$ such that $(x,y), (y,z) \in R$ and $(x,z) \notin R$. If all three or any two of the elements of $\left\{ x, y, z \right\}$ are equal, then this is a contradiction. If $x$, $y$, and $z$ are distinct, then we have $(z,x) \in R$ since $R$ is connective, and $(y,x), (z,y) \notin R$ since $R$ is antisymmetric. Thus $\left\{ x, y, z \right\}$ has no minimum, which is a contradiction.
\eproof

\bdefi
\label{defi bound sup inf}
\index{Upper bound}
\index{Bound!upper}
\index{Lower bound}
\index{Bound!lower}
\index{Supremum}
\index{Least upper bound}
\index{Bound!least upper}
\index{Infimum}
\index{Greatest lower bound}
\index{Bound!greatest lower}
Let $(X,\prec)$ be a pre-ordered space and $A \subset X$. A member $x \in X$ is called an {\bf upper bound of}~$A$ if $y \prec x$ for every $y \in A \!\setminus\! \left\{ x \right\}$. A member $x \in X$ is called a {\bf lower bound of}~$A$ if $x \prec y$ for every $y \in A \!\setminus\! \left\{ x \right\}$. A member $x \in X$ is called a {\bf supremum of}~$A$ or a {\bf least upper bound of}~$A$ if it is a minimum of the set of all upper bounds of~$A$. A member $x \in X$ is called an {\bf infimum of}~$A$ or {\bf a greatest lower bound of}~$A$ if it is a maximum of the set of all lower bounds of~$A$. If $A$ has a unique supremum (infimum), then it is denoted by $\sup A$~($\inf A$).
\edefi

Again, the notions defined in Definition~\ref{defi bound sup inf} are invariant under a change from the original relation to the relations defined in Lemma~\ref{lemm both orderings}.

\midvspace

\blemm
\label{lemm bound sup inv}
Let $(X,R)$ be a relational space, $T \in \left\{ R \cup \Delta,\, R \!\setminus\! \Delta \right\}$, $A \subset X$, and $x \in X$. $x$~is an upper bound, lower bound, supremum, or infimum of~$A$ with respect to~$T$ iff it has the respective property with respect to~$R$.
\elemm

\bproof
Exercise.
\eproof

\blede
\label{lede ordered sup unique}
Let $(X,\prec)$ be an ordered space. Every $A \subset X$ has at most one supremum and at most one infimum. Given a set $Y$, a subset $B \subset Y$, and a function $f : Y \longrightarrow X$, the supremum of the set $f \left[ B \right] = \left\{ f(y) \, : \, y \in B \right\}$ is also denoted by $\sup_{y \in B} f(y)$ and its infimum by $\inf_{y \in B} f(y)$.
\elede

\bproof
Let $C$ be the set of all upper bounds of~$A$. If $C$ has a minimum, then this minimum is unique by Remarks~\ref{rema inv rel prop rel} and~\ref{rema max unique}. Therefore $A$ has at most one supremum. The proof regarding the minimum is similar.
\eproof

The following is a property that, for example, the real numbers have as demonstrated in Lemma~\ref{lemm least upper bound reals} below.

\midvspace

\bdefi
\label{defi least upper bound property}
\index{Least upper bound property}
Let $(X,\prec)$ be a pre-ordered space. We say that $\prec$ {\bf has the least upper bound property} if every set $A \subset X$ with $A \neq \O$ which has an upper bound has a supremum.
\edefi

The least upper bound property is equivalent to the intuitively reversed property as stated in the following Theorem. In the proof we follow~\cite{Kelley}, p.~14.

\midvspace

\btheo
\label{theo least upper greatest lower}
Let $(X,\prec)$ be a pre-ordered space. $\prec$~has the least upper bound property iff every set $A \subset X$ with $A \neq \O$ which has a lower bound has an infimum.
\etheo

\bproof
First assume that $\prec$ has the least upper bound property. Let $A \subset X$ such that $A \neq \O$ and $A$ has a lower bound. Further let $B$ be the set of all lower bounds of~$A$. Let $x \in A$. It follows that, for every $y \in B$, we have $y = x$ or $y \prec x$. Hence $x$ is an upper bound of~$B$. Therefore all members of~$A$ are upper bounds of~$B$. By assumption $B$ has a supremum, say~$y$. Since $y$ is the minimum of all upper bounds of~$B$, we have $y \prec x$ for every $x \in A \!\setminus\! \left\{ y \right\}$. Thus $y$ is a lower bound of~$A$. In order to see that $y$ is the greatest lower bound of~$A$, let $z$ be a lower bound of~$A$, i.e.\ $z \in B$. Since $y$ is an upper bound of~$B$ we have $z \prec y$ or $z = y$.

The converse can be proven similarly.
\eproof

\bexam
\label{exam sup inf}
Let $X$ be a set. Then $({\mc P}(X),\subset)$ is an ordered space. Further let ${\mc A} \subset {\mc P}(X)$. Then $X$ is an upper bound of~$\mc A$ and $\O$ is a lower bound of~$\mc A$. If ${\mc A} \neq \O$, then $\bigcup {\mc A}$ is the supremum of~$\mc A$, and $\bigcap {\mc A}$ is the infimum of~$\mc A$. Thus the relation $\subset$ on ${\mc P}(X)$ has the least upper bound property.
\eexam

\blede
\label{defi structure relation}
\index{Structure relation}
\index{Relation!structure}
Given a set $X$, we define
\begin{eqnarray*}
{\mc Q}(X) \!\!\!\! & = & \!\!\!\! \big\{ \, (x,A) \in X \!\times {\mc P}(X) \, : \, x \in A \, \big\}\\ 
 & = & \!\!\! \bigcup \big\{ \left\{ x \right\} \!\times\! \big({\mc P}(X)(x)\big) \, : \, x \in X \, \big\} \, \subset \, X \!\times {\mc P}(X)
\end{eqnarray*}

A relation $R \subset {\mc Q}(X)$ is called a {\bf structure relation on}~$X$. The relation $\leq$ on $R$ defined by
\[
(y,B) \leq (x,A) \quad \Longleftrightarrow \quad x = y \; \wedge \; A \subset B
\]

is an ordering in the sense of~"$\leq$". If $X$ has more than one member, this ordering is not connective.
\elede

\bproof
Exercise.
\eproof

\section{Functions}
\label{functions}

In this Section we introduce the important concept of function. We analyse various fundamental properties of functions, in particular the interplay of functions with unions and intersections of sets as well as with  pre-orderings.

\midvspace

\bdefi
\index{Functional relation}
\index{Relation!functional}
\index{Function}
\index{Map}
\index{Image}
\index{Inverse}
Let $X$ and $Y$ be two sets. A {\bf functional relation} $f$ is a relation $f \subset X \!\times Y$ such that for every $x \in X$ there exists at most one $y \in Y$ with $(x,y) \in f$. Let $D$ be the domain of $f$. Then $f$ is called a {\bf function from} $D$ {\bf to} $Y$. We use the standard notation $f: D \longrightarrow Y$. A function is also called {\bf map} in the sequel. For every $x \in D$ we denote by $f(x)$ or $f_x$ the member $y \in Y$ such that $(x,y) \in f$. $f(x)$~is called the {\bf value of~$f$ at~$x$}. For every $A \subset X$ we call $f \left[ A \right] = \left\{ f(x) \, : \, x \in A \right\}$ the {\bf image of} $A$ {\bf under}~$f$. For every $B \subset Y$ the set $f^{-1} \left[ B \right] = \left\{ x \in X \, : \, f(x) \in B \right\}$ is called the {\bf inverse of} $B$ {\bf under}~$f$. For a system ${\mc A} \subset {\mc P}(X)$ the system $f \, \llbracket {\mc A} \, \rrbracket = \left\{ f \left[ A \right] \, : \, A \in {\mc A} \right\}$ is called {\bf image of} ${\mc A}$ {\bf under}~$f$. For a system ${\mc B} \subset {\mc P}(Y)$ the system $f^{-1} \, \llbracket {\mc B} \, \rrbracket = \left\{ f^{-1} \left[ B \right] \, : \, B \in {\mc B} \right\}$ is called {\bf inverse of} ${\mc B}$ {\bf under}~$f$.
\edefi

\bdefi
\label{defi function set}
Let $X$ and $Y$ be two sets. The set of all functions from $X$ to $Y$ is denoted by~$Y^X$.
\edefi

Notice that slightly different definitions are used if $X$ is a natural number (see Definition~\ref{defi exp nat number}) or if $Y$ is a relation (see Definition~\ref{def rel npower}). Generally there is no risk of confusion.

\midvspace

\bdefi
Let $X$, $Y$, and $Z$ be sets, $f : X \!\times Y \longrightarrow Z$ a function, and $x \in X$, $y \in Y$. Then we also write $f(x, y)$ instead of $f((x, y))$ for the value of~$f$ at~$(x, y)$.
\edefi

\blede
\index{Surjective}
\index{Function!surjective}
\index{Injective}
\index{Function!injective}
\index{Bijective}
\index{Function!bijective}
\index{Inverse}
\index{Identity map}
Given two sets $X$ and $Y$, and a function $f: X \longrightarrow Y$, $f$ is called {\bf surjective} if $f \left[ X \right] = Y$, i.e.\ the range of $f$ is~$Y$. $f$~is called {\bf injective} if $f^{-1} \left\{ y \right\}$ contains at most one member for each $y \in Y$. $f$~is called {\bf bijective} if $f$ is both surjective and injective.

If $f$ is bijective, then the inverse relation $f^{-1}$ is a functional relation with domain~$Y$, i.e.\ $f^{-1} : Y \longrightarrow X$. $f^{-1}$~is called {\bf inverse function of}~$f$, or short, {\bf inverse of}~$f$. We have $f^{-1} \left(f(x)\right) = x$ for every $x \in X$. $f^{-1}$ is bijective.

If $X = Y$ and $f = \Delta$, then $f$ is called the {\bf identity map on} $X$, and also denoted by~$\mathrm{id}_X$, or, when the set is evident from the context, by~$\mathrm{id}$.
\elede

\bproof
Exercise.
\eproof

\brema
\label{rema bijection systems}
Let $X$, $Y$ be two sets, $f: X \longrightarrow Y$ a bijection, ${\mc A} \subset {\mc P}(X)$, and ${\mc B} = f \, \llbracket {\mc A} \, \rrbracket \,$. Then the map $F : {\mc A} \longrightarrow {\mc B}$, $F(A) = f \left[ A \right]$ is a bijection too.
\erema

\bdefi
\index{Fixed point}
Given a set $X$ and a map $f : X \longrightarrow X$, a member $x \in X$ is called {\bf fixed point of}~$f$ if $f(x) = x$.
\edefi

\bdefi
\index{Restriction!function}
\index{Function!restriction}
Given two sets $X$, $Y$, a function $f: X \longrightarrow Y$, and a set $A \subset X$, the functional relation $\left\{ (x,y) \in f \, : \, x \in A \right\}$ is called {\bf restriction of $f$ to}~$A$. It is denoted by~$f \, | \, A$.
\edefi

\blede
\index{Composition}
\index{Function!composition}
Given sets $X$, $Y$ and $Z$, and functions $f: X \longrightarrow Y$ and $g: Y \longrightarrow Z$, the product $gf$ of $g$ and $f$ as defined in Definition \ref{defi rel inv prod} is also denoted by $g \circ f$. It is a function from $X$ to $Z$, i.e.\ $g \circ f : X \longrightarrow Z$. It is also called the {\bf composition of $f$ and~$g$}. We have $g(f(x)) = (gf)(x)$ for every $x \in X$. We also write $gf(x)$ for $(gf)(x)$.
\elede

\bproof
Exercise.
\eproof

\bdefi
\index{Projection}
\index{Projective}
\index{Function!projective}
Given a set $X$ and a map $f : X \longrightarrow X$, $f$ is called a {\bf projection} or {\bf projective}, if $f \circ f = f$.
\edefi

Using the notion of a function, a system of sets and the union and intersection of a system as defined in Lemma and Definition~\ref{lede great union} and Definition~\ref{defi great intersection}, respectively, can be written in a different form as follows.

\midvspace

\bdefi
\label{defi index set}
\index{Index set}
\index{Set!index}
A set $I$ is called an {\bf index set} if $I \neq \O$. Given a non-empty system $\mc A$, an index set~$I$, and a function $A : I \longrightarrow {\mc A}$, we define the following notations:
\[
\bigcup_{i \in I} A_i \, = \, \bigcup {\mc B}, \quad \quad \quad \bigcap_{i \in I} A_i \, = \, \bigcap {\mc B}
\]

where ${\mc B} = A \left[ I \right]$. If $A$ is surjective, then it follows that
\[
\bigcup_{i \in I} A_i \, = \, \bigcup {\mc A}, \quad \quad \quad \bigcap_{i \in I} A_i \, = \, \bigcap {\mc A}
\]

\edefi

We mainly use the notion "index set" for a set that is the non-empty domain of a function to a system of sets as in Definition~\ref{defi index set}, but not for arbitrary non-empty sets; of course, formally also the system $\mc A$ is an index set. With the notation of Definition~\ref{defi index set} we clearly have ${\mc A} = \left\{ A_i \, : \, i \in I \right\}$ if $A$ is surjective. It is often more convenient to use the index notations than an abstract letter for the system of sets. Notice that there is a slight difference between the two notations because we may have $A_i = A_j$ for $i, j \in I$ with $i \neq j$. This happens if the map $A$ is not injective. However, in most cases this turns out to be irrelevant. When using index notation, we often do neither explicitly introduce a letter for the range system~(e.g.~$\mc A$) nor a letter for the function (e.g.~$A$). Instead we only introduce an index set~$I$ and the sets $A_i$ ($i \in I$) that specify the values of the function and that are precisely the members of the range system. In particular, a definition of the index set~$I$ and the sets $A_i$ ($i \in I$) does not tacitly imply that the letter $A$ without subscript stands for the corresponding function unless this is explicitly said; we may even use the letter~$A$ for other purposes, for example we may define $A = \bigcup_{i \in I} A_i$\,.

The following identities are the analogues of Lemmas~\ref{lemm set equalities great} and~\ref{lemm demorgan formulas great}.

\midvspace

\blemm
\label{lemm set equalities index}
Let $I$ and $J$ be index sets and $A_i$ ($i \in I$), $B_j$ ($j \in J$) sets. Then the following equalities hold:
\benum
\item \label{lemm set equalities index 1} $\big( \bigcup_{i \in I} A_i \big) \cap \big( \bigcup_{j \in J} B_j \big) = \, {\displaystyle \bigcup} \, \big\{ A_i \cap B_j \, : \, i \in I,\, j \in J \big\}$
\item \label{lemm set equalities index 2} $\big( \bigcap_{i \in I} A_i \big) \cup \big( \bigcap_{j \in J} B_j \big) = \, {\displaystyle \bigcap} \, \big\{ A_i \cup B_j \, : \, i \in I,\, j \in J \big\}$
\item \label{lemm set equalities index 3} $\big( \bigcap_{i \in I} A_i \big)^c = \, \bigcup_{i \in I} A_i^c$
\item \label{lemm set equalities index 4} $\big( \bigcup_{i \in I} A_i \big)^c = \, \bigcap_{i \in I} A_i^c$
\eenum

\elemm

\bproof
Exercise.
\eproof

The following two Lemmas show how the image and the inverse under~$f$ behave together with intersections and unions.

\midvspace

\blemm
\label{lemm func relations}
Given a function $f: X \longrightarrow Y$, an index set~$I$, and sets $A_i \subset X$ ($i \in I$), the following statements hold:
\benum
\item \label{lemm func relations 1} $f \left[ \, \bigcup_{i \in I} A_i \, \right] \, = \, \bigcup_{i \in I} f \left[ A_i \right]$
\item \label{lemm func relations 2} $f \left[ \, \bigcap_{i \in I} A_i \, \right] \, \subset \, \bigcap_{i \in I} f \left[ A_i \right]$
\eenum

\elemm

\bproof
Exercise.
\eproof

\blemm
Given a function $f: X \longrightarrow Y$, an index set~$I$, a set $A \subset Y$, and sets $A_i \subset Y$ ($i \in I$), the following relations hold:
\benum
\item \label{inverse1} $f^{-1} \left[ \, \bigcup_{i \in I} A_i \, \right] \, = \, \bigcup_{i \in I} f^{-1} \left[ A_i \right]$
\item \label{inverse2} $f^{-1} \left[ \, \bigcap_{i \in I} A_i \, \right] \, = \, \bigcap_{i \in I} f^{-1} \left[ A_i \right]$
\item \label{inverse3} $f^{-1} \left[ A^c \right] \, = \left( f^{-1} \left[ A \right] \right)^c$
\eenum

where the complement refers to~$Y$.
\elemm

\bproof
Exercise.
\eproof

The following Definition generalizes the Cartesian product of two sets as defined in Definition~\ref{defi cartesian prod pair}.

\midvspace

\bdefi
\label{defi Cartesian product}
\index{Cartesian product}
\index{Projection}
Let $I$ be an index set, $\mc A$ a non-empty system, $A : I \longrightarrow {\mc A}$ a map, and \mbox{$B = \bigcup_{i \in I} A_i$}. We define the {\bf Cartesian product of~$A$} as follows:
\[
\bigtimes_{\!\! i \in I}\, A_i \, = \, \left\{ f \in B^I \; : \; \forall i \in I \;\; f(i) \in A_i \right\}
\]

For each $i \in I$, the map $p_i : \bigtimes_{\!\! i \in I}\, A_i \longrightarrow A_i$, $p_i(f) = f(i)$, is called the {\bf projection on}~$A_i$.
\edefi

Notice that in our definition of the Cartesian product we use index notation, which allows identical factors. Of course, using index notation for the projections we formally have to think of a surjective map $p : I \longrightarrow \left\{ p_i \, : \, i \in I \right\}$.

It is a consequence of the Choice axiom that the Cartesian product is not always empty. The following Remark is a repetition of Lemma and Definition~\ref{lede choice function}, now using functional notation.

\midvspace

\brema
\label{rema choice}
Let $X$ be a set and ${\mc A} \subset {\mc P}(X)$ with $\O \notin {\mc A}$. Then there exists a function $f : {\mc A} \longrightarrow X$ such that $f(A) \in A$ for every $A \in {\mc A}$. $f$~is a choice function.
\erema

\bcoro
With definitions as in Definition~\ref{defi Cartesian product}, $\O \notin {\mc A}$ implies that $\bigtimes_{\!\! i \in I}\, A_i \neq \O$. 
\ecoro

\bproof
This is a consequence of Remark~\ref{rema choice}.
\eproof

\brema
With definitions as in Definition~\ref{defi Cartesian product}, the following statements hold:
\benum
\item If $A_i = A$ ($i \in I$) for some set~$A$, then $\bigtimes_{\!\! i \in I}\, A_i = A^I$.
\item If $A_i = \O$ for some $i \in I$, then $\bigtimes_{\!\! i \in I}\, A_i = \O$.
\eenum
\erema

If the index set in Definition~\ref{defi Cartesian product} is a singleton, the Cartesian product can obviously be identified with the single factor set in the following manner.

\midvspace

\brema
\label{rema singleton cartesian}
Let $X$ and $a$ be two sets and $I = \left\{ a \right\}$. We define the map $f : X^I \longrightarrow X$, $f(h) = h(a)$. Then $f$ is bijective.
\erema

The following Remark says that Definitions~\ref{defi cartesian prod pair} and~\ref{defi Cartesian product} are in agreement with each other.

\midvspace

\brema
\label{rema square set}
Let $X$, $Y$, $a$, and $b$ be sets with $a \neq b$. We further define the sets $I = \left\{ a, b \right\}$, $X_a = X$, $X_b = Y$, and the function $f : \bigtimes_{\!\! i \in I}\, X_i \longrightarrow X \!\times Y$, $f(h) = \left( h(a), h(b) \right)$. Then $f$ is bijective. In particular, this gives us a bijection from $X \!\times X$ to~$X^I$.
\erema

The following result says that an iterated Cartesian product can be identified with a simple Cartesian product.

\midvspace

\brema
\label{rema it Cartesian product}
Let $\mc J$ be a non-empty system of disjoint index sets, $J : I \longrightarrow {\mc J}$ a bijection where $I$ is an index set, $K = \bigcup {\mc J}$, $\mc A$~a non-empty system, and $F : K \longrightarrow {\mc A}$ a map. Then the Cartesian product $\bigtimes_{\!\! j \in K}\, F_j$ is well-defined. Further let $G : I \longrightarrow {\mc P} \left( K \!\times\! {\mc A} \right)$ be the map such that, for every $i \in I$, $G(i)$ is a functional relation with domain $J_i$, i.e.\ $G_i : J_i \longrightarrow {\mc A}$, and $G_i (j) = F(j)$. Thus, for every $i \in I$, the Cartesian product $\bigtimes_{\!\! j \in J(i)}\, G_i (j)$ is well-defined. Now let $A = \bigcup {\mc A}$ and $H : I \longrightarrow {\mc P}^2 \left( K \!\times\! A \right)$, $H(i) = \bigtimes_{\!\! j \in J(i)}\, G_i (j)$. We define
\begin{center}
\begin{tabular}{l}
$f : \bigtimes_{\!\! j \in K}\, F_j \longrightarrow \bigtimes_{\!\! i \in I}\, H_i \;$,\\[.8em]
$\Big(\big(f(h)\big)(i)\Big)(j) = h(j)$ \quad for every $i \in I$ and $j \in J_i \;$
\end{tabular}
\end{center}

Then $f$ is a bijection.
\erema

Remarks~\ref{rema singleton cartesian}, \ref{rema square set}, and~\ref{rema it Cartesian product} can be combined in different ways. The following is a useful example.

\midvspace

\brema
\label{rema cartesian mixed}
Let $I$ be an index set, $X_i$ ($i \in I$) sets, $j \in I$, and $J = I \!\setminus\! \left\{j \right\}$. If $J \neq \O$, then there is a bijection from $\bigtimes_{\!\! i \in I}\, X_i$ to $\left( \bigtimes_{\!\! i \in J}\, X_i \right) \!\times\! X_j$ by Remarks~\ref{rema singleton cartesian}, \ref{rema square set}, and~\ref{rema it Cartesian product}.
\erema

\bdefi
\index{Distinguishes points}
Let $X$ be a set and $I$ an index set. Further let $Y_i$ ($i \in I$) be sets and $f_i : X \longrightarrow Y_i$ maps. We say that $\left\{ f_i \, : \, i \in I \right\}$ {\bf distinguishes points} if for every $x, y \in X$ with $x \neq y$ there is $i \in I$ such that $f_i(x) \neq f_i(y)$.
\edefi

\brema
With definitions as in Definition~\ref{defi Cartesian product}, the set of functions $\left\{ p_i \, : \, i \in I \right\}$ distinguishes points.
\erema

\bdefi
\index{Bounded function}
\index{Function!bounded}
\index{Unbounded function}
\index{Function!unbounded}
\index{Bounded from below}
\index{Function!bounded from below}
\index{Bounded from above}
\index{Function!bounded from above}
Let $X$ be a set, $(Y,\prec)$ a pre-ordered space, and $f: X \longrightarrow Y$ a function. $f$~is called {\bf bounded} if there exist $x, y \in Y$ such that $f \left[ X \right] \subset \, \left] x, y \right[ \, \cup \left\{ x, y \right\}$. Otherwise $f$ is called {\bf unbounded}. $f$~is called {\bf bounded from below} if there is $x \in Y$ such that $f \left[ X \right] \subset \, \left] x, \infty \right[ \, \cup \left\{ x \right\}$. $f$~is called {\bf bounded from above} if there is $y \in Y$ such that $f \left[ X \right] \subset \, \left] -\infty, y \right[ \, \cup \left\{ y \right\}$.
\edefi

\blemm
\label{lemm two func inf sup}
Let $(X,\prec)$ be an ordered space where $\prec$ has the least upper bound property, $Y$ a non-empty set, and $f : Y \longrightarrow X$, $g : Y \longrightarrow X$ two functions such that $f(y) \prec g(y)$ for every $y \in Y$. The following two statements hold:
\benum
\item \label{lemm two func inf sup 1} If $f$ is bounded from below, then \, $\inf f \left[ Y \right] \prec \inf g \left[ Y \right]$ \, or \\ $\inf f \left[ Y \right] = \inf g \left[ Y \right]$.
\item \label{lemm two func inf sup 2} If $g$ is bounded from above, then \, $\sup f \left[ Y \right] \prec \sup g \left[ Y \right]$ \, or \\ $\sup f \left[ Y \right] = \sup g \left[ Y \right]$.
\eenum
\elemm

\bproof
In order to prove~(\ref{lemm two func inf sup 1}), let $L_f$ be the set of all lower bounds of~$f \left[ Y \right]$ and $L_g$ the set of all lower bounds of~$g \left[ Y \right]$. Under the stated conditions we have $L_f \subset L_g$. Since $f \left[ Y \right]$ has a lower bound, it has an infimum by Theorem~\ref{theo least upper greatest lower}. Moreover, the infimum of $f \left[ Y \right]$ is unique since $\prec$ is an ordering. Similarly, also $\inf g \left[ Y \right]$ exists. The claim now follows by the fact that $L_f \subset L_g$.

The proof of~(\ref{lemm two func inf sup 2}) is similar.
\eproof

\bdefi
\label{defi monotonic function}
\index{Monotonically increasing}
\index{Function!monotonically increasing}
\index{Increasing}
\index{Function!increasing}
\index{Non-decreasing}
\index{Function!non-decreasing}
\index{Monotonically decreasing}
\index{Function!monotonically decreasing}
\index{Decreasing}
\index{Function!decreasing}
\index{Non-increasing}
\index{Function!non-increasing}
\index{Strictly increasing}
\index{Function!strictly increasing}
\index{Strictly decreasing}
\index{Function!strictly decreasing}
\index{Strictly monotonic}
\index{Function!strictly monotonic}
Let $(X,\prec)$ and $(Z,\prec)$ be ordered spaces, $A \subset X$, $B \subset Z$, and $f : A \longrightarrow B$ a function. $f$~is called {\bf monotonically increasing} or {\bf increasing} or {\bf non-decreasing} if, for every $x, y \in A$, $x \prec y$ implies $f(x) \prec f(y)$ or $f(x) = f(y)$. $f$~is called {\bf monotonically decreasing} or {\bf decreasing} or {\bf non-increasing} if, for every $x, y \in A$, $x \prec y$ implies $f(y) \prec f(x)$ or $f(x) = f(y)$. $f$~is called {\bf monotonic} if it is either increasing or decreasing.

Further, $f$~is called {\bf strictly increasing} if, for every $x, y \in A$ with $x \neq y$, $x \prec y$ implies $f(x) \prec f(y)$ and $f(x) \neq f(y)$. $f$~is called {\bf strictly decreasing} if, for every $x, y \in A$ with $x \neq y$, $x \prec y$ implies $f(y) \prec f(x)$ and $f(x) \neq f(y)$. $f$~is called {\bf strictly monotonic} if it is either strictly increasing or strictly decreasing.
\edefi

The following result shows that the notions defined in Definition~\ref{defi monotonic function} are invariant under the change from the original orderings to the orderings in the sense of~"$<$" and~"$\leq$" as defined in Lemma~\ref{lemm both orderings}, both in the domain and in the range space.

\midvspace

\blemm
\label{lemm mon invar}
Let $(X,R)$ and $(Y,S)$ be ordered spaces, and $f : X \longrightarrow Y$ a function. Further let \mbox{$T \in \left\{ R \cup \Delta,\, R \setminus \Delta \right\}$} and \mbox{$U \in \left\{ S \cup \Delta,\, S \setminus \Delta \right\}$}. $f$~is increasing, decreasing, strictly increasing, or strictly decreasing with respect to the orderings $T$ on~$X$ and $U$ on~$Y$ iff it has the respective property with respect to the orderings $R$ and~$S$.
\elemm

\bproof
Exercise.
\eproof

\bdefi
\index{Binary function}
\index{Function!binary}
\index{Associative!binary function}
\index{Function!associative binary}
\index{Commutative!binary function}
\index{Function!commutative binary}
\index{Group}
\index{Abelian group}
\index{Group!Abelian}
Given a set $X$, a function $f : X \!\times\! X \longrightarrow X$ is called {\bf binary function on}~$X$. The symbol~$f$ is also denoted by~$\bullet$ \,, and for every $x, y \in X$ we also write $x \bullet y$ for $f(x, y)$. If the equality $x \bullet (y \bullet z) = (x \bullet y) \bullet z$ holds for every $x, y, z \in X$, the function is called {\bf associative}. If the equality $x \bullet y = y \bullet x$ holds for every $x, y \in X$, the function is called {\bf commutative}.

The triple $(X, f, e)$---or $(X, \bullet, e)$---with $e \in X$ is called a {\bf group} if the following statements hold:
\benum
\item The function $\bullet$ is associative.
\item For every $x \in X$, we have $x \bullet e = x$. $e$~is called {\bf neutral element}.
\item For every $x \in X$, there is $y \in X$ such that $x \bullet y = e$. $y$~is called {\bf inverse of}~$x$.
\eenum

If $\bullet$ is commutative, the group is called~{\bf Abelian}.
\edefi

In the remainder of the text various binary functions are introduced and different symbols are defined, for example $+$ is used instead of $\bullet$ for the addition of natural numbers.

\midvspace

\btheo
\label{lemm group properties}
Given a group $(X, \bullet, e)$, the following statements hold:
\benum
\item \label{lemm group properties 1} For every $x \in X$, we have $e \bullet x = x$.
\item \label{lemm group properties 2} There is no $d \in X \setminus\! \left\{ e \right\}$ such that $x \bullet d = x$ for every $x \in X$.
\item \label{lemm group properties 3} For every $x \in X$, there is a unique $y \in X$ such that $x \bullet y = y \bullet x = e$.
\eenum

\etheo

\bproof
We first show that $x \bullet y = e$ implies $y \bullet x = e$ for every $x, y \in X$. Assume $x \bullet y = e$. There is $z \in X$ such that $y \bullet z = e$. It follows that $y \bullet x = y \bullet (x \bullet e) = y \bullet (x \bullet (y \bullet z)) = y \bullet ((x \bullet y) \bullet z) = y \bullet (e \bullet z) = (y \bullet e) \bullet z = y \bullet z = e$.

To prove~(\ref{lemm group properties 1}), let $x \in X$. There is $y \in X$ such that $x \bullet y = e$. Thus $e \bullet x = (x \bullet y) \bullet x = x \bullet (y \bullet x) = x$.

To prove~(\ref{lemm group properties 2}), assume that such a member $d$ exists. Then $d = e \bullet d$ by~(\ref{lemm group properties 1}) and $e \bullet d = e$ by assumption. Thus $d = e$.

To show the uniqueness in~(\ref{lemm group properties 3}), let $x, y, z \in X$ such that $x \bullet y = x \bullet z = e$. It follows that $y \bullet x = z \bullet x = e$. Therefore $y = y \bullet e = y \bullet (x \bullet z) = (y \bullet x) \bullet z = e \bullet z = z$.
\eproof

\section{Relations and maps}
\label{relations and maps}

In this Section we analyse how relations behave under maps. This is used subsequently for various purposes.

\midvspace

\blede
\label{lede gen relation prop}
Let $(X,R)$ and $(Y,S)$ be two relational spaces, and $f : X \longrightarrow Y$ a map. We use the same symbol for the function $f : X \!\times\! X \longrightarrow Y \!\times\! Y$, $f (x, z) = (f(x), f(z))$, as the two functions can be distinguished by their arguments. We have $f \left[ R \right] = \left\{ \, \left(f(x), f(z)\right) \, : \, (x, z) \in R \, \right\}$, which is a relation on~$Y$, and $f^{-1} \left[ S \right] = \left\{ \, (x, z) \in X \!\times\! X \, : \, \left(f(x), f(z)\right) \in S \, \right\}$, which is a relation on~$X$. The following statements hold:
\benum
\item If $S$ is transitive, then $f^{-1} \left[ S \right]$ is transitive.
\item If $S$ is reflexive, then $f^{-1} \left[ S \right]$ is reflexive.
\item If $S$ is antisymmetric and $f$ is injective, then $f^{-1} \left[ S \right]$ is antisymmetric.
\item If $S$ is antireflexive, then $f^{-1} \left[ S \right]$ is antireflexive.
\eenum

\elede

\bproof
Exercise.
\eproof

\bexam
\label{exam product pre-orderings}
Let $(X_i,R_i)$ ($i \in I$) be pre-ordered spaces, where $I$ is an index set, and $X = \bigtimes_{\!\! i \in I}\, X_i$. Then ${\mc R} = \left\{ p_i^{-1} \left[ R_i \right] \, : \, i \in I \right\}$ is a system of pre-orderings on~$X$.
\eexam

\brema
\label{rema generated relation}
Let $X$ be a set, $\mc R$ a system of relations on~$X$, and $S = \bigcap {\mc R}$. Then the following statements hold:
\benum
\item \label{rema generated relation 1} If every $R \in {\mc R}$ is transitive, then $S$ is transitive.
\item \label{rema generated relation 2} If every $R \in {\mc R}$ is reflexive, then $S$ is reflexive.
\item \label{rema generated relation 3} If there is $R \in {\mc R}$ such that $R$ is antisymmetric, then $S$ is antisymmetric.
\item \label{rema generated relation 4} If there is $R \in {\mc R}$ such that $R$ is antireflexive, then $S$ is antireflexive.
\eenum

In other words, if every $R \in {\mc R}$ is a pre-ordering, then $S$ is a pre-ordering. Moreover, if the members of ${\mc R}$ are pre-orderings and at least one member is an ordering, then $S$ is an ordering. Finally notice that the above also states conditions under which $S$ is an ordering in the sense of "$<$" or~"$\leq$".
\erema

\bexam
\label{exam product section ordering}
Let $(X_i,R_i)$ ($i \in I$) be pre-ordered spaces, where $I$ is an index set, and $X = \bigtimes_{\!\! i \in I}\, X_i$. Then $S = \bigcap \left\{ p_i^{-1} \left[ R_i \right] \, : \, i \in I \right\}$ is a pre-ordering on~$X$.
\eexam

\blede
\label{product directed set}
\index{Product directed space}
\index{Directed space!product}
Let $(X_i,R_i)$ ($i \in I$) be directed spaces, where $I$ is an index set, and $X = \bigtimes_{\!\! i \in I}\, X_i$. Then $S = \bigcap \left\{ p_i^{-1} \left[ R_i \right] \, : \, i \in I \right\}$ is a direction on~$X$. $(X,S)$ is called {\bf product directed space}.
\elede

\bproof
Exercise.
\eproof

The following notion is used in Section~\ref{interval topology} where we consider interval topologies.

\midvspace 

\bdefi
\index{Upwards independent}
\index{Relations!upwards independent}
\index{Downwards independent}
\index{Relations!downwards independent}
\index{Independent}
\index{Relations!independent}
Let $X$ be a set, ${\mc R} = \left\{ \prec_i \; : \; i \in I \right\}$ a system of pre-orderings on~$X$ where $I$ is an index set, and $S = \bigcap {\mc R}$. The pre-ordering $S$ is also denoted by~$\prec$. $\mc R$~is called {\bf upwards independent} if for every $i \in I$, every $s \in X$, and every $x \in X$ with $x \prec_i s$, there exists $y \in X$ such that $x \prec y$ and $y \prec_i s$.

$\mc R$~is called {\bf downwards independent} if for every $i \in I$, every $r \in X$, and every $x \in X$ with $r \prec_i x$, there exists $y \in X$ such that $y \prec x$ and $r \prec_i y$.

$\mc R$ is called {\bf independent} if it is both upwards and downwards independent.
\edefi

\blemm
\label{lemm independent induced rel}
Let $X$ be a set, ${\mc R} = \left\{ \prec_i \; : \; i \in I \right\}$ a system of pre-orderings on~$X$ where $I$ is an index set, and $S = \bigcap {\mc R}$. The pre-ordering $S$ is also denoted by~$\prec$. Intervals with respect to the pre-ordering~$\prec_i$ are denoted by subscript~$i$, those with respect to the pre-ordering~$\prec$ are denoted without subscript.
\benum
\item \label{lemm independent induced rel 1} If $\mc R$ is upwards independent, then we have for every $i \in I$ and $s \in X$ 
\[
\, \left] -\infty, s \right[_{\, i} \, = \, \bigcup \big\{ \, \left] -\infty, x \right[ \; : \; x \in X, \; x \prec_i s \big\}
\]

\item \label{lemm independent induced rel 2} If $\mc R$ is downwards independent, then we have for every $i \in I$ and $r \in X$ 
\[
\, \left] r, \infty \right[_{\, i} \, = \, \bigcup \big\{ \, \left] x, \infty \right[ \; : \; x \in X, \; r \prec_i x \big\}
\]

\item \label{lemm independent induced rel 3} If $\mc R$ is independent, then we have for every $i \in I$ and $r, s \in X$ 
\[
\, \left] r, s \right[_{\, i} \, = \, \bigcup \big\{ \, \left] x, y \right[ \; : \; x, y \in X, \; x \prec y, \; r \prec_i x, \; y \prec_i s \big\}
\]
\eenum

\elemm

\bproof
To see (\ref{lemm independent induced rel 1}), assume the stated condition and let $i \in I$ and $s \in X$. We have
\begin{eqnarray*}
\, \left] -\infty, s \right[_{\, i}
 \!\! & = & \!\! \big\{ z \in X \, : \, z \prec_i s \big\}\\
 & = & \!\! \big\{ z \in X \, : \, \exists x \in X \quad z \prec x, \; x \prec_i s \big\}\\
 & = & \!\! \bigcup \Big\{ \big\{ z \in X \, : \, z \prec x \big\} \, : \, x \in X, \; x \prec_i s \Big\}
\end{eqnarray*}

The proof of~(\ref{lemm independent induced rel 2}) is similar.

To show~(\ref{lemm independent induced rel 3}), assume the stated condition and let $i \in I$ and $r, s \in X$. We have
\begin{eqnarray*}
\, \left] r, s \right[_{\, i} \!\! & = & \! \left] r, \infty \right[_{\, i} \cap \, \left] -\infty, s \right[_{\, i}\\
 & = & \!\! \bigcup \big\{ \, \left] x, \infty \right[ \; : \; x \in X, \; r \prec_i x \big\} \; \cap \; \bigcup \big\{ \, \left] -\infty, y \right[ \; : \; y \in X, \; y \prec_i s \big\}\\
 & = & \!\! \bigcup \big\{ \, \left] x, y \right[ \; : \; x, y \in X, \; x \prec y, \; r \prec_i x, \; y \prec_i s \big\}
\end{eqnarray*}

where the second equation follows by~(\ref{lemm independent induced rel 1}) and~(\ref{lemm independent induced rel 2}), and the third equation by Lemma~\ref{lemm set equalities great}~(\ref{lemm set equalities great 1}).
\eproof

In the following Definition we introduce a notation that is convenient for the analysis of set functions in Section~\ref{set systems}.

\midvspace

\bdefi
\label{defi gen relation}
Given a relational space $(X,R)$ and a function $f : X \longrightarrow X$, the relation $f^{-1} \left[ R \right]$ is also denoted by~$R_f$. If $R$ is a pre-ordering, we also write $\, x \prec_f y \,$ for $(x,y) \in R_f$.
\edefi

The notation defined in the above Definition is meaningful since $R_f$ is a pre-ordering if $R$ is a pre-ordering.

\midvspace

\bdefi
\index{R-increasing}
Given a relational space $(X,R)$, a function $f : X \longrightarrow X$ is called $R$-{\bf increasing}, if $(x,y) \in R$ implies $\left(x,f(y)\right), \left(f(x),f(y)\right) \in R$ for every $x,y \in X$.
\edefi

\brema
Let $(X,R)$ be a relational space, and $f : X \longrightarrow X$ and $g : X \longrightarrow X$ two $R$-increasing maps. Then $g \circ f$ is $R$-increasing.
\erema

\blemm
\label{lemm increasing projection}
Given a set $X$, a reflexive pre-ordering $\prec$ on $X$, and an $\prec$-increasing projective map $f : X \longrightarrow X$, we have
\[
x \prec_f y \quad \Longleftrightarrow \quad x \prec f(y)
\]
\elemm

\bproof
Fix $x, y \in X$. We have $x \prec x$, and therefore $x \prec f(x)$. Assume $x \prec_f y$. It follows that $f(x) \prec f(y)$, and thus $x \prec f(y)$, since $f$ is transitive. Now assume instead that $x \prec f(y)$. Since $f$ is $\prec$-increasing and projective, we obtain $f(x) \prec f(y)$.
\eproof


\chapter{Numbers I}
\label{numbers i}
\setcounter{equation}{0}

\pagebreak

\section{Natural numbers, induction, recursion}
\label{natural numbers}

In Definition~\ref{defi naturalnumbers} we have defined the set $\naturalnumbers$ of natural numbers. In this Section we derive two important Theorems: the Induction principle for natural numbers and the Recursion theorem for natural numbers. Based on these Theorems we define and analyse the addition, multiplication, and exponentiation on the natural numbers. The natural numbers are the starting point for the construction of the other number systems below in this Chapter.

We first introduce the conventional symbols for four specific sets, that are natural numbers.

\midvspace

\bdefi
We define the sets $0 = \O$, $1 = \left\{ 0 \right\}$, $2 = \left\{ 0, 1 \right\}$, and $3 = \left\{ 0, 1, 2 \right\}$. Furthermore, we define the function $\sigma : \naturalnumbers \longrightarrow \naturalnumbers \setminus\! \left\{ 0 \right\}$, $\sigma (m) = m \cup \left\{ m \right\}$.
\edefi

We clearly have $\sigma(0) = 1$, $\sigma(1) = 2$, and $\sigma(2) = 3$. Notice that $\sigma$ is well-defined since $\naturalnumbers$ is inductive and $m \cup \left\{ m \right\}$ is non-empty for every $m \in \naturalnumbers$.

\midvspace

\btheo [Induction principle for natural numbers]
\label{lemm induction}
\index{Induction principle!natural numbers}
Let $A \subset \naturalnumbers$. If $0 \in A$ and if $\sigma (n) \in A$ for every $n \in A$, then $A = \naturalnumbers$.
\etheo

\bproof
Assume $A$ satisfies the conditions. Then $A$ is inductive. It follows that $\naturalnumbers \subset A$ by Definition~\ref{defi naturalnumbers}.
\eproof

\btheo
\label{theo naturals}
The natural numbers have the following properties:
\benum
\item \label{theo naturals 1} $\forall m \in \naturalnumbers \quad m \subset \naturalnumbers$
\item \label{theo naturals 2} $\forall n \in \naturalnumbers \quad \forall m \in n \quad m \in \naturalnumbers \, \wedge \, m \subset n$
\item \label{theo naturals 3} $\forall n \in \naturalnumbers \quad \forall m \in n \quad \sigma(m) \in \sigma(n)$
\item \label{theo naturals 4} $\forall m, n \in \naturalnumbers \quad m = n \, \vee \, m \in n \, \vee \, n \in m$
\item \label{theo naturals 5} $\forall m, n, p \in \naturalnumbers \quad m \in n \, \wedge \, n \in p \; \Longrightarrow \; m \in p$
\item \label{theo naturals 6} $\neg \exists m \in \naturalnumbers \quad m \notin m$
\item \label{theo naturals 7} $\forall m \in \naturalnumbers \quad \neg \exists n \in \naturalnumbers \quad m \in n \in m \cup \left\{ m \right\}$
\eenum

\etheo

\bproof
To see~(\ref{theo naturals 1}), let $A = \left\{ m \in \naturalnumbers \, : \, m \subset \naturalnumbers \right\}$. Clearly, $0 \in A$. Now assume that $m \in A$ for some $m \in \naturalnumbers$. Then we have $m \subset \naturalnumbers$. Therefore $\sigma(m) \subset \naturalnumbers$, and thus $\sigma(m) \in A$. It follows that $A = \naturalnumbers$ by the Induction principle.

To show~(\ref{theo naturals 2}), let $A = \left\{ n \in \naturalnumbers \, : \, \forall m \in n \; \left( m \in \naturalnumbers \wedge m \subset n \right) \right\}$. Clearly, $0 \in A$. Now assume that $n \in A$ for some $n \in \naturalnumbers$. Let $m \in \sigma(n)$. We have either $m \in n$ or $m = n$. In the first case, we obtain $m \in \naturalnumbers$ and $m \subset n \subset \sigma(n)$ by assumption. In the second case, we obviously have $m \in \naturalnumbers$ and $m \subset \sigma(n)$. We obtain $A = \naturalnumbers$ by the Induction principle.

To show~(\ref{theo naturals 3}) we again apply the Induction principle. The claim is trivially true for $n = 0$ and every $m \in n$. Assume the claim is true for some $n \in \naturalnumbers$ and every $m \in n$. Let $m \in \sigma(n)$. If $m \in n$, then $\sigma(m) \in \sigma(n)$ by assumption. If $m = n$, then $\sigma(m) = \sigma(n)$. It follows that $\sigma(m) \in \sigma(\sigma(n))$.

To prove~(\ref{theo naturals 4}), we use the Induction principle with respect to~$m$. First let $m = 0$. If $n = 0$, then we have $m = n$. If $n \neq 0$, then $0 \in n$, which is easily shown by the Induction principle. Thus the claim is true for $m = 0$. Assume the claim holds for some $m \in \naturalnumbers$ and every $n \in \naturalnumbers$. Fix $n \in \naturalnumbers$. If $n \in m$ or $n = m$, then $n \in \sigma(m)$. If $m \in n$, then either $\sigma(m) = n$ or $\sigma(m) \in n$ by~(\ref{theo naturals 3}).

To see~(\ref{theo naturals 5}), notice that $n \subset p$ under the stated conditions by~(\ref{theo naturals 2}).

(\ref{theo naturals 6}) is a consequence of Lemma~\ref{lemm element min}~(\ref{lemm element min 1}).

To see~(\ref{theo naturals 7}), notice that if such $n$ exists, then we have either $n = m$ or $n \in m$ both of which is excluded by Lemma~\ref{lemm element min}.
\eproof

\blede
\label{lede naturalnumbers ordering}
We define a total ordering in the sense of~"$<$" on the natural numbers by:
\[
m < n \quad \Longleftrightarrow \quad m \in n
\]

For every $m \in \naturalnumbers$, $\sigma(m)$ is the unique successor of~$m$, and we have $m = \left\{ n \in \naturalnumbers \, : \, n < m \right\}$.

We further define $\leq$ to be the total ordering in the sense of~"$\leq$" on the natural numbers by the method of Lemma~\ref{lemm both orderings}.
\elede

\bproof
This follows by Theorem~\ref{theo naturals}.
\eproof

Notice that, in particular, $1$ is the successor of~$0$~etc. 

\midvspace

\bdefi
\label{defi convention ordering natural numbers}
We adopt the convention that all notions related to orderings on~$\naturalnumbers$ refer to the ordering~$<$ as defined in Lemma and Definition~\ref{lede naturalnumbers ordering} unless otherwise specified.
\edefi

Note that in many contexts it is irrelevant whether the ordering~$<$ or the ordering~$\leq$ on~$\naturalnumbers$ is considered since most properties related to orderings are invariant, cf.\ Lemmas~\ref{lemm max min inv}, \ref{lemm bound sup inv}, and~\ref{lemm mon invar}.

The following version of the Induction principle allows us to prove statements inductively for all natural numbers that are larger than a fixed number.

\midvspace

\bcoro
\label{coro induction naturals version}
Let $A \subset \naturalnumbers$ and $m \in \naturalnumbers$. If $\sigma(m) \in A$ and if $\sigma(n) \in A$ for every $n \in A$ with $n > m$, then we have $\left\{ n \in \naturalnumbers \, : \, n > m \right\} \subset A$.
\ecoro

\bproof
We show that, under the stated conditions, $n > m$ implies $n \in A$ for every $n \in \naturalnumbers$ by the Induction principle. This implication is trivially true for $n = 0$. Now assume that it is true for some $n \in \naturalnumbers$. We distinguish the cases $n < m$, $n = m$, and $n >m$ by Theorem~\ref{theo naturals}~(\ref{theo naturals 4}). If $n < m$, then either $\sigma(n) < m$ or $\sigma(n) = m$ by Theorem~\ref{theo naturals}~(\ref{theo naturals 3}), and thus the implication again holds trivially for~$\sigma(n)$. If $n = m$, then $\sigma(n) \in A$ as this is amongst the conditions. If $n > m$, then $n \in A$ by assumption, and thus $\sigma(n) \in A$ since this is amongst the conditions.
\eproof

We also refer to Corollary~\ref{coro induction naturals version} as the Induction principle.

\midvspace

\btheo
\label{theo sigma bijective}
$\sigma : \naturalnumbers \longrightarrow \naturalnumbers \setminus\! \left\{ 0 \right\}$ is a bijection.
\etheo

\bproof
To see that $\sigma$ is surjective, notice that $1 \in \sigma \left[ \naturalnumbers \right]$ and $\sigma(m) \in \sigma \left[ \naturalnumbers \right]$ whenever $m \in \sigma \left[ \naturalnumbers \right]$. It follows that $\sigma \left[ \naturalnumbers \right] = \naturalnumbers \setminus\! \left\{ 0 \right\}$ by the Induction principle in the form of Corollary~\ref{coro induction naturals version}.

To show that $\sigma$ is injective, let $m, n \in \naturalnumbers$ such that $\sigma(m) = \sigma(n)$. Hence we have $m \cup \left\{ m \right\} = n \cup \left\{ n \right\}$. This implies
\[
m = n \vee \big( m \in n \wedge n \in m \big)
\]

By Theorem~\ref{theo naturals} it follows that $m = n$.
\eproof

\bcoro
\label{lemm naturals predecessor}
Every $m \in \naturalnumbers \setminus\! \left\{ 0 \right\}$ has a unique predecessor with respect to the ordering~$<$.
\ecoro

\bproof
This follows by Lemma and Definition~\ref{lede naturalnumbers ordering} and Theorem~\ref{theo sigma bijective}.
\eproof

The next definition is a modification of Definition~\ref{defi function set} for the case where the superscript is a natural number larger than~$0$.

\midvspace

\bdefi
\label{defi exp nat number}
Let $X$ be a set and $m \in \naturalnumbers$, $m > 0$. The system of functions $X^I$ where $I = \sigma(m) \!\setminus\! \left\{ 0 \right\}$ is also denoted by~$X^m$.
\edefi

Notice that this deviates from Definition~\ref{defi function set} where the superscript would be the domain and thus would contain $0$ but not~$m$. By Remark~\ref{rema square set} we may write members of~$X^2$ as ordered pairs as follows.

\midvspace

\bdefi
Given a set $X$, we also write $\left( f(1), f(2) \right)$ for $f \in X^2$.
\edefi

\bdefi
Let $m, n \in \naturalnumbers$ with $m < n$. Further let $I = \sigma(n) \!\setminus\! m$ and $A_i$ ($i \in I$) be sets. We define
\begin{center}
$\bigcup_{i = m}^n A_i \, = \, \bigcup_{i \in I} A_i \, , \quad \quad
\bigcap_{i = m}^n A_i \, = \, \bigcap_{i \in I} A_i \, , \quad \quad
\bigtimes_{\!\! i = m}^{\!\! n} \, A_i \, = \, \bigtimes_{\!\! i \in I} \, A_i$
\end{center}

\edefi

\bdefi
Let $m \in \naturalnumbers$, $I = \naturalnumbers \!\setminus\! m$, and $A_i$ ($i \in I$) be sets. We define
\begin{center}
$\bigcup_{i = m}^{\infty} A_i \, = \, \bigcup_{i \in I} A_i \, , \quad \quad
\bigcap_{i = m}^{\infty} A_i \, = \, \bigcap_{i \in I} A_i \, , \quad \quad
\bigtimes_{\!\! i = m}^{\!\! \infty} \, A_i \, = \, \bigtimes_{\!\! i \in I} \, A_i$
\end{center}

\edefi

The following Theorem states that one may define a function from $\naturalnumbers$ to a set~$X$ recursively.

\midvspace

\btheo[Recursion for natural numbers]
\label{theo recursive def}
Given a set $X$, a point $x \in X$ and a function $f : X \longrightarrow X$, there exists a unique function $g : \naturalnumbers \longrightarrow X$ with the following properties:
\benum
\item \label{theo recursive def 1} $g(0) = x$
\item \label{theo recursive def 2} $g \left( \sigma(n) \right) = f \left( g(n) \right)$ \quad for every $n \in \naturalnumbers$
\eenum 

\etheo

\bproof
For every $p \in \naturalnumbers$ there exists a map $G : \sigma (p) \longrightarrow X$ with the following properties:
\benum
\item $G(0) = x$
\item $G \left( \sigma (n) \right) = f \left( G(n) \right)$ \quad for every $n \in p$
\eenum\
\hspace{0.05\textwidth}
\parbox{0.95\textwidth}
{[This is clear for $p = 0$. Assume there exists such a function $G$ for $p \in \naturalnumbers$. We define $H : \sigma (\sigma (p)) \longrightarrow X$, $H|\sigma(p) = G$, $H(\sigma (p)) = f(G(p))$. The assertion follows by the Induction principle.]}
\\

We call a function $G$ with these properties a "debut of size~$\sigma (p)$". Let $G$ and $H$ be two debuts of sizes $\sigma(p)$ and $\sigma(q)$, respectively, where $p, q \in \naturalnumbers$. We may assume that $p < q$ or $p = q$. We have $G = H|\sigma(p)$.
\\

\hspace{0.05\textwidth}
\parbox{0.95\textwidth}
{[Clearly, $G(0) = x = H(0)$. Now let $n \in \naturalnumbers$, $n < p$ and assume that $G(n) = H(n)$. Then also $G(\sigma(n)) = H(\sigma(n))$ holds. The assertion follows by the Induction principle.]}
\\

Now, for every $n \in \naturalnumbers$, let $g(n) = G(n)$ where $G$ is the debut of size~$\sigma (n)$. Clearly, $g$ satisfies (\ref{theo recursive def 1}) and~(\ref{theo recursive def 2}) of the claim.
\\

\hspace{0.05\textwidth}
\parbox{0.95\textwidth}
{[(\ref{theo recursive def 1}) is satisfied since $g(0) = x$. (\ref{theo recursive def 2})~is satisfied for $n = 0$ since $g(1) = G(1) = f(G(0)) = f(x) = f(g(0))$ where $G$ is the debut of size~$2$. Now assume that~(\ref{theo recursive def 2}) is true for some $n \in \naturalnumbers$, that is $g(\sigma(n)) = f(g(n))$. Then we have $g \big( \sigma(\sigma(n)) \big) = G \big( \sigma(\sigma(n)) \big) = f \big( G(\sigma(n)) \big) = f \big( H(\sigma(n)) \big) = f \big( g(\sigma(n)) \big)$ where $G$ is the debut of size~$\sigma \big( \sigma(\sigma(n)) \big)$ and $H$ the debut of size~$\sigma ( \sigma (n) )$. Thus (\ref{theo recursive def 2}) is true for every $n \in \naturalnumbers$ by the Induction principle.]}
\\

To see that $g$ is unique, assume that also $h : \naturalnumbers \longrightarrow X$ satisfies (\ref{theo recursive def 1}) and~(\ref{theo recursive def 2}) of the claim. Obviously, $g(0) = x = h(0)$. Moreover, $g(n) = h(n)$ for some $n \in \naturalnumbers$ implies $g(\sigma(n)) = h(\sigma(n))$. It follows by the Induction principle that $g = h$.
\eproof

Theorem~\ref{theo recursive def} can be used to define powers (i.e.\ multiple products) of a relation on a given set~$X$ as follows.

\midvspace

\blede
\label{def rel npower}
Given a relational space $(X,R)$, we define, for every $m \in \mathbb{N}$, $0 < m$, a relation $R^m$ on $X$ by
\benum
\item $R^1 = R$
\item $R^{\sigma(m)} = R^m R$
\eenum

\elede

\bproof
The existence and uniqueness of $R^m$ for every $m \in \mathbb{N}$, $0 < m$ follows by Theorem~\ref{theo recursive def}.
\eproof

Notice that this definition deviates from Definition~\ref{defi function set}. Since Lemma and Definition~\ref{def rel npower} is valid only for relations, there is generally no risk of confusion.

As a further immediate consequence of Theorem~\ref{theo recursive def} we obtain the following result that allows us to define binary functions on~$\naturalnumbers$ recursively.

\midvspace

\bcoro
\label{coro recursive operation}
Given two maps $e : \naturalnumbers \longrightarrow \naturalnumbers$ and $f : \naturalnumbers \times \naturalnumbers \longrightarrow \naturalnumbers$, there exists a unique function $g : \naturalnumbers \times \naturalnumbers \longrightarrow \naturalnumbers$ with the following properties:
\benum
\item $g(m, 0) = e(m)$ \quad for every $m \in \naturalnumbers$
\item $g\left(m, \sigma (n)\right) = f\left(m, g(m, n)\right)$ \quad for every $m, n \in \naturalnumbers$
\eenum 

\ecoro

\bproof
For each $m \in \naturalnumbers$ there is a unique function $g_m : \naturalnumbers \longrightarrow \naturalnumbers$ by Theorem~\ref{theo recursive def} with the following properties:
\benum
\item $g_m(0) = e(m)$
\item $g_m\left(\sigma(n)\right) = f\left(m, g_m(n)\right)$ \quad for every $n \in \naturalnumbers$
\eenum 

We may then define $g$ by $g(m,n) = g_m(n)$.
\eproof 

In the remainder of this Section we introduce three binary functions on~$\naturalnumbers$, viz.\ addition, multiplication, and exponentiation.

\midvspace

\blede
\index{Addition!natural numbers}
\index{Sum!natural numbers}
There is a unique binary function $+$ on $\naturalnumbers$ such that for every $m, n \in \naturalnumbers$ we have

\benum
\item \label{natural addition zero} $m + 0 = m$
\item \label{natural addition greater} $m + \sigma(n) = \sigma(m + n)$
\eenum

The function $+$ is called {\bf addition on}~$\naturalnumbers$. For every $m, n \in \naturalnumbers$, the expression $m + n$ is called the {\bf sum of $m$ and~$n$}.

$+$~is associative and commutative. Furthermore, for every $m, n, p \in \naturalnumbers$, the following implication holds:
\[
m < n \quad \Longrightarrow \quad m + p < n + p
\]

\elede

\bproof
The existence and uniqueness of the function~$+$ follows by Corollary~\ref{coro recursive operation}.

$+$~is commutative since it follows by the Induction principle that for every $p \in \naturalnumbers$ the following equations hold for $m, n \in \naturalnumbers$, $m + n = p$:

\benum
\item \label{equ comm 1} $m + n = n + m$
\item \label{equ comm 2} $\sigma(n) + m = n + \sigma(m)$
\eenum

\hspace{0.05\textwidth}
\parbox{0.95\textwidth}
{[First notice that $m + n = 0$ implies $m = n = 0$. Thus (\ref{equ comm 1}) holds for $m, n \in \naturalnumbers$, $m + n = 0$. Moreover, we have $0 + m = m$ for every $m \in \naturalnumbers$. In fact, this equation clearly holds for $m = 0$ and, assuming it holds for some $m \in \naturalnumbers$, it also holds for $\sigma(m)$ since $0 + \sigma(m) = \sigma(0 + m) = \sigma(m)$. As a special case, we obtain (\ref{equ comm 2}) for $m, n \in \naturalnumbers$, $m + n = 0$, viz.\ $1 + 0 = 0 + 1$. Now assume that (\ref{equ comm 1}) and (\ref{equ comm 2}) hold for some $p \in \naturalnumbers$ and for every $m, n \in \naturalnumbers$ with $m + n = p$. Let $m, n \in \naturalnumbers$ such that $m + n = \sigma(p)$. If $n = 0$, then (\ref{equ comm 1}) holds as shown above. If $0 < n$, then let $q$ be the predecessor of~$n$. We then have $m + n = m + \sigma(q) = \sigma(m + q) = \sigma(q + m) = q + \sigma(m) = \sigma(q) + m = n + m$. Hence (\ref{equ comm 1}) also holds for $0 < n$. If $m = 0$, then (\ref{equ comm 2}) clearly holds. If $0 < m$, then let $q$ be the predecessor of~$m$. We then obtain $\sigma(n) + m = \sigma(n) + \sigma(q) = \sigma(\sigma(n) + q) = \sigma(n + \sigma(q)) = \sigma(n + m) = n + \sigma(m)$, which is equation~(\ref{equ comm 2}).]}
\\

To see that the addition is associative, let $m, p \in \naturalnumbers$ and
\[
A = \big\{ n \in \naturalnumbers \, : \, (m + n) + p = m + (n + p) \big\}
\]

Clearly, $0 \in A$ by commutativity. Now assume that $n \in A$ for some $n \in \naturalnumbers$ . We have $(m + \sigma(n)) + p = \sigma(m+n) + p = p + \sigma(m+n) = \sigma(p+(m+n)) = \sigma((m+n)+p) = \sigma(m+(n+p)) = m + \sigma(n+p) = m + \sigma(p+n) = m + (p + \sigma(n))$.

Also the last claim can be shown by the Induction principle. It clearly holds for $m, n \in \naturalnumbers$ and $p = 0$. Assume it holds for some $p \in \naturalnumbers$ and every $m, n \in \naturalnumbers$. Let $m, n \in \naturalnumbers$ with $m < n$. Then $m + \sigma(p) = \sigma(m + p) < \sigma(n + p) = n + \sigma(p)$.
\eproof

The associativity of the addition allows us to write multiple sums without brackets, i.e.\ $m + n + p$ instead of $(m + n) + p$ or $m + (n + p)$, and similarly for sums with more than three terms.

\midvspace

\blemm
\label{lemm natural diff}
For every $m, n \in \naturalnumbers$ with $m < n$ there is $p \in \naturalnumbers$ such that $m + p = n$.
\elemm

\bproof
Let $m \in \naturalnumbers$. Then the claim follows by the Induction principle in the form of Corollary~\ref{coro induction naturals version} as follows. The claim clearly holds for $n = \sigma(m)$. Assuming it holds for some $n \in \naturalnumbers$ with $n > m$, we may choose $p \in \naturalnumbers$ such that $m + p = n$. It follows that $m + \sigma(p) = \sigma(n)$.
\eproof

\blemm
\label{proposition relation power}
Let $(X,R)$ be a relational space, and $m, n \in \mathbb{N} \setminus\! \left\{ 0 \right\}$. Then the following statements hold:
\benum
\item \label{proposition relation power 1} $R^{m + n} = R^m R^n = R^n R^m$
\item \label{proposition relation power 2} $\Delta \subset R \; \Longrightarrow \; R^m \subset R^{m + n}$
\eenum

\elemm

\bproof
The first equation in~(\ref{proposition relation power 1}) clearly holds for $m \in \naturalnumbers \setminus\! \left\{ 0 \right\}$ and $n = 1$. Assume it holds for every $m \in \naturalnumbers \setminus\! \left\{ 0 \right\}$ and some $n \in \naturalnumbers \setminus\! \left\{ 0 \right\}$. Then we obtain, for every $m \in \naturalnumbers \setminus\! \left\{ 0 \right\}$, $R^{m + \sigma(n)} = R^{m + n} R = \left( R^m R^n \right) R = R^m \left( R^n R \right) = R^m R^{\sigma(n)}$. Thus the first equation is true for every $m, n \in \naturalnumbers \setminus\! \left\{ 0 \right\}$ by the Induction principle, Corollary~\ref{coro induction naturals version}. The second equation is a consequence of the first one and the commutativity of addition.

Also (\ref{proposition relation power 2}) can be shown by the Induction principle. The claim clearly holds for $n = 1$ and every $m \in \naturalnumbers$. Assume it holds for some $n \in \naturalnumbers$ and every $m \in \naturalnumbers$. We have, for every $m \in \naturalnumbers$, $R^m \subset R^{m + n} \subset R^{(m+n)+1} = R^{m + \sigma(n)}$.
\eproof

Using arbitrary powers of relations one can construct a pre-ordered space from an arbitrary relational space such that the set is the same and the pre-ordering contains the original relation. The following result also says that the constructed pre-ordering is minimal.

\midvspace

\blemm
Let $(X,R)$ be a relational space. Then $S = \bigcup \left\{ R^n \, : \, n \in \mathbb{N},\, n > 0 \right\}$ is a pre-ordering on $X$. If $T$ is a pre-ordering on $X$ with $R \subset T$, then $S \subset T$.
\elemm

\bproof
Let $(x,y), (y,z) \in S$. Then there exist $m, n \in \mathbb{N} \setminus\! \left\{ 0 \right\}$ such that $(x,y) \in R^m$ and $(y,z) \in R^n$. Hence $(x,z) \in R^{m+n} \subset S$.

Now let $T$ be a pre-ordering on $X$ with $R \subset T$. Then $R^m \subset T$ for every $m \in \naturalnumbers \setminus\! \left\{ 0 \right\}$ by the Induction principle. It follows that $S \subset T$.
\eproof

\blede
\index{Multiplication!natural numbers}
\index{Product!natural numbers}
There is a unique binary function $\cdot$ on $\naturalnumbers$ such that for every $m, n \in \naturalnumbers$ we have

\benum
\item $m \cdot 0 = 0$
\item $m \cdot \sigma(n) = (m \cdot n) + m$
\eenum

The function $\cdot$ is called {\bf multiplication on}~$\naturalnumbers$. For every $m, n \in \naturalnumbers$, the expression $m \cdot n$ is called the {\bf product of $m$ and~$n$}, also written $m \, n$. The multiplication on~$\naturalnumbers$ is commutative and associative. Furthermore, for every $m, n, p \in \naturalnumbers$, $0 < p$, we have the distributive law
\[
(m + n) \cdot p \, = \, m \cdot p + n \cdot p
\]

and the following implication holds:
\[
m < n \quad \Longrightarrow \quad m \cdot p < n \cdot p
\]

We define that, in the absence of brackets, products are evaluated before sums. Thus we may write $m \cdot n + p$ instead of $(m \cdot n) + p$, and $m + n \cdot p$ instead of $m + (n \cdot p)$ without ambiguity.
\elede

\bproof
The existence and uniqueness of the function follows by Corollary~\ref{coro recursive operation}.

In order to prove that the multiplication is commutative and that the distributive law holds, first notice that, for every $m, n \in \naturalnumbers$, we have $\sigma(m) \, n = m \, n + n$.
\\

\hspace{0.05\textwidth}
\parbox{0.95\textwidth}
{[For every $m \in \naturalnumbers$, we have $\sigma(m) \cdot 0 = 0 = m \cdot 0 + 0$. Assuming the claim is true for some $n \in \naturalnumbers$ and every $m \in \naturalnumbers$, we have $\sigma(m) \, \sigma(n) = \sigma(m) \, n + \sigma(m) = m \, n + n + m + 1 = m \, n + m + \sigma(n) = m \, \sigma(n) + \sigma(n)$.]}
\\

It follows that $m \, n = n \, m$ for every $m, n \in \naturalnumbers$.
\\

\hspace{0.05\textwidth}
\parbox{0.95\textwidth}
{[We show that, for every $p \in \naturalnumbers$, the equation $m \, n = n \, m$ holds for every $m, n \in \naturalnumbers$ with $m + n = p$. Clearly, this is true for $p = 0$ because this implies $m = n = 0$. Assuming it is true for $p \in \naturalnumbers$, let $m, n \in \naturalnumbers$ with $m + n = \sigma(p)$. We may assume $0 < m$. Let $q$ be the predecessor of $m$, that is $\sigma(q) = m$. Then we have $q + n = p$. It follows that $m \, n = \sigma(q) \, n = q \, n + n = n \, q + n = n \, \sigma(q) = n \, m$.]}
\\

It also follows that $(m + n) \, p = m \, p + n \, p$ for every $m, n, p \in \naturalnumbers$.
\\

\hspace{0.05\textwidth}
\parbox{0.95\textwidth}
{[For every $m, p \in \naturalnumbers$ this equation is clearly satisfied for $n = 0$. Now assume this equation holds for some $n \in \naturalnumbers$ and for every $m, p \in \naturalnumbers$. We then have $\left( m + \sigma(n) \right) \, p = \sigma(m + n) \, p = (m + n) \, p + p = m \, p + n \, p + p = m \, p + \sigma(n) \, p$.]}
\\

We next show that the multiplication is associative. We have $(m \cdot 0) \cdot p = 0 = m \cdot (0 \cdot p)$ for every $m, p \in \naturalnumbers$. Assume that $(m \, n) \, p = m \, (n \, p)$ for some $n \in \naturalnumbers$ and every $m, p \in \naturalnumbers$. Then $(m \, \sigma(n)) \, p = (m \, n + m) \, p = (m \, n) \, p + m \, p = m \, (n \, p) + m \, p = m \, (n \, p + p) = m \, (\sigma(n) \, p)$.

To see the last claim, notice that it clearly holds for $p = 1$ and every $m, n \in \naturalnumbers$. Assume it holds for some $p \in \naturalnumbers$, $0 < p$, and every $m, n \in \naturalnumbers$. Let $m, n \in \naturalnumbers$ and $m < n$. Then $m \, \sigma(p) = m \, p + m < n \, p + n = n \, \sigma(p)$.
\eproof

As in the case of addition, we may write multiple products without brackets because of the associativity of the multiplication, i.e.\ we may write $m \, n \, p$ instead of $(m \, n) \, p$ or $m \, (n \, p)$. Clearly, also the distributive law $p \, (m + n) = p \, m + p \, n$ holds, since addition and multiplication are commutative.

\midvspace

\bdefi
\index{Even}
\index{Odd}
Let $m \in \naturalnumbers$. $m$~is called {\bf even} if there is $n \in \naturalnumbers$ such that $m = 2n$, otherwise it is called {\bf odd}.
\edefi

The following Lemma and Corollary are immediate consequences of this Definition.

\midvspace

\blemm
Let $m \in \naturalnumbers$. If $m$ is odd, then there is $n \in \naturalnumbers$ such that $m = 2n + 1$.
\elemm

\bproof
We prove that for every $m \in \naturalnumbers$ there is $n \in \naturalnumbers$ such that $m = 2n$ or $m = 2n + 1$. This is clearly true for $m = 0$. Assuming it is true for some $m \in \naturalnumbers$, we have either $m + 1 = 2n + 1$ or $m + 1 = 2n + 1 + 1 = 2(n+1)$.
\eproof

\bcoro
Let $m, n, p, q \in \naturalnumbers$ where $m$ and $n$ are even, and $p$ and $q$ are odd. Then $m+n$ and $p+q$ are even, and $m+p$ is odd.
\ecoro

\bproof
We may choose $m_0, n_0, p_0, q_0 \in \naturalnumbers$ such that $m = 2m_0$, $n = 2n_0$, $p = 2p_0 + 1$, and $q = 2q_0 + 1$. It follows that $m + n = 2(m_0 + n_0)$, $p + q = 2(p_0 + q_0 + 1)$, and $m + p = 2(m_0 + p_0) + 1$.

\eproof

The following proposition is applied in Section~\ref{cardinality}.

\midvspace

\bprop
\label{prop gauss integer}
For every $m \in \naturalnumbers$ there is $n \in \naturalnumbers$ such that $2 n = m \, \sigma(m)$.
\eprop

\bproof
The claim clearly holds for $m = 0$. Assume that it holds for some $m \in \naturalnumbers$. Let $n \in \naturalnumbers$ such that $2 n = m \, \sigma(m)$. Then we have $\sigma(m) \, \sigma(\sigma(m)) = \sigma(m) \, \sigma(m) + \sigma(m) = \sigma(m) \, m + \sigma(m) + \sigma(m) = 2 \left( n + \sigma(m) \right)$.
\eproof

\blede
\label{lede exponentiation natural numbers}
\index{Exponentiation!natural numbers}
There is a unique binary function $h$ on $\naturalnumbers$ such that for every $m, n \in \naturalnumbers$ we have

\benum
\item $h(m, 0) = 1$
\item $h(m, \sigma(n)) = h(m, n) \cdot m$
\eenum

This function is called {\bf exponentiation on}~$\naturalnumbers$. We also write $m^n$ for~$h(m, n)$. We define that, in the absence of brackets, exponentiation is evaluated before sums or products. Thus, for $m, n, p \in \naturalnumbers$, we may write $m^n + p$ instead of $\left( m^n \right) + p$, and $m^n \cdot p$ instead of $\left( m^n \right) \cdot p$. We have, for every $m, n, p \in \naturalnumbers$
\[
m^{n+p} = m^n \, m^p, \quad \quad \left( m^n \right)^p = m^{n \, p}, \quad \quad \left(m \, n\right)^{p} = m^p \, n^p
\]

and the implication
\[
\left( m < n \right) \, \wedge \, \left( 1 < p \right) \quad \Longrightarrow \quad \left( m^p < n^p \right) \, \wedge \, \left( p^m < p^n \right)
\]
\elede

\bproof
The existence and uniqueness of the function follows from Corollary~\ref{coro recursive operation}.

Further notice that the three equations clearly hold for $p = 0$ and every $m, n \in \naturalnumbers$. Now assume they hold for some $p \in \naturalnumbers$ and every $m, n \in \naturalnumbers$. We then have $m^{n + \sigma(p)} = m^{n + p} \, m = m^n \, m^p \, m = m^n \, m^{\sigma(p)}$, which proves the first equation. To show the second equation, notice that $\big( m^n \big)^{\sigma(p)} = \big( m^n \big)^p \, m^n = m^{n p} \, m^n = m^{n p + n} = m^{n \, \sigma(p)}$. Finally, we have $\left(m \, n\right)^{\sigma(p)} = \left(m \, n\right)^p \, m \, n = m^p \, n^p \, m \, n = m^{\sigma(p)} \, n^{\sigma(p)}$, which proves the third equation.

We now show that $\left( m < n \right) \wedge \left( 1 < p \right)$ implies $m^p < n^p$ for every $m, n, p \in \naturalnumbers$. Let $m, n \in \naturalnumbers$. The implication clearly holds for $p = 2$. Further, assuming it holds for some $p > 1$, it also holds for $\sigma(p)$.

Finally we show that $m < n$ implies $p^m < p^n$ for every $m, n, p \in \naturalnumbers$ where $p > 1$. Let $m, p \in \naturalnumbers$ with $p > 1$. The claim is clearly true for $n = \sigma(m)$. Assuming it holds for some $n \in \naturalnumbers$ with $n \geq \sigma(m)$, we have $p^{\sigma(n)} = p^n \, p > p^m \, p > p^m$.
\eproof

Notice that exponentiation is neither associative nor commutative.

\section{Ordinal numbers}
\label{ordinal numbers}

In this Section we introduce the notion of ordinal numbers. Ordinals are a natural extension of the natural numbers. In this text they are used to analyse the Choice axiom in Section~\ref{choice} and as a basis of the concept of cardinality in Section~\ref{cardinality}, which corresponds, in a sense, to the "number of elements" of a set.

In order to define the ordinal numbers we begin with the following definition following~\cite{Kelley}.

\midvspace

\bdefi
\index{Full}
\index{Set!full}
A set $X$ is called {\bf full} if $A \in X$ implies $A \subset X$.
\edefi

In the literature full sets are also called transitive (cf.~\cite{Jech}). However we do not adopt this notion here because we use it as a property that relations may have.

We now define a relation that corresponds to the property of a set to be member of another, but is however restricted to a specified set. Only when restricted to a given set it is a relation according to our Definition~\ref{defi rel inv prod}.

\midvspace

\bdefi
\index{Element relation}
\index{Relation!element}
Given a set $X$, the relation $R$ on $X$ defined by
\[
(x,y) \in R \quad \Longleftrightarrow \quad x \in y
\]

is called {\bf element relation on}~$X$. We also write $x \in_X y$ for $(x,y) \in R$.
\edefi

\brema
Given a set $X$, the element relation $\in_X$ is irreflexive. If $Y \subset X$, the restriction of $\in_X$ to $Y$ is $\in_Y$.
\erema

\bdefi
\label{defi ordinal}
\index{Ordinal number}
\index{Numbers!ordinal}
\index{Ordinal}
A non-empty set $X$ is called {\bf ordinal number} or {\bf ordinal} if $X$ is full and $\in_X$ has the minimum property. We also define that $\O$ is an ordinal.
\edefi

We mainly use Greek letters for ordinals---as we do for positive reals. However, we always say explicitly when a variable is meant to be an ordinal. This property is never automatically implied by the mere usage of a Greek letter.

\midvspace

\blemm
Every natural number is an ordinal.
\elemm

\bproof
Let $m \in \naturalnumbers$. If $m = 0$, then $m$ is an ordinal by definition. If $0 < m$, it follows by Theorem~\ref{theo naturals}~(\ref{theo naturals 2}) that $m$ is full. We show by the Induction principle that $m$ has the minimum property for every $m \in \naturalnumbers$, $0 < m$. $\in_m$ clearly has the minimum property if $m = 1$. Now assume that $\in_m$, where $m \in \naturalnumbers$ with $m > 0$, has the minimum property, and let $A \subset \sigma(m)$ with $A \neq \O$. If $A \cap m \neq \O$, the minimum of $A$ is the same as the minimum of $A \!\setminus\! \left\{ m \right\}$. If $A \cap m = \O$, then $m$ is the only member and hence the minimum of~$A$.
\eproof

\brema
\label{rema member ordinal}
Let $\alpha$ be a non-empty ordinal and $\beta \in \alpha$. Then $\beta \subset \alpha$ and, if $\beta \neq \O$, then $\in_{\beta}$ has the minimum property.
\erema

\blemm
\label{lemm ordinal well-ordering}
If $\alpha$ is a non-empty ordinal, then $\in_{\alpha}$ is a well-ordering in the sense of~"$<$".
\elemm

\bproof
Notice that $\in_{\alpha}$ is connective by Remark~\ref{prop minimum property connective} since it has the minimum property.

It remains to show that $\in_{\alpha}$ is transitive. Let $\beta, \gamma, \delta \in \alpha$ with $\beta \in \gamma$ and $\gamma \in \delta$. Since $\in_{\alpha}$ is connective, we have $\beta = \delta$, $\delta \in \beta$, or $\beta \in \delta$. The first two cases are excluded by Lemma~\ref{lemm element min}.
\eproof

\blemm
\label{lemm ordinals properties}
The following statements hold for ordinals:

\benum
\item \label{lemm ordinals properties 1} If $\alpha$ is an ordinal and $\beta \in \alpha$, then $\beta$ is an ordinal.
\item \label{lemm ordinals properties 2} If $\alpha$ and $\beta$ are ordinals with $\alpha \subset \beta$ and $\alpha \neq \beta$, then $\alpha \in \beta$.
\item \label{lemm ordinals properties 3} If $\alpha$ and $\beta$ are ordinals, then $\alpha \subset \beta$ or $\beta \subset \alpha$.
\item \label{lemm ordinals properties 4} If $\alpha$ is an ordinal, then $\alpha \cup \left\{ \alpha \right\}$ is an ordinal.
\item \label{lemm ordinals properties 5} If $\alpha$ is an ordinal, then there is no ordinal $\beta$ such that $\alpha \in \beta \in \alpha \cup \left\{ \alpha \right\}$.
\eenum

\elemm

\bproof
To see~(\ref{lemm ordinals properties 1}), assume the stated conditions. If $\beta = \O$, then $\beta$ is an ordinal. Now assume that $\beta \neq \O$. Then $\in_{\beta}$ has the minimum property by Remark~\ref{rema member ordinal}. To see that $\beta$ is full, let $\delta \in \gamma \in \beta$. Since $\alpha$ is full, we have $\gamma, \delta \in \alpha$. It follows that $\delta \in \beta$ because $\in_{\alpha}$ is transitive by Lemma~\ref{lemm ordinal well-ordering}. Thus $\gamma \subset \beta$.

To prove~(\ref{lemm ordinals properties 2}) let, under the stated conditions, $\gamma$ be the minimum of $\beta \!\setminus\! \alpha$. We clearly have $\gamma \subset \alpha$. This shows the claim for $\alpha = \O$. Now assume that $\alpha \neq \O$. Let $\delta \in \alpha$. Since $\in_{\beta}$ is connective, we have either $\gamma \in \delta$ or $\delta \in \gamma$. The case $\gamma \in \delta$ is excluded, because this implies $\gamma \in \alpha$, since $\alpha$ full. It follows that $\alpha \subset \gamma$, and thus $\alpha = \gamma \in \beta$.

To see~(\ref{lemm ordinals properties 3}), notice that for ordinals $\alpha$ and $\beta$, $\alpha \cap \beta$ clearly is an ordinal, say~$\gamma$. By~(\ref{lemm ordinals properties 2}) it follows that $\gamma \in \alpha$ or $\gamma = \alpha$. Similarly, it follows that $\gamma \in \beta$ or $\gamma = \beta$. It is not possible that $\gamma \in \alpha \cap \beta = \gamma$. Hence we have either $\gamma = \alpha$ or $\gamma = \beta$.

(\ref{lemm ordinals properties 4}) follows by Definition~\ref{defi ordinal}.

To see~(\ref{lemm ordinals properties 5}), assume there are ordinals $\alpha$ and $\beta$ such that $\alpha \in \beta \in \alpha \cup \left\{ \alpha \right\}$. It follows that $\alpha \subset \beta \subset \alpha \cup \left\{ \alpha \right\}$, and thus $\beta = \alpha$ or $\beta = \alpha \cup \left\{ \alpha \right\}$, which is a contradiction.
\eproof

Note that~(\ref{lemm ordinals properties 2}) implies that every ordinal that is not the empty set contains the empty set.

\midvspace

\blemm
\label{lemm ordinals total ord}
Let $A$ be a non-empty set of ordinals. Then $\in_A$ is a total ordering in the sense of~"$<$".
\elemm

\bproof
Let $\alpha, \beta, \gamma \in A$. If $\alpha \in \beta \in \gamma$, then $\beta \subset \gamma$, and therefore $\alpha \in \gamma$. Thus $\in_A$ is transitive.

Moreover $\in_A$ is connective by Lemma~\ref{lemm ordinals properties} (\ref{lemm ordinals properties 2}) and~(\ref{lemm ordinals properties 3}).
\eproof

\blemm
\label{lemm ordinals minimum}
Let $A$ be a non-empty set of ordinals. $A$ has a minimum with respect to the ordering~$\in_A$. $\in_A$~is a well-ordering.
\elemm

\bproof
We may choose $\alpha \in A$ such that $\alpha \cap A = \O$ by the Regularity axiom. For every $\beta \in A \!\setminus\! \left\{ \alpha \right\}$, we have either $\alpha \in \beta$ or $\beta \in \alpha$ since $\in_A$ is connective by Lemma~\ref{lemm ordinals total ord}. The latter is a contradiction. Therefore $\alpha$ is a minimum of~$A$.

The second claim is a consequence of the first one.
\eproof

\blemm
\label{lemm ordinals union intersection}
Let $A$ be a non-empty set of ordinals. Then $\bigcap A$ and $\bigcup A$ are ordinals too. In particular, $\naturalnumbers$ is an ordinal.
\elemm

\bproof
Exercise.
\eproof

\brema
\label{rema well-ordering naturals}
Notice that the ordering $<$ on $\naturalnumbers$ as defined in Lemma and Definition~\ref{lede naturalnumbers ordering} is identical to the well-ordering $\in_{\naturalnumbers}$ on $\naturalnumbers$ when $\naturalnumbers$ is considered as ordinal.
\erema

\brema
\label{rema ordinal vs natural}
Let $\alpha$ be an ordinal such that $\alpha \notin \naturalnumbers$ and $\alpha \neq \naturalnumbers$. Then we have $m \in \alpha$ for every $m \in \naturalnumbers$ and $\naturalnumbers \in \alpha$.
\erema

\blemm
\label{lemm no set ordinals}
There is no set that contains every ordinal number.
\elemm

\bproof
Assume $A$ is such a set. We define $\alpha = \bigcup A$ and $\beta = \alpha \cup \left\{ \alpha \right\}$. $\alpha$~and $\beta$ are ordinals by Lemmas~\ref{lemm ordinals union intersection} and~\ref{lemm ordinals properties}~(\ref{lemm ordinals properties 4}), and thus $\beta \in A$. Moreover we have $A \subset \beta$. It follows that $\beta \in \beta$, which is a contradiction.
\eproof

We have seen in Lemma~\ref{lemm ordinal well-ordering} that for every ordinal $\alpha$ the relation $\in_{\alpha}$ is a well-ordering in the sense of "$<$". We now show that well-orderings as defined on ordinals are essentially the only well-orderings in the sense of~"$<$" that exist. To this end we first need the notion of isomorphism between two pre-ordered spaces. We then establish some important results between well-orderings that do not explicitly refer to ordinals. Thereby we essentially follow~\cite{Jech}.

\midvspace

\bdefi
\index{Order preserving}
\index{Function!order preserving}
\index{Order isomorphic}
\index{Order isomorphism}
\index{Function!order isomorphic}
\index{Isomorphism}
\index{Isomorphic}
Let $(X,R)$ and $(Y,S)$ be ordered spaces, and $f : X \longrightarrow Y$ a map. $f$~is called {\bf order preserving}, if $(x,y) \in R$ implies $\big( f(x), f(y) \big) \in S$. If $f$ is bijective and $f$ as well as $f^{-1}$ are order preserving, then $f$ is called an {\bf order isomorphism}, or short an {\bf isomorphism}. If such an isomorphism exists, the ordered spaces $(X,R)$ and $(Y,S)$ are called {\bf order isomorphic}, or short {\bf isomorphic}.
\edefi

\brema
\label{rema totally isomorphism}
Let $(X,R)$ be a totally ordered space and $(Y,S)$ be an ordered space, and \mbox{$f : X \longrightarrow Y$} a map. If $f$ is bijective and order preserving, then $f$ is an order isomorphism.
\erema

We recall that the symbol~$<$ always denotes an ordering in the sense of~"$<$".

\midvspace

\blemm
\label{lemm order-pres inequality}
Given a well-ordered space $(X,<)$ and an order preserving map $f : X \longrightarrow X$, we have $x < f(x)$ or $x = f(x)$ for every $x \in X$.
\elemm

\bproof
Let $y$ be the minimum of $A = \left\{ x \in X \, : \, f(x) < x \right\}$. Further let $z = f(y)$. Since $f$ is order preserving, $f(y) < y$ implies $f(z) < z$. Hence $z \in A$ and $z = f(y) < y$, which is a contradiction.
\eproof

\bcoro
\label{coro well-order id}
Given a well-ordered space $(X,<)$, the identity is the only isomorphism from $X$ to itself.
\ecoro

\bproof
This follows by Lemma~\ref{lemm order-pres inequality}.
\eproof

\bcoro
\label{coro well-ordered isomorphic unique}
Given two isomorphic well-ordered spaces $(X,<)$ and $(Y,<)$, there is a unique isomorphism $f : X \longrightarrow Y$.
\ecoro

\bproof
This is a direct consequence of Corollary~\ref{coro well-order id}.
\eproof

\bcoro
\label{coro isom lower segment}
Let $(X,<)$ be a well-ordered space. Then for every subset $A \subset X$, the restriction of the ordering $<$ to~$A$ is a well-ordering in the sense of "$<$" on~$A$, and denoted by~$<$ too. There exists no point $x \in X$ such that $(\, \left] -\infty, x \right[ \,, <)$ and $(X,<)$ are isomorphic.
\ecoro

\bproof
Assume there is such a point $x \in X$ and $f : X \longrightarrow \, \left] -\infty, x \right[ \,$ is an isomophism. Then $f(x) < x$, which is a contradiction by Lemma~\ref{lemm order-pres inequality}.
\eproof

\bprop
\label{prop well-ord restr}
Let $(X,<)$ and $(Y,<)$ be two well-ordered spaces, $f : X \longrightarrow Y$ an isomorphism, and $x \in X$. Then $f( \, \left] -\infty, x \right[ \,) = \, \left] -\infty, f(x) \right[ \,$, and 
\[
f \, | \, \left] -\infty, x \right[ \, \; : \; \left] -\infty, x \right[ \; \longrightarrow \; \left] -\infty, f(x) \right[
\]

is an isomorphism.
\eprop

\bproof
Exercise.
\eproof

\btheo
\label{theo well-ordering isom}
Given two well-ordered spaces $(X,<)$ and $(Y,<)$, exactly one of the following statements is true:

\benum
\item \label{theo well-ordering isom 1} $(X,<)$ and $(Y,<)$ are isomorphic.
\item \label{theo well-ordering isom 2} There is $y \in Y$ such that $(X,<)$ and $(\, \left] -\infty, y \right[ \,, <)$ are isomorphic.
\item \label{theo well-ordering isom 3} There is $x \in X$ such that $(\, \left] -\infty, x \right[ \,, <)$ and $(Y,<)$ are isomorphic.
\eenum

\etheo

\bproof
Following~\cite{Jech} we define the relation $f$ on $X \!\times Y$ by
\[
(x,y) \in f \; \quad \Longleftrightarrow \; \quad (\, \left] -\infty, x \right[ \,, <) \;\, \mathrm{and} \;\, (\, \left] -\infty, y \right[ \,, <) \;\, \mathrm{are \; isomorphic}
\]

Let $D \subset X$ be the domain of $f$ and $R \subset Y$ the range of $f$. Then $f$ is a bijection from $D$ to $R$ by Corollary~\ref{coro isom lower segment}. Let $x_1, x_2 \in D$ with $x_1 < x_2$ and let
\[
g : \; \left] -\infty, x_2 \right[ \; \longrightarrow \; \left] -\infty, f(x_2) \right[ \,
\]

be an isomorphism. Then the restriction
\[
g \, | \, \left] -\infty, x_1 \right[ \, \; : \; \left] -\infty, x_1 \right[ \; \longrightarrow \; \left] -\infty, g(x_1) \right[ \,
\]

is an isomorphism by Proposition~\ref{prop well-ord restr}. Hence $f(x_1) = g(x_1) < f(x_2)$. Therefore $f$ is order preserving. Thus $f$ is an isomorphism from $D$ to $R$ by Remark~\ref{rema totally isomorphism}.

Clearly at most one of the three statements (\ref{theo well-ordering isom 1}) to (\ref{theo well-ordering isom 3}) is true by Corollary~\ref{coro isom lower segment} and Proposition~\ref{prop well-ord restr}.

Now note that either $D = X$ or $D = \, \left] -\infty, x \right[ \;$ for some $x \in X$. Similarly, we either have $R = Y$ or $R = \, \left] -\infty, y \right[ \;$ for some $y \in Y$.
\\

\hspace{0.05\textwidth}
\parbox{0.95\textwidth}
{[If $D \neq X$, then $X \!\setminus\! D$ has a minimum, say~$x$. For every $z \in X$ with $z < x$, we have $z \in D$. On the other hand, for every $z \in X$ with $z \geq x$ we have $ z \notin D$ by Proposition~\ref{prop well-ord restr}. Therefore $D = \, \left] -\infty, x \right[ \;$. The proof for $R$ is similar.]}
\\

Now assume that $D = \, \left] -\infty, x \right[ \;$ for some $x \in X$, and $R = \, \left] -\infty, y \right[ \;$ for some $y \in Y$. Since $D$ and $R$ are isomorphic, we have $(x,y) \in f$ by definition of~$f$, which is a contradiction.
\eproof

\bcoro
\label{coro isom ordinals}
Let $\alpha$ and $\beta$ be two ordinals such that $(\alpha,<)$ and $(\beta,<)$ are isomorphic. Then $\alpha = \beta$.
\ecoro

\bproof
This follows from Lemma~\ref{lemm ordinals properties}~(\ref{lemm ordinals properties 3}) and Theorem~\ref{theo well-ordering isom}.
\\

\hspace{0.05\textwidth}
\parbox{0.95\textwidth}
{[Assume that, under the stated conditions, we have $\alpha \subset \beta$ and $\alpha \neq \beta$. Then we have $\alpha \in \beta$ and $\alpha = \, \left] -\infty, \alpha \right[ \;$ where the interval refers to the ordering~$\in_{\beta}$. Thus the identity is an isomorphism from $\alpha$ to $\, \left] -\infty, \alpha \right[ \;$. It follows that $\beta$ and $\, \left] -\infty, \alpha \right[ \;$ are isomorphic, which is a contradiction by Corollary~\ref{coro isom lower segment}.]}
\eproof

We now establish two important results: first it is shown that for any given set, there is a larger ordinal. Second, for any given well-ordered space that is ordered in the sense of~"$<$", there is a specific ordinal that is isomorphic. We proceed similiarly to~\cite{Ebbinghaus}.

\midvspace

\btheo
\label{theo ordinal injection}
Given a set $X$, there is an ordinal $\alpha$ such that there exists no injection \mbox{$f : \alpha \longrightarrow X$}.
\etheo

\bproof
Assume that $X$ is a set such that for every ordinal $\alpha$ there exists an injection $f : \alpha \longrightarrow X$. We define
\[
{\mc R} = \big\{ R \subset X \!\times\! X \, : \, \mbox{$R$ \,\!\! is a well-ordering on $(\mathrm{field} \, R)$}\big\}
\]

Note that, for every ordinal $\alpha$, there is $A \subset X$ and a well-ordering $R \in {\mc R}$ in the sense of "$<$" on $A$ such that $(\alpha,<)$ and $(A,R)$ are isomorphic.
\\

\hspace{0.05\textwidth}
\parbox{0.95\textwidth}
{[By assumption there is $A \subset X$ and a bijection $f : \alpha \longrightarrow A$. Then
\[
R = \left\{ \big( f(\beta), f(\gamma) \big) \, : \, \beta, \gamma \in \alpha \; \wedge \; \beta < \gamma \right\}
\]

is a well-ordering in the sense of~"$<$" on~$A$. Moreover $(\alpha,<)$ and $(A,R)$ are isomorphic.]}
\\

Now consider the following statement, which we write down partially in the informal language, but we are aware that it could be also written entirely in our formal language with $R$, $X$, and $z$ as free variables:
\\

\hspace{0.05\textwidth}
\parbox{0.95\textwidth}
{If $R \in {\mc R}$ and if there is an ordinal $\alpha$ such that $(\alpha,<)$ and $\left(\mathrm{field} \, R, R \right)$ are isomorphic, then $z = \alpha$, else $z = \O$.}
\\

Clearly, for every set $R$, this statement is true for exactly one set~$z$ by Corollary~\ref{coro isom ordinals}. Thus, by the Replacement schema~\ref{axio repl schema}, the following set exists:
\begin{center}
\begin{tabular}{l}
$W = \big\{ \alpha \, : \, \alpha$ is ordinal,\\[.3em]
 \quad \quad \quad $\exists R \in {\mc R} \quad (\alpha,<)$ and $\left( \mathrm{field} \, R, R \right)$ are isomorphic$\big\}$
\end{tabular}
\end{center}

Now, $W$ contains every ordinal as member by the above result, which is a contradiction to Lemma~\ref{lemm no set ordinals}.
\eproof

\btheo
For every well-ordered space $(X,<)$ there exists a unique isomorphic ordinal $\alpha$.
\etheo

\bproof
Let $\alpha$ be an ordinal such that there is no injection $f : \alpha \longrightarrow X$ by Theorem~\ref{theo ordinal injection}. By Theorem~\ref{theo well-ordering isom} there is $\beta < \alpha$ such that $(\beta,<)$ and $(X,<)$ are isomorphic.

The uniqueness follows by Corollary~\ref{coro isom ordinals}.
\eproof

The following Theorem is a generalization of the Induction principle, Theorem~\ref{lemm induction} (which applies to the set of natural numbers) to any given ordinal number.

\midvspace

\btheo[Induction principle for ordinal numbers]
\label{theo induction ordinals}
\index{Induction principle!ordinal numbers}
Let $\alpha$ be an ordinal and $A \subset \alpha$. If, for every ordinal $\beta$ with $\beta < \alpha$, $\beta \subset A$ implies $\beta \in A$, then $A = \alpha$.
\etheo

\bproof
Assume that, under the stated conditions, $\alpha \!\setminus\! A \neq \O$. Let $\gamma$ be the minimum of $\alpha \!\setminus\! A$. It follows that $\gamma \subset A$, and therefore $\gamma \in A$, which is a contradiction.
\eproof

We remark that again we take a prudent approach in Theorem~\ref{theo induction ordinals} by formulating an Induction principle for an arbitrary ordinal number but not for the collection of all ordinals. This approach allows us to state the Theorem in terms of sets---while the collection of all ordinals is not a set. One may derive a similar theorem for the class of all ordinals as done in~\cite{Jech}, or in terms of a formula holding for every ordinal as done in~\cite{Ebbinghaus}.

The following result is a generalization of the Recursion theorem for natural numbers, Theorem~\ref{theo recursive def}, to any given ordinal number.

\midvspace

\btheo[Local recursion]
\label{theo recursion}
\index{Local recursion}
Let $\alpha$ be an ordinal, $X$ a set, $D = \bigcup \left\{ X^{\beta}  : \, \beta < \alpha \right\}$, and $F : D \longrightarrow X$ a map. There is a unique function $f : \alpha \longrightarrow X$ such that $f(\beta) = F (f\,|\,\beta)$ for every $\beta < \alpha$.
\etheo

\bproof
Notice that if such a function exists, then it is unique.
\\

\hspace{0.05\textwidth}
\parbox{0.95\textwidth}
{[Assume $f$ and $g$ are two such functions and $f \neq g$. Let $A = \left\{ \beta \, : \, \beta < \alpha,\; f(\beta) \neq g(\beta) \right\}$ and $\gamma = \mathrm{min} \, A$. Then we have $f \, | \, \gamma = g \, | \, \gamma$, and therefore $f(\gamma) = F(f \, |\,\gamma) = F(g \, | \, \gamma) = g(\gamma)$, which is a contradiction.]}
\\

Now assume there exists no such function. Let $M \subset \alpha \cup \left\{ \alpha \right\}$ be the set of those ordinals $\gamma < \alpha \cup \left\{ \alpha \right\}$ for which there is no function $f : \alpha \longrightarrow X$ such that $f(\beta) = F(f \, | \, \beta)$ for every $\beta < \gamma$. We have $\alpha \in M$ by assumption. Let $\delta$ be the minimum of~$M$. Then for every $\varepsilon < \delta$ there exists a function $f : \alpha \longrightarrow X$ such that $f(\xi) = F(f \, | \, \xi)$ for every $\xi < \varepsilon$. For every $\varepsilon < \delta$ we may denote the system of such functions by~$F_{\varepsilon}$. Moreover, if $\varepsilon_1, \varepsilon_2 < \delta$ with $\varepsilon_1 < \varepsilon_2$ or $\varepsilon_1 = \varepsilon_2$, and $f_1 \in F_{\varepsilon(1)}$, $f_2 \in F_{\varepsilon(2)}$, then we have $f_1 \, |\, \varepsilon_1 = f_2 \, |\, \varepsilon_1$ by the result above.

We now distinguish two cases. If $\delta$ has no predecessor, then we may choose a point $x \in X$ and define the map $g : \alpha \longrightarrow X$ by
\[
g(\xi) = \left\{
\begin{array}{ll}
f(\xi) \quad & \mathrm{if} \quad \xi < \delta, \quad \mathrm{where} \; f \in F_{\varepsilon}, \; \xi < \varepsilon < \delta\\[.5em]
x \quad & \mathrm{if} \quad \xi = \delta \;\; \mathrm{or} \;\; \delta < \xi
\end{array}
\right.
\]

Now, if $\delta$ has a predecessor, say $\varepsilon$, we may choose $x \in X$, $f \in F_{\varepsilon}$, and define $g : \alpha \longrightarrow X$ by
\[
g(\xi) = \left\{
\begin{array}{ll}
f(\xi) \quad & \mathrm{if} \quad \xi < \varepsilon\\[.5em]
F(f \, | \, \varepsilon) \quad & \mathrm{if} \quad \xi = \varepsilon\\[.5em]
x \quad & \mathrm{if} \quad \varepsilon < \xi
\end{array}
\right.
\]

In both cases we find that $g(\xi) = F(g \, | \, \xi)$ for every $\xi < \delta$, which is a contradiction.
\eproof

\section{Choice}
\label{choice}

Several important consequences of the Choice axiom, Axiom~\ref{axio choice}, that are required subsequently are proven in this Section. Since our aim is not an exhaustive discussion of~ZFC, but to explain its consequences for the foundations of analysis, we do not prove the equivalence of the different forms of the Choice axiom neither provide a comprehensive list of known equivalent forms. Instead, we always assume the Choice axiom and derive some of its relevant implications. Remember that the Choice axiom implies the existence of a choice function as stated in Lemma and Definition~\ref{lede choice function} and Remark~\ref{rema choice}.

\midvspace

\btheo[Well-ordering principle]
\label{theo aoc well-ordering}
\index{Well-ordering principle}
For every set $X$ there exists a well-ordering on~$X$.
\etheo

\bproof
We may choose an ordinal $\beta$ such that there is no injection from $\beta$ to $X$ by Theorem~\ref{theo ordinal injection}. Moreover, let $g : {\mc P}(X) \!\setminus\! \left\{ \O \right\} \longrightarrow X$ a choice function and $y$ a point such that $y \notin X$. Further let
\[
Y = X \cup \left\{ y \right\}, \quad D = \bigcup \left\{ Y^{\gamma} : \, \gamma < \beta \right\}
\]

We define the map $F : D \longrightarrow Y$ by
\[
F(h) = \left\{
\begin{array}{ll}
g \big( X \!\setminus\! \mathrm{ran} \, h \big) \quad & \mathrm{if} \quad X \not\subset \mathrm{ran} \, h \\[.5em]
y \quad & \mathrm{if} \quad X \subset \mathrm{ran} \, h \end{array}
\right.
\]

By Theorem~\ref{theo recursion} there exists a unique function $f : \beta \longrightarrow Y$ such that $f(\gamma) = F (f \, | \, \gamma)$ for every $\gamma < \beta$.

$f^{-1} \left[ X \right]$ is an ordinal, since it is full and the relation $<$ on this set has the minimum property by Lemma~\ref{lemm ordinals minimum}.
\\

\hspace{0.05\textwidth}
\parbox{0.95\textwidth}
{[Let $\delta < \gamma \in f^{-1} \left[ X \right]$. Then $f(\gamma) \in X$, and thus $X \not\subset f \left[ \gamma \right]$. It follows that $X \not\subset f \left[ \delta \right]$. Therefore $f(\delta) \in X$.]}
\\

We define $\alpha = f^{-1} \left[ X \right]$.

For every $x \in X$ the set $f^{-1} \left\{ x \right\}$ is a singleton or empty.
\\

\hspace{0.05\textwidth}
\parbox{0.95\textwidth}
{[Assume $\gamma, \delta \in f^{-1} \left\{ x \right\}$ with $\gamma < \delta$. Then we have $f(\delta) \in X$, and thus $f(\delta) = g \big( X \!\setminus\! f \left[ \delta \right] \big) \neq x$, which is a contradiction.]}
\\

Moreover $X \subset f \left[ \beta \right]$. 
\\

\hspace{0.05\textwidth}
\parbox{0.95\textwidth}
{[Assume there is $x \in X$ such that $x \notin f \left[ \beta \right]$. Let $\gamma < \beta$. Then $X \not\subset f \left[ \gamma \right]$, and therefore $f(\gamma) = g \big( X \!\setminus\! f \left[ \gamma \right] \big) \in X$. Thus $f$ is an injection that maps $\beta$ on~$X$, which is a contradiction.]}
\\

It follows that the map $t : \alpha \longrightarrow X$, $t = f \, | \, \alpha$ is bijective.

We now define a relation $<$ on $X$ by
\[
x < y \quad \Longleftrightarrow \quad t^{-1} \left\{ x \right\} < t^{-1} \left\{ y \right\}
\]

This relation is clearly a well-ordering on~$X$.
\eproof

Notice that Theorem~\ref{theo aoc well-ordering} and Lemma~\ref{lemm both well-orderings} imply that for every set $X$ there exists a well-ordering in the sense of "$<$" and a well-ordering in the sense of "$\leq$". Theorem~\ref{theo aoc well-ordering} is used to prove the pseudo-metrization theorem.

\midvspace

\bcoro
\label{coro ordinal bijection}
For every set $X$ there exists an ordinal~$\alpha$ and a bijection $t : \alpha \longrightarrow X$.
\ecoro

\bproof
Such a bijection is explicitily constructed in the proof of Theorem~\ref{theo aoc well-ordering}.
\eproof

The next result is another important implication of the Choice axiom. Our proof does not make use of the Well-ordering principle.

\midvspace

\btheo[Zorn's Lemma]
\label{theo aoc Zorn}
\index{Zorn's Lemma}
Let $(X,R)$ be an ordered space. If every chain has an upper bound, then $X$ has a weak maximum.
\etheo

\bproof
We may choose an ordinal $\beta$ such that there is no injection from $\beta$ to~$X$ by Theorem~\ref{theo ordinal injection}. Let $g : {\mc P}(X) \!\setminus\! \left\{ \O \right\} \longrightarrow X$ be a choice function, and ${\mc A} \subset {\mc P}(X)$ the system of all chains. Moreover let
\[
H : {\mc A} \longrightarrow {\mc P}(X), \quad H(A) = \big\{ z \in X \!\setminus\! A \; : \; A \cup \left\{ z \right\} \in {\mc A} \big\}
\]

Further let
\[
x \in X, \quad y \notin X, \quad Y = X \cup \left\{ y \right\}, \quad D = \bigcup \left\{ Y^{\gamma} : \, \gamma < \beta \right\},
\]

We define the map \mbox{$F : D \longrightarrow Y$} by
\[
F(h) = \left\{
\begin{array}{ll}
x \quad & \mathrm{if} \quad h = \O\\[.5em]
g \, \big( H ( \mathrm{ran} \, h) \big) \quad & \mathrm{if} \quad h \neq \O,\;\; \mathrm{ran} \, h \in {\mc A},\;\; H (\mathrm{ran} \, h) \neq \O\\[.5em]
y \quad & \mathrm{if} \quad h \neq \O,\;\; \mathrm{ran} \, h \in {\mc A},\;\; H (\mathrm{ran} \, h) = \O\\[.5em]
y \quad & \mathrm{if} \quad h \neq \O,\;\; \mathrm{ran} \, h \notin {\mc A}
\end{array}
\right.
\]

By the Local recursion theorem there is a unique function $f : \beta \longrightarrow Y$ such that $f(\gamma) = F( f \, | \, \gamma)$ for every $\gamma < \beta$. We define $\alpha = f^{-1} \left[ X \right]$. $\alpha$~is an ordinal since it is full and the relation $<$ is a well-ordering on~$\alpha$ by Lemma~\ref{lemm ordinals minimum}.
\\

\hspace{0.05\textwidth}
\parbox{0.95\textwidth}
{[Let $\delta < \gamma < \alpha$. Then $f(\gamma) \in X$, and thus $f \left[ \gamma \right] \in {\mc A}$ and $H (f \left[ \gamma \right]) \neq \O$. It follows that $f \left[ \delta \right] \in {\mc A}$ and $H (f \left[ \delta \right]) \neq \O$. Therefore $f(\delta) \in X$.]}
\\

For every $z \in X$ the set $f^{-1} \left\{ z \right\}$ is a singleton or empty.
\\

\hspace{0.05\textwidth}
\parbox{0.95\textwidth}
{[Assume $\gamma, \delta \in f^{-1} \left\{ z \right\}$ with $\gamma < \delta$. Then we have $f(\delta) \in X$, and thus $f(\delta) = g \big( H (f \left[ \delta \right] )\big) \neq z$, which is a contradiction.]}
\\

Therefore we have $\alpha < \beta$ by the choice of~$\beta$. We define $B = f \left[ \alpha \right]$. Then $B$ is a chain.
\\

\hspace{0.05\textwidth}
\parbox{0.95\textwidth}
{[Let $\gamma, \delta \in \alpha$ with $\gamma < \delta$. Then $f(\delta) = g \big( H (f \left[ \delta \right]) \big)$, and therefore $f \left[ \delta \right] \cup \left\{ f(\delta) \right\}$ is a chain. We have $\gamma \in \delta$, and thus $f(\gamma) \in f \left[ \delta \right]$.]}
\\

Furthermore, $H(B) = \O$.
\\

\hspace{0.05\textwidth}
\parbox{0.95\textwidth}
{[Assume that $H(B) \neq \O$. It follows that $f(\alpha) = g(H(B)) \in X$, and hence $\alpha \in \alpha$, which is a contradiction.]}
\\

Let $b \in X$ be an upper bound of~$B$, which exists by assumption. Then $b$ is a weak maximum of~$X$.
\\

\hspace{0.05\textwidth}
\parbox{0.95\textwidth}
{[Assume that $(b,c) \in R$. Then $B \cup \left\{ b, c \right\}$ is a chain, and thus $b, c \in B$. It follows that $(c, b) \in R$.]}
\eproof

Finally we generalize Zorn's Lemma to pre-ordered spaces, which is applied in the proof of the existence of an ultrafilter base finer than a given filter base in Theorem~\ref{theo existence ultrafilter base}.

\midvspace

\btheo
\label{theo Zorn}
\index{Zorn's Theorem}
Let $(X,R)$ be a pre-ordered space. If every chain has an upper bound, then $X$ has a weak maximum.
\etheo

\bproof
Assume the stated condition. Let $(Y,S)$ be the ordered space constructed from~$(X,R)$ as in Lemma~\ref{lemm pre-ordering ordering}. Every chain in~$Y$ has an upper bound.
\\

\hspace{0.05\textwidth}
\parbox{0.95\textwidth}
{[Let $A \subset Y$ be a chain. Then $B = \bigcup A$ is chain too. Therefore $B$ has an upper bound, say~$x$. Then $\left[ x \right]$ is an upper bound of~$A$.]}
\\

Let $b$ be the weak maximum of~$Y$ by Theorem~\ref{theo aoc Zorn}. Then every point $x \in b$ is a weak maximum of~$X$ (exercise).
\eproof

\section{Cardinality}
\label{cardinality}

In this Section we define several notions in order to describe what could be considered as the size of a set.

\midvspace

\bdefi
\index{Cardinality}
\index{Finite}
\index{Infinite}
\index{Countable}
\index{Uncountable}
Two sets $X$ and $Y$ are said to be {\bf of the same cardinality}, also written $X \sim Y$, if there exists a bijection $f : X \longrightarrow Y$. Let $X$ be a set. $X$~is called {\bf finite} if there is $m \in \naturalnumbers$ such that $m \sim X$, else it is called {\bf infinite}. If $X \sim \naturalnumbers$ or $X$ is finite, then $X$ is called {\bf countable}. If $X$ is not countable, it is called {\bf uncountable}.
\edefi

\brema
\label{rema ordinal same card}
For every set $X$ there is an ordinal~$\alpha$ such that $X \sim \alpha$ by Corollary~\ref{coro ordinal bijection}.
\erema

\blemm
Let $X$ be a set. If $X$ is infinite, then for every $m \in \naturalnumbers$ there is a subset $Y \subset X$ such that $Y \sim m$. If $X$ is uncountable, then there is a subset $Y \subset X$ such that $Y \sim \naturalnumbers$.
\elemm

\bproof
This follows by Remarks~\ref{rema ordinal vs natural} and~\ref{rema ordinal same card}.
\eproof

\blemm
\label{lemm tot ord subset min max}
Let $(X,\prec)$ be a connective pre-ordered space and $A \subset X$ where $A$ is finite and $A \neq \O$. Then $A$ has a minimum and a maximum. If $\prec$ is a total ordering, then the minimum and the maximum of~$A$ are unique. 
\elemm

\bproof
If $A \sim 1$, then $A$ is a singleton. Thus it has a unique minimum and a unique maximum.

To see the first claim, assume that every $A \subset X$ with $A \sim m$ for some $m \in \naturalnumbers$, $m > 0$, has a minimum and a maximum. Let $B \subset X$ with $B \sim \sigma(m)$. There is a bijection $f : \sigma(m) \longrightarrow B$. By assumption $B \setminus\! \left\{ f(m) \right\}$ has a minimum, say~$x$, and a maximum, say~$y$. Then the set $\left\{ x, f(m) \right\}$ has a minimum, which is a minimum of~$B$, and the set $\left\{ y, f(m) \right\}$ has a maximum, which is a maximum of~$B$.

The second claim follows by Remark~\ref{rema max unique}.
\eproof

\bcoro
Let $X$ be a set. If $X \sim \naturalnumbers$, then $X$ is infinite.
\ecoro

\bproof
It is enough to show that $\naturalnumbers$ is infinite. Let $m \in \naturalnumbers$, $m > 0$, and $f : m \longrightarrow \naturalnumbers$ be an injection. Then the set $f \left[ m \right]$ is finite and not empty, and therefore has a unique maximum by Lemma~\ref{lemm tot ord subset min max}. Thus $f$ is not surjective.
\eproof

\bprop
\label{prop finite predecessor}
Let $X$ be a finite non-empty set, $x \in X$, and $m \in \naturalnumbers$. Then $X \sim \sigma(m)$ implies $X \!\setminus\! \left\{ x \right\} \sim m$.
\eprop

\bproof
Exercise.
\eproof

\blemm
\label{lemm finite sets}
Let $X$ and $Y$ be finite sets and $Z \subset X$. Then $Z$, $X \cup Y$, and $X \!\times Y$ are finite.
\elemm

\bproof
We first show that $Z$ is finite. This is clearly true if $X \sim 0$. Now let $m \in \naturalnumbers$. Assume that the claim is true for every $X$ with $X \sim m$. Let $U$ be a set with $U \sim \sigma(m)$. Further let $u \in U$, $V = U \!\setminus\! \left\{ u \right\}$, and $Z \subset U$. Then we have $V \sim m$ by Proposition~\ref{prop finite predecessor}, and hence $Z \cap V$ is finite by assumption. If $Z \subset V$, then $Z$ is finite. If $Z \not\subset V$, we may choose $n \in \naturalnumbers$ and a bijection $g : n \longrightarrow Z \cap V$. We define a bijection $h : \sigma(n) \longrightarrow Z$ by $h \, | \, n = g$ and $h(n) = u$, and therefore $Z$ is finite.

Next we prove that $X \cup Y$ is finite. Since in the case $X = Y = \O$ the claim is obvious, we show that $X \sim m$ and $Y \sim n$, where $m, n \in \naturalnumbers$, $0 < m$, implies that there is an injection $h : X \cup Y \longrightarrow m + n$. Let $m \in \naturalnumbers$, $0 < m$. Then this implication is clearly true for $n = 0$. Assuming that it is true for some $n \in \naturalnumbers$, let $X$ and $Y$ be sets with $X \sim m$, $Y \sim \sigma(n)$. Further let $y \in Y$, and $V = Y \!\setminus\! \left\{ y \right\}$. Then we have $V \sim n$, and there is an injection $f : X \cup V \longrightarrow m + n$ by assumption. We define the injection $h : X \cup Y \longrightarrow m + \sigma(n)$ by $h \, |\left(X \cup V\right) = f$ and, if $y \notin X$, then $h(y) = m + n$.

To see that $X \!\times Y$ is finite, first notice that this is clear if $X = \O$ or $Y = \O$. Now assume $X \neq \O$ and $Y \neq \O$. We show that $X \sim m$ and $Y \sim n$ ($m, n \in \naturalnumbers \setminus\! \left\{ 0 \right\}$) implies that $X \!\times Y \sim m \, n$. Let $m \in \naturalnumbers$, $0 < m$. Then this implication is clearly true for $n = 1$. Assuming that it is true for some $n \in \naturalnumbers$, $0 < n$, let $X$ and $Y$ be sets with $X \sim m$, $Y \sim \sigma(n)$. Further let $y \in Y$, $V = Y \!\setminus\! \left\{ y \right\}$, and $h : m \longrightarrow X$ a bijection. Then there is a bijection $f : X \!\times V \longrightarrow m \, n$ by assumption. We define a bijection $s : X \!\times Y \longrightarrow m \, \sigma(n)$ by $s \, |\left(X \times V\right) = f$ and $s\big((h(k), y)\big) = m \, n + \sigma(k)$ for every $k \in \naturalnumbers$, $k < m$.
\eproof

The following Lemma extends the results of Lemma~\ref{lemm finite sets} to arbitrary finite unions and products.

\midvspace

\blemm
\label{lemm finite union product}
Let $I$ be a finite index set, and for every $i \in I$, let $X_i$ be a finite non-empty set. Then $\bigcup_{i \in I} X_i$ and $\bigtimes_{\!\! i \in I}\, X_i$ are finite.
\elemm

\bproof
We first show that $\bigcup_{i \in I} X_i$ is finite. This is clear if $I \sim 1$. Now let $m \in \naturalnumbers$, $0 < m$, and assume that the claim holds for every index set $I$ with $I \sim m$. Let $J$ be a set with $J \sim \sigma(m)$. Further let $k \in J$, $K = J \!\setminus\! \left\{ k \right\}$, and $X_j$ a finite non-empty set for every $j \in J$. Then $K \sim m$ by Proposition~\ref{prop finite predecessor}, and $\bigcup_{j \in K} X_j$ is finite by assumption. It follows that $\bigcup_{j \in J} X_j = \bigcup_{j \in K} X_j \cup X_k$ is finite by Lemma~\ref{lemm finite sets}.

We now prove that $\bigtimes_{\!\! i \in I}\, X_i$ is finite. The claim clearly holds if $I \sim 1$. Now let $m \in \naturalnumbers$, $0 < m$, and assume that the claim holds for every index set $I$ with $I \sim m$. Let $J$ be a set with $J \sim \sigma(m)$. Further let $k \in J$, $K = J \!\setminus\! \left\{ k \right\}$, and $X_j$ a finite non-empty set for every $j \in J$. Then $\bigtimes_{\!\! j \in K}\, X_j$ is finite by assumption. It follows that $\bigtimes_{\!\! j \in K}\, X_j \times X_k$ is finite by Lemma~\ref{lemm finite sets}. Thus $\bigtimes_{\!\! j \in J}\, X_j$ is finite by Remark~\ref{rema cartesian mixed}.
\eproof

\brema
Let $X$ and $Y$ be finite non-empty sets. Then $X^Y$ is finite.
\erema

\bprop
\label{prop count subset}
Let $X$ be a countable set. There is a set $Y$ such that $X \subset Y$ and $Y \sim \naturalnumbers$.
\eprop

\bproof
If $X$ is infinite, the claim clearly holds. If $X$ is finite, then there is $m \in \naturalnumbers$ and a bijection $f : m \longrightarrow X$. If the set $X \cap \naturalnumbers$ is non-empty, then let $n$ be its unique maximum by Lemma~\ref{lemm tot ord subset min max}. In this case we may choose a number $p \in \naturalnumbers$ such that $n < p$ and $m < p$. If $X \cap \naturalnumbers$ is empty, then let $p = \sigma(m)$. We define $Y = X \cup \left\{ q \in \naturalnumbers \, : \, q \geq p \right\}$ and the map $g : \naturalnumbers \longrightarrow Y$ by
\[
g(r) = \left\{\begin{array}{ll}
f(r) & \mbox{if $r < m$}\\
p & \mbox{if $r = m$}\\
p + s & \mbox{if $r > m$}
\end{array}\right.
\]

where, in the last case, $s$ is the number such that $m + s = r$ by Lemma~\ref{lemm natural diff}. Then $g$ is bijective.
\eproof

\blemm
\label{lemm countable}
A set $X$ is countable iff there is an injection $f : X \longrightarrow \naturalnumbers$.
\elemm

\bproof
Assume there is an injection $f : X \longrightarrow \naturalnumbers$. Let $A = \mathrm{ran} \, f$. If $A$ is finite, then $X$ is finite. If $A$ is infinite, we may recursively define the function $g : \naturalnumbers \longrightarrow A$ by Lemma~\ref{theo recursive def} and Lemma~\ref{lemm ordinals union intersection} as follows: Let $g(0)$ be the minimum of~$A$, and, for every $m \in \naturalnumbers$, let $g(\sigma(m))$ be the minimum of $\left\{ n \in A \, : \, n > g(m) \right\}$. Then $g$ is strictly increasing.
\\

\hspace{0.05\textwidth}
\parbox{0.95\textwidth}
{[We show by the Induction principle that $g$ is strictly increasing, i.e.\ $m < n$ implies $g(m) < g(n)$ for every $m, n \in \naturalnumbers$. Let $m \in \naturalnumbers$. Then the claim clearly holds for $n = \sigma(m)$. Assume it holds for some $n \in \naturalnumbers$ with $m < n$. Then $g(\sigma(n)) > g(n)$, and therefore the claim also holds for $\sigma(n)$.]}
\\

Since $g$ is strictly increasing, it is injective. To see that $g$ is surjective, assume there exists $m \in A \!\setminus\! \mathrm{ran} \, g$. Then we have \,\! $\mathrm{ran} \, g \subset m$.
\\

\hspace{0.05\textwidth}
\parbox{0.95\textwidth}
{[Let $B = \left\{ n \in \naturalnumbers \, : \, g(n) > m \right\}$, and assume that $B \neq \O$. Let $p$ be the minimum of~$B$ by Lemma~\ref{lemm ordinals union intersection}. We have $p > 0$ by the definition of~$g(0)$. Let $q$ be the predecessor of~$p$. Then $g(q) < m$. It follows that $g(\sigma(q)) \leq m$ since $m \in A$, which is a contradiction.]}
\\

Hence $\mathrm{ran} \, g$ is finite, which is a contradiction. Thus $g$ is bijective, and therefore $X$ is countable.

The converse follows by Theorem~\ref{theo naturals}~(\ref{theo naturals 1}).
\eproof

\bprop
\label{prop naturalssqu}
We have $\naturalnumbers^2 \sim \naturalnumbers$.
\eprop

\bproof
We define the functions

\begin{center}
\begin{tabular}{ll}
$s : \naturalnumbers \longrightarrow \naturalnumbers$ , & $s(m) = m (m + 1) / 2$ ;\\[.8em]
$h : \naturalnumbers^2 \longrightarrow \naturalnumbers$ , & $h(p, q) = s(p + q) + q$
\end{tabular}
\end{center}

Note that $s$~is well-defined by Proposition~\ref{prop gauss integer}.

We first show that $h$ is injective. For given $r \in \naturalnumbers$ there is at most one pair $(m,q) \in \naturalnumbers^2$ such that $q \leq m$ and $s(m) + q = r$.
\\

\hspace{0.05\textwidth}
\parbox{0.95\textwidth}
{[For $i \in \left\{ 1, 2 \right\}$, assume that $m_i, q_i \in \naturalnumbers$ with $q_i \leq m_i$ and $s(m_i) + q_i = r$. Further we may assume that $m_1 < m_2$. Since $s$ is increasing, we have
\[
s(m_2) + q_1 \geq s(m_1 + 1) + q_1 = s(m_1) + m_1 + 1 + q_1 = r + m_1 + 1 > r + q_1
\]

It follows that $s(m_2) > r$, which is a contradiction.]}
\\

Thus there is at most one pair $(p,q) \in \naturalnumbers^2$ such that $s(p + q) + q = r$.

Thus $\naturalnumbers^2$ is countable by Lemma~\ref{lemm countable}. This set is clearly not finite, so it is of the same cardinality than~$\naturalnumbers$.
\eproof

\blemm
\label{lemm countable sets}
Let $X$ and $Y$ be countable sets and $Z \subset X$. Then $Z$ and $X \!\times Y$ are countable.
\elemm

\bproof
$Z$ is countable by Lemma~\ref{lemm countable}.

In order to show that $X \times Y$ is countable, we may choose two injections \mbox{$f : X \longrightarrow \naturalnumbers$} and \mbox{$g : Y \longrightarrow \naturalnumbers$} by Lemma~\ref{lemm countable}, and an injection \mbox{$h: \naturalnumbers^2 \longrightarrow \naturalnumbers$} by Proposition~\ref{prop naturalssqu}. We define the function
\[
t : X \!\times Y \longrightarrow \naturalnumbers, \quad t(x,y) = h \big( f(x), g(y) \big)
\]

$t$ is clearly an injection. The claim follows by Lemma~\ref{lemm countable}.
\eproof

\blemm
\label{lemm count times count}
\index{Countable}
\index{Set!countable}
Let $I$ be a countable index set, and for each $i \in I$ let $X_i$ be a countable set. Then $X = \bigcup_{i \in I} X_i$ is countable.
\elemm

\bproof
We may assume that $I \sim \naturalnumbers$ and $X_i \sim \naturalnumbers$ for every $i \in I$.
\\

\hspace{0.05\textwidth}
\parbox{0.95\textwidth}
{[If $I$ is countable, then there is a set~$J$ and a bijection $f : J \longrightarrow \naturalnumbers$ such that $I \subset J$ by Proposition~\ref{prop count subset}. For every $j \in J$, we may choose a set $Y_j$ and a bijection $f_j : Y_j \longrightarrow \naturalnumbers$ such that $X_j \subset Y_j$ ($j \in I$) by the same Proposition. Then $X \subset Y$ where $Y = \bigcup_{j \in J} Y_j$. If there is an injection $t : Y \longrightarrow \naturalnumbers$, then $t \, | \, X$ is an injection.]}
\\

Let $g : \naturalnumbers \longrightarrow I$ and $g_i : \naturalnumbers \longrightarrow X_i$ ($i \in I$) be bijections. We may choose a bijection $h : \naturalnumbers \longrightarrow \naturalnumbers^2$ by Lemma~\ref{lemm countable sets}. For $i \in \left\{ 1, 2 \right\}$, let $p_i : \naturalnumbers^2 \longrightarrow \naturalnumbers$ be the projections on the coordinates. We define
\[
G : \naturalnumbers \longrightarrow X, \quad \quad G(m) = g_{g \, p_1 h (m)} \big( p_2 \, h (m) \big)
\]

$G$~is clearly surjective. Let $H : X \longrightarrow \naturalnumbers$ where, for every $x \in X$, $H(x)$ is the minimum of~$G^{-1} \left\{ x \right\}$, which exists since $<$ is a well-ordering on~$\naturalnumbers$ by Remark~\ref{rema well-ordering naturals}. Then $H$ is injective. The claim follows by Lemma~\ref{lemm countable}.
\eproof

\blemm
\label{lemm fin count}
\index{Countable}
\index{Set!countable}
Let $X$ be a countable non-empty set and ${\mc A} = \left\{ A \subset X \, : \, \mbox{$A$ is finite} \right\}$. Then $\mc A$ is countable.
\elemm

\bproof
We have ${\mc A} = \bigcup_{m \in \naturalnumbers} {\mc A}_m$ where ${\mc A}_m = \left\{ A \subset X \, : \, A \sim m \right\}$ ($m \in \naturalnumbers$).

First we show by the Induction principle that, for every $m \in \naturalnumbers$, ${\mc A}_m$ is countable. Clearly, ${\mc A}_0$ is finite. Now assume that ${\mc A}_m$ is countable for some $m \in \naturalnumbers$. We have
\begin{eqnarray*}
{\mc A}_{\sigma(m)} \!\!\! & = & \!\! \big\{ A \subset X \, : \, \exists x \in A \quad A \!\setminus\! \left\{ x \right\} \in {\mc A}_m \big\}\\[.5em]
 & \subset & \!\!\! \bigcup_{x \in X} \big\{ A \cup \left\{ x \right\} \, : \, A \in {\mc A}_m \big\}
\end{eqnarray*}

which is countable by Lemma~\ref{lemm count times count}.

Now it follows that $\mc A$ is countable by the same Lemma.
\eproof


\chapter{Numbers II}
\label{numbers ii}
\setcounter{equation}{0}

\pagebreak

\section{Positive dyadic rational numbers}
\label{positive dyadic rational numbers}

In this Section we define the positive dyadic rational numbers, and in the next Section the positive real numbers. It is then possible to construct the full system of real numbers from its positive counterpart. The set of (positive and negative) integers and the set of (positive and negative) dyadic rationals can finally be identified with subsets of the reals.

\midvspace

\blede
\label{lede pos dyadic rational numbers}
\index{Positive dyadic rational numbers}
\index{Numbers!positive dyadic rational}
We define an equivalence relation $Q$ on $\naturalnumbers^2$ by
\[
\big((m,u), (n,v)\big) \in Q \quad \Longleftrightarrow \quad m \, 2^v = n \, 2^u
\]

and ${\mathbb D}_+ = \naturalnumbers^2 / Q$. For every $m, u \in \naturalnumbers$, the equivalence class of $(m,u)$ is denoted by $\lfloor m, u \rfloor$. The members of ${\mathbb D}_+$ are called {\bf positive dyadic rational numbers}.

Furthermore, we define the relation $<$ on ${\mathbb D}_+$ as follows:
\[
\lfloor m, u \rfloor < \lfloor n, v \rfloor \quad \Longleftrightarrow \quad m \, 2^v < n \, 2^u
\]

This is a total ordering in the sense of "$<$" on~${\mathbb D}_+$.
Moreover, we define $\leq$ to be the total ordering in the sense of "$\leq$" on~${\mathbb D}_+$ obtained from the ordering~$<$ by the method of Lemma~\ref{lemm both orderings}.
\elede

\bproof
Clearly, $Q$ is an equivalence relation. To see that the relation $<$ on ${\mathbb D}_+$ is well defined, let $\lfloor m, u \rfloor < \lfloor n, v \rfloor$, and $(p,w) \in \lfloor m, u \rfloor$, $(q,r) \in \lfloor n, v \rfloor$. Then we have $p \, 2^u = m \, 2^w$ and $q \, 2^v = n \, 2^r$. It follows that
\[
p \, 2^{u+v+r} = m \, 2^{w+v+r} < n \, 2^{u+w+r} = q \, 2^{v+u+w}
\]

and thus $p \, 2^r < q \, 2^w$. To see that it is transitive let $\lfloor m, u \rfloor, \lfloor n, v \rfloor, \lfloor p, w \rfloor \in {\mathbb D}_+$ with $\lfloor m, u \rfloor < \lfloor n, v \rfloor < \lfloor p, w \rfloor$. Then we have $m \, 2^v < n \, 2^u$ and $n \, 2^w < p \, 2^v$. It follows that $m \, 2^{v+w} < n \, 2^{u+w} < p \, 2^{v+u}$, and therefore $m \, 2^w < p \, 2^u$. Thus we obtain $\lfloor m, u \rfloor < \lfloor p, w \rfloor$. Moreover, it is obviously antireflexive and connective.
\eproof

\bdefi
\label{defi convention ordering pos dyadic numbers}
We adopt the convention that all notions related to orderings on~${\mathbb D}_+$\,, in particular intervals, refer to the ordering in the sense of~"$<$" as defined in Lemma and Definition~\ref{lede pos dyadic rational numbers} unless otherwise specified.
\edefi

Notice that this convention agrees with the one in the context of natural numbers, cf.\ Definition~\ref{defi convention ordering natural numbers}. Note again that in many cases it is irrelevant whether the ordering~$<$ or the ordering~$\leq$ on~${\mathbb D}_+$ is considered as most order properties are invariant, cf.\ Lemmas~\ref{lemm max min inv}, \ref{lemm bound sup inv}, and~\ref{lemm mon invar}. However, the choice of the ordering {\it is} relevant for intervals.

\midvspace

\bcoro
\label{coro dyadic countable}
The set ${\mathbb D}_+$ is countable.
\ecoro

\bproof
This follows from Lemma~\ref{lemm countable sets}.
\eproof

\blemm
\label{lemm dyadic order dense}
${\mathbb D}_+$ is $<$-dense.
\elemm

\bproof
Let $\lfloor m, u \rfloor, \lfloor n, v \rfloor \in {\mathbb D}_+$ with $\lfloor m, u \rfloor < \lfloor n, v \rfloor$. We define $m_0 = m \, 2^{v + 1}$, and $n_0 = n \, 2^{u + 1}$. Then we have $\lfloor m, u \rfloor = \lfloor m_0, u + v + 1 \rfloor$ and $\lfloor n, v \rfloor = \lfloor n_0, u + v + 1 \rfloor$. Since $m_0$ and $n_0$ are even, there is $k \in \naturalnumbers$ such that $m_0 < k < n_0$. Thus $\lfloor m, u \rfloor < \lfloor k, u + v + 1 \rfloor < \lfloor n, v \rfloor$. 
\eproof

\blede
\index{Addition!positive dyadic rational numbers}
\index{Sum!positive dyadic rational numbers}
Let $+$ be the binary function on ${\mathbb D}_+$ defined by
\[
\lfloor m, u \rfloor + \lfloor n, v \rfloor = \lfloor m \, 2^v + n \, 2^u, \; u + v \rfloor
\]

This function is called {\bf addition on}~${\mathbb D}_+$. The expression $\lfloor m, u \rfloor + \lfloor n, v \rfloor$ is called the {\bf sum of $\lfloor m, u \rfloor$ and $\lfloor n, v \rfloor$}. $+$~is commutative and associative. Moreover, for every $d, e, f \in {\mathbb D}_+$, we have
\[
d < e \quad \Longrightarrow \quad d + f < e + f
\]

\elede

\bproof
To see that $+$ is well-defined, let $\lfloor m, u \rfloor = \lfloor p, w \rfloor \in {\mathbb D}_+$ and $\lfloor n, v \rfloor = \lfloor q, r \rfloor \in {\mathbb D}_+$. Then we have
\begin{eqnarray*}
\lfloor m, u \rfloor + \lfloor n, v \rfloor \!\! & = &  \!\! \lfloor m \, 2^v + n \, 2^u, \; u + v \rfloor\\
& = & \!\! \lfloor m \, 2^{v + w + r} + n \, 2^{u + w + r}, \; u + v + w + r \rfloor\\
& = & \!\! \lfloor p \, 2^{v + u + r} + q \, 2^{u + w + v}, \; u + v + w + r \rfloor\\
& = & \!\! \lfloor p \, 2^r + q \, 2^w, \; w + r \rfloor \;\, = \; \lfloor p, w \rfloor + \lfloor q, r \rfloor
\end{eqnarray*}

$+$ is clearly commutative.

Now let $\lfloor m, u \rfloor, \lfloor n, v \rfloor, \lfloor p, w \rfloor \in {\mathbb D}_+$. To see that $+$ is associative, notice that
\begin{eqnarray*}
\big(\lfloor m, u \rfloor + \lfloor n, v \rfloor\big) + \lfloor p, w \rfloor \!\! & = & \!\! \lfloor m \, 2^v + n \, 2^u , \; u + v \rfloor + \lfloor p, w \rfloor \\
 & = & \!\! \lfloor m \, 2^{v + w} + n \, 2^{u + w} + p \, 2^{u + v}, \; u + v + w \rfloor\\
 & = & \!\! \lfloor m, u \rfloor + \big(\lfloor n, v \rfloor + \lfloor p, w \rfloor \big)
\end{eqnarray*}

To see the last assertion, notice that
\[
\lfloor m, u \rfloor + \lfloor p, w \rfloor \, = \, \lfloor m \, 2^w + p \, 2^u, \; u + w \rfloor \, = \, \lfloor m \, 2^{v + w} + p \, 2^{u + v}, \; u + v + w \rfloor
\]

and
\[
\lfloor n, v \rfloor + \lfloor p, w \rfloor \, = \, \lfloor n \, 2^w + p \, 2^v, \; v + w \rfloor \, = \, \lfloor n \, 2^{u + w} + p \, 2^{u + v}, \; u + v + w \rfloor
\]

\eproof

As in the case of natural numbers, associativity of $+$ for positive dyadic rationals allows one to write multiple sums without brackets, i.e.\ for every $d, e, f \in {\mathbb D}_+$ we may write $d + e + f$ instead of $(d + e) + f$ or $d + (e + f)$, and similarly for sums of more than three terms.

\midvspace

\bprop
\label{prop dyadic diff}
Let $d, e \in {\mathbb D}_+$. If $d < e$, then there is $f \in {\mathbb D}_+$ such that $d + f = e$.
\eprop

\bproof
Let $\lfloor m, u \rfloor, \lfloor n, v \rfloor \in {\mathbb D}_+$ with $\lfloor m, u \rfloor < \lfloor n, v \rfloor$. There is $p \in \naturalnumbers$ such that $m \, 2^v + p = n \, 2^u$ by Lemma~\ref{lemm natural diff}. It follows that
\begin{eqnarray*}
\lfloor m, u \rfloor + \lfloor p, \; u + v \rfloor \!\! & = & \!\! \lfloor m \, 2^v, \; u + v \rfloor + \lfloor p, \; u + v \rfloor\\
 & = & \!\! \lfloor m \, 2^v + p, \; u + v \rfloor\\
 & = & \!\! \lfloor n \, 2^u, \; u + v \rfloor \; = \; \lfloor n, v \rfloor
\end{eqnarray*}

\eproof

\blede
\index{Multiplication!positive dyadic rational numbers}
\index{Product!positive dyadic rational numbers}
Let $\cdot$ be the binary function on ${\mathbb D}_+$ defined by $\lfloor m, u \rfloor \cdot \lfloor n, v \rfloor = \lfloor m \, n, \; u + v \rfloor$. This function is called {\bf multiplication on}~${\mathbb D}_+$. The expression $\lfloor m, u \rfloor \cdot \lfloor n, v \rfloor$ is called the {\bf product of $\lfloor m, u \rfloor$ and $\lfloor n, v \rfloor$}. For every $a, b \in {\mathbb D}_+$ we also write $a \, b$ for $a \cdot b$. $\cdot$~is commutative and associative. For every $d, e, f \in {\mathbb D}_+$ the following distributive law holds:
\[
(d + e) \cdot f \, = \, (d \cdot f) + (e \cdot f)
\]

Moreover, we have
\[
\left( d < e \right) \, \wedge \, \left( 0 < f \right) \quad \Longrightarrow \quad d \cdot f < e \cdot f
\]

We define that in the absence of brackets products are evaluated before sums.
\elede

\bproof
To see that $\cdot$ is well-defined, let $\lfloor m, u \rfloor = \lfloor p, w \rfloor \in {\mathbb D}_+$ and $\lfloor n, v \rfloor = \lfloor q, r \rfloor \in {\mathbb D}_+$. Then we have
\begin{eqnarray*}
\lfloor m, u \rfloor \cdot \lfloor n, v \rfloor \!\! & = & \!\! \lfloor m \, n, \; u + v \rfloor\\
 & = & \!\! \lfloor m \, 2^w \, n \, 2^r, \; u + v + w + r \rfloor\\
 & = & \!\! \lfloor p \, 2^u \, q \, 2^v, \; u + v + w + r \rfloor\\
 & = & \!\! \lfloor p \, q ,\; w + r \rfloor \; = \; \lfloor p, w \rfloor \cdot \lfloor q, r \rfloor
\end{eqnarray*}

$\cdot$ is clearly commutative.

Let $\lfloor m, u \rfloor, \lfloor n, v \rfloor, \lfloor p, w \rfloor \in {\mathbb D}_+$. To see that $\cdot$ is associative, notice that
\begin{eqnarray*}
\big(\lfloor m, u \rfloor \cdot \lfloor n, v \rfloor\big) \cdot \lfloor p, w \rfloor \!\! & = & \!\! \lfloor m \, n \, p, \; u + v + w \rfloor\\
 & = & \!\! \lfloor m, u \rfloor \cdot \big(\lfloor n, v \rfloor \cdot \lfloor p, w \rfloor\big)
\end{eqnarray*}

The distributive law is seen by the following calculation:
\begin{eqnarray*}
 \big(\lfloor m, u \rfloor + \lfloor n, v \rfloor \big) \cdot \lfloor p, w \rfloor \!\!
 & = & \!\! \lfloor m \, 2^v + n \, 2^u, \; u + v \rfloor \cdot \lfloor p, w \rfloor\\
 & = & \!\! \lfloor m \, p \, 2^v + n \, p \, 2^u, \; u + v + w \rfloor\\
 & = & \!\! \lfloor m \, p \, 2^{v + w} + n \, p \, 2^{u + w}, \; u + v + 2 w \rfloor\\ 
 & = & \!\! \lfloor m \, p, \; u + w \rfloor + \lfloor n \, p, \; v + w \rfloor\\
 & = & \!\! \lfloor m, u \rfloor \cdot \lfloor p, w \rfloor + \lfloor n, v \rfloor \cdot \lfloor p, w \rfloor
\end{eqnarray*}

Now assuming that $\lfloor m, u \rfloor < \lfloor n, v \rfloor$, we have $m \, 2^v < n \, 2^u$ by definition. Thus $m \, p \, 2^{v + w} < n \, p \, 2^{u + w}$, and therefore $\lfloor m \, p, \; u + w \rfloor < \lfloor n \, p, \; v + w \rfloor$. This shows the asserted implication.
\eproof

Again, associativity allows one to write multiple products without brackets.

\midvspace

\blemm
\label{lemm nat pos dyad embed}
Let $g : \naturalnumbers \longrightarrow {\mathbb D}_+$, $g(m) = \lfloor m, 0 \rfloor$. Then $g$ is injective. For every $m, n \in \naturalnumbers$ we have 
\benum
\item $m < n \;\; \Longleftrightarrow \;\; g(m) < g(n)$
\item $g(m + n) = g(m) + g(n)$
\item $g(m \cdot n) = g(m) \cdot g(n)$
\eenum

\elemm

\bproof
Exercise.
\eproof

The injection from~$\naturalnumbers$ to~${\mathbb D}_+$ in Lemma~\ref{lemm nat pos dyad embed} preserves the ordering in the sense of~$<$ (and that in the sense of~$\leq$ too) as well as the binary functions addition and multiplication. This justifies the usage of the same symbols $<$, $\leq$, $+$, and $\cdot$. Furthermore this allows the deliberate usage of mixed notations such as $\lfloor m, u \rfloor + n$ for $\lfloor m, u \rfloor + g(n)$, $n < \lfloor m, u \rfloor$ for $g(n) < \lfloor m, u \rfloor$, etc.\ where $m, n, u \in \naturalnumbers$. In each case such mixed notation is understood as shorthand notation for the full expression including the required injections. Moreover we may write $0$ for $g(0)$ and $1$ for~$g(1)$. Similarly, if $A \subset \naturalnumbers$ and $B \subset {\mathbb D}_+$, we may write $A \cap B$ instead of $g \left[ A \right] \cap B$ without ambiguity. Occasionally, given $d \in {\mathbb D}_+$, we may even write $d \in \naturalnumbers$ instead of $d \in g \left[ \naturalnumbers \right]$, or, given $A \subset {\mathbb D}_+$, we may write $A \subset \naturalnumbers$ instead of $A \subset g \left[ \naturalnumbers \right]$.

\midvspace

\bprop
\label{prop multipl dense}
Let $d, e, f \in {\mathbb D}_+$ with $d < f$ and $e \neq 0$. There is $g \in {\mathbb D}_+$ such that $d < e \, g < f$.
\eprop

\bproof
Let $\lfloor m, u \rfloor, \lfloor n, v \rfloor, \lfloor p, w \rfloor \in {\mathbb D}_+$ where $\lfloor n, v \rfloor < \lfloor p, w \rfloor$ and $m > 0$. We may choose $r \in \naturalnumbers$ such that $m < 2^r$. Then we have $\lfloor m, r \rfloor < 1$. Further we may choose $s \in \naturalnumbers$ such that $v + w < u + s$. We have $\lfloor m, u \rfloor \cdot \lfloor q, \; r + s \rfloor = \lfloor m, r \rfloor \cdot \lfloor q, \; u + s \rfloor$ for every $q \in \naturalnumbers$. We define $n_0 = n \, 2^w$ and $p_0 = p \, 2^v$. There exists $q \in \naturalnumbers$ such that $\lfloor n, v \rfloor = \lfloor n_0, \; v + w \rfloor < \lfloor m, r \rfloor \cdot \lfloor q, \; u + s \rfloor < \lfloor p_0, \; v + w \rfloor = \lfloor p, w \rfloor$.
\\

\hspace{0.05\textwidth}
\parbox{0.95\textwidth}
{[Let $q_0 \in \naturalnumbers$ be the maximum natural number such that $\lfloor m, r \rfloor \cdot \lfloor q_0, \; u + s \rfloor \leq \lfloor n_0, \; v + w \rfloor$. Then
\begin{eqnarray*}
\lfloor m, r \rfloor \cdot \lfloor q_0+1, \; u + s \rfloor \!\! & = & \!\! \lfloor m, r \rfloor \cdot \big( \lfloor q_0, \; u + s \rfloor + \lfloor 1, \; u + s \rfloor \big)\\
& < & \!\! \lfloor n_0, \; v + w \rfloor + \lfloor 1, \; v + w \rfloor\\
& \leq & \!\! \lfloor p_0, \; v + w \rfloor
\end{eqnarray*}]}
\eproof

\bprop
\label{prop dyadic prod inequ}
Let $d, e, f, g \in {\mathbb D}_+$. If $d < e$ and $f < g$, then we have $d \, g + e \, f < d \, f + e \, g$.
\eprop

\bproof
If the stated condition holds, then there are $a, b \in {\mathbb D}_+$ such that $d + a = e$ and $f + b = g$ by Proposition~\ref{prop dyadic diff}. Then we have
\begin{eqnarray*}
d \, g + e \, f \!\! & = & \!\! d \, (f + b) + (a + d) \, f\\
 & < & \!\! (d + a) \, (f + b) + d \, f\\
 & = & \!\! e \, g + d \, f
\end{eqnarray*}

\eproof

\bprop
\label{prop sum dyadic}
Let $a, b, c \in {\mathbb D}_+$ with $a < b + c$. Then there are $f, g \in {\mathbb D}_+$ such that
\[
f < b \, , \quad g < c \, , \quad a < f + g
\]

\eprop

\bproof
We may choose $d \in {\mathbb D}_+$, $d > 0$, such that $a + 2 d < b + c$ by Proposition~\ref{prop dyadic diff}. If $b < d$, then we may choose $f \in {\mathbb D}_+$ with $f < b$ by Lemma~\ref{lemm dyadic order dense}. If $b \geq d$, then there is $h$ with $h + d = b$, and we may choose $f \in {\mathbb D}_+$ such that $h < f < b$ by Lemma~\ref{lemm dyadic order dense}. In both cases we have $b < f + d$. In a similar way we may choose a number $g \in {\mathbb D}_+$ such that $c < g + d$. We then have $a + 2 d < f + g + 2 d$, and hence $a < f + g$.
\eproof

\bprop
\label{prop factor dyadic}
Let $a, b, c \in {\mathbb D}_+$ with $a < b \, c$. Then there are $f, g \in {\mathbb D}_+$ such that
\[
f < b \, , \quad g < c \, , \quad a < f \, g
\]

\eprop

\bproof
We may choose $d, e \in {\mathbb D}_+$ such that $a < d < e < b \, c$ by Lemma~\ref{lemm dyadic order dense}. There is $m \in \naturalnumbers$ such that $d + \lfloor 1, m \rfloor < e$ by Proposition~\ref{prop dyadic diff}. Further there is $k \in \naturalnumbers$ such that 
$(b + c) \lfloor 1, k \rfloor < \lfloor 1, 2m \rfloor$ by Proposition~\ref{prop multipl dense}. We define $n$ to be the maximum of $\left\{ m , k \right\}$. We may choose $f, g \in {\mathbb D}_+$ such that
\[
f < b \, , \quad b < f + \lfloor 1, n \rfloor \, , \quad g < c \, , \quad c < g + \lfloor 1, n \rfloor
\]

It follows that $(f + g) \, \lfloor 1, n \rfloor < (b + c) \, \lfloor 1, k \rfloor < \lfloor 1, 2m \rfloor$. Thus we have
\begin{eqnarray*}
b \, c \!\! & < & \!\! f \, g + (f + g) \, \lfloor 1, n \rfloor + \lfloor 1, 2 n \rfloor\\
 & < & \!\! f \, g + \lfloor 1, 2 m \rfloor + \lfloor 1, 2 m \rfloor \; = \; f \, g + \lfloor 1, m \rfloor
\end{eqnarray*}

\eproof

\section{Positive real numbers}
\label{positive real numbers}

In this Section we introduce the positive real numbers, its orderings, as well as addition, multiplication, and exponentiation on positive real numbers. We also show how to identify the positive dyadic rationals as a subset of the positive reals.

\midvspace

\blede
\label{lede pos real num}
\index{Positive real numbers}
\index{Numbers!positive real}
We define ${\mathbb D}_0 = \left\{ \; \left] -\infty, d \right[ \; : \, d \in {\mathbb D}_+ \right\}$ where the lower segments refer to the relation $<$ on ${\mathbb D}_+$, and
${\mathbb R}_+ = \big\{ \, {\textstyle \bigcup} {\mc A} \, : \, {\mc A} \subset {\mathbb D}_0, \, {\mc A} \neq \O \big\} \setminus \left\{ {\mathbb D}_+ \right\}$. The members of ${\mathbb R}_+$ are called {\bf positive real numbers}. We further define a total ordering in the sense of~"$\leq$" on ${\mathbb R}_+$ by
\[
\alpha \leq \beta \quad \Longleftrightarrow \quad \alpha \subset \beta
\]

Moreover, we define $<$ to be the total ordering in the sense of~"$<$" on~${\mathbb R}_+$ obtained from the ordering~$\leq$ by the method of Lemma~\ref{lemm both orderings}.
\elede

\bproof
It follows from Lemma~\ref{lede subset ordering} that $\leq$ is an ordering in the sense of~"$\leq$".

To see that $\leq$ is connective let $\alpha, \beta \in {\mathbb R}_+$ and assume that $\alpha \leq \beta$ does not hold. Then there is $d \in {\mathbb D}_+$ such that $d \in \alpha$, $d \notin \beta$. It follows that $e < d$ for every $e \in \beta$, and thus $\beta \subset \alpha$.
\eproof

\brema
Notice that $\lfloor 0, 0 \rfloor \in {\mathbb D}_+$, and $\O = \; \left] -\infty, \lfloor 0, 0 \rfloor \right[ \; \in {\mathbb D}_0 \subset {\mathbb R}_+$.
\erema

Remember that the usage of the symbol $-\infty$ in an interval denotes a lower segment but does generally not imply that the interval has no lower bound. On the contrary, the members of~${\mathbb D}_0$ all have a lower bound, viz.~$\O$.

\midvspace

\blemm
\label{lemm D0 dense}
${\mathbb D}_0$ is $<$-dense in~${\mathbb R}_+$. ${\mathbb R}_+$ is $<$-dense.
\elemm

\bproof
To see the first claim, let $\alpha, \beta \in {\mathbb R}_+$ with $\alpha < \beta$. Then there exists $d \in \beta \!\setminus\! \alpha$. Furthermore, there exists $e \in {\mathbb D}_+$ such that $d \in \;\left] -\infty, e \right[\; \subset \beta$. Thus we have $\alpha \leq \;\left] -\infty, d \right[\; < \;\left] -\infty, e \right[\; \leq \beta$. The claim follows by Lemma~\ref{lemm dyadic order dense}.

The second claim is a consequence of the first one.
\eproof

\blemm
\label{lemm least upper bound pos reals}
The ordered space $({\mathbb R}_+,<)$ has the least upper bound property. Specifically, if $A \subset {\mathbb R}_+$, $A \neq \O$, and $A$ has an upper bound, then \,\! $\sup A = \bigcup A$.
\elemm

\bproof
Assume the stated conditions and let $\alpha = \bigcup A$. We may choose $\beta \in {\mathbb R}_+$ such that $\gamma < \beta$ for every $\gamma \in A \!\setminus\! \left\{ \beta \right\}$. It follows that $\alpha \subset \beta$, and thus $\alpha \in {\mathbb R}_+$. Moreover, we have $\gamma \leq \alpha$ for every $\gamma \in A$. It is also clear that $\alpha$ is the least upper bound.
\eproof

\midvspace

\blede
\label{lede add mult pos reals}
\index{Addition!positive real numbers}
\index{Sum!positive real numbers}
\index{Multiplication!positive real numbers}
\index{Product!positive real numbers}
We define two binary functions $+$ (called {\bf addition}) and $\cdot$ (called {\bf multiplication}) on ${\mathbb R}_+$ by
\begin{eqnarray*}
\alpha + \beta \!\! & = & \!\! \bigcup \big\{ \left] -\infty, d + e \right[ \; : \, d \in D, \, e \in E \big\}\\[.2em]
\alpha \cdot \beta \!\! & = & \!\! \bigcup \big\{ \left] -\infty, d \cdot e \right[ \; : \, d \in D, \, e \in E \big\}
\end{eqnarray*}

where $D, E \subset {\mathbb D}_+$ such that $D, E \neq \O$, $\alpha = \bigcup \left\{ \; \left] -\infty, d \right[ \; : \, d \in D \right\}$ and $\beta = \bigcup \left\{ \; \left] -\infty, e \right[ \; : \, e \in E \right\}$. $\alpha + \beta$ and $\alpha \cdot \beta$ are called the {\bf sum} and the {\bf product of $\alpha$ and $\beta$}, respectively. We also write $\alpha \, \beta$ for $\alpha \cdot \beta$. Both functions are commutative and associative, and the distributive law
\[
(\alpha + \beta) \cdot \gamma \, = \, (\alpha \cdot \gamma) + (\beta \cdot \gamma)
\]

holds for $\alpha, \beta, \gamma \in {\mathbb R}_+$.

We define that in the absence of brackets products are evaluated before sums.
\elede

\bproof
We first show that $+$ and $\cdot$ are well-defined.

We may choose upper bounds $d_0, e_0 \in {\mathbb D}_+$ of $D$ and $E$, respectively. Then $d + e \leq d_0 + e_0$ and $d \, e \leq d_0 \, e_0$ for every $d \in D$, $e \in E$. Thus $(\alpha + \beta), (\alpha \cdot \beta) \neq {\mathbb D}_+$.

In order to see that the definitions of~$+$ and $\cdot$ do not depend on the choice of the index sets $D$ and~$E$, let $\alpha, \beta \in {\mathbb R}_+$ and, for $i \in \left\{ 1, 2 \right\}$, let $D_i, E_i \subset {\mathbb D}_+$ such that $D_i, E_i \neq \O$, $\alpha = \bigcup \left\{ \; \left] -\infty, d \right[ \; : \, d \in D_i \right\}$ and $\beta = \bigcup \left\{ \; \left] -\infty, e \right[ \; : \, e \in E_i \right\}$. Further let $d_1 \in D_1$ and $e_1 \in E_1$, and assume that not both $d_1 = 0$ and $d_2 = 0$.

To show the claim for~$+$, let $f \in {\mathbb D}_+$ with $f < d_1 + e_1$. First we consider the case $d_1 > 0$, $e_1 > 0$. We may choose $a, b \in {\mathbb D}_+$ such that $a < d_1$, $b < e_1$, and $f < a + b$ by Proposition~\ref{prop sum dyadic}. Since $a \in \alpha$, there is $d_2 \in D_2$ such that $a < d_2$. Similarly, since $b \in \beta$, there is $e_2 \in E_2$ such that $b < e_2$. It follows that $f < d_2 + e_2$. Second consider the case $d_1 = 0$, $e_1 > 0$. We may choose any $d_2 \in D_2$, and $b \in {\mathbb D}_+$ such that $f < b < e_1$ by Proposition~\ref{lemm dyadic order dense}, and $e_2 \in E_2$ such that $b < e_2$. It follows that $f < d_2 + e_2$. The case $d_1 > 0$, $e_1 = 0$ is handled similarly.

To show the claim for~$\cdot$, assume that $d_1 > 0$ and $e_1 > 0$, and let $f \in {\mathbb D}_+$ with $f < d_1 \, e_1$. We may choose $a, b \in {\mathbb D}_+$ such that $a < d_1$, $b < e_1$, and $f < a \, b$ by Proposition~\ref{prop factor dyadic}. There is $d_2 \in D_2$ such that $a < d_2$. Further there is $e_2 \in E_2$ such that $b < e_2$. It follows that $f < d_2 \, e_2$.

The commutativity and associativity of $+$ and $\cdot$ on ${\mathbb R}_+$ is a consequence of the respective properties of $+$ and $\cdot$ on~${\mathbb D}_+$.

In order to prove the distributive law we define the following sets:
\begin{center}
\begin{tabular}{ll}
$\alpha \, = \, {\displaystyle \bigcup} \, \big\{ \left] -\infty, d \right[ \; : \, d \in D \big\}$\,, \quad \quad &
$R \, = \, \left\{ d + e \, : \, d \in D, \, e \in E \right\}$\,,\\[.8em]
$\beta \, = \, {\displaystyle \bigcup} \, \big\{ \left] -\infty, e \right[ \; : \, e \in E \big\}$\,, \quad \quad &
$S \, = \, \left\{ d \, f \, : \, d \in D, \, f \in F \right\}$\,,\\[.8em]
$\gamma \, = \, {\displaystyle \bigcup} \, \big\{ \left] -\infty, f \right[ \; : \, f \in F \big\}$\,, \quad \quad &
$T \, = \, \left\{ e \, f \, : \, e \in E, \, f \in F \right\}$
\end{tabular}
\end{center}\

The distributive law then follows from the following calculation:
\begin{eqnarray*}
\left( \alpha + \beta \right) \, \gamma \!\!\! & = & \!\! 
\bigcup \big\{ \left] -\infty, d + e \right[ \; : \, d \in D, \, e \in E \big\} \, \cdot \;
\bigcup \big\{ \left] -\infty, f \right[ \; : \, f \in F \big\}\\[.2em]
 & = & \!\! \bigcup \big\{ \left] -\infty, r \right[ \; : \, r \in R \big\} \, \cdot \;
\bigcup \big\{ \left] -\infty, f \right[ \; : \, f \in F \big\}\\[.2em]
 & = & \!\! \bigcup \big\{ \left] -\infty, r \, f \right[ \; : \, r \in R, \, f \in F \big\}\\[.2em]
 & = & \!\! \bigcup \big\{ \left] -\infty, (d + e) \, f \right[ \; : \, d \in D, \, e \in E, \, f \in F \big\}\\[.2em]
 & = & \!\! \bigcup \big\{ \left] -\infty, d \, f + e \, f \right[ \; : \, d \in D, \, e \in E, \, f \in F \big\}\\[.2em]
 & = & \!\! \bigcup \big\{ \left] -\infty, d \, f + e \, h \right[ \; : \, d \in D, \, e \in E, \, f \in F, \, h \in F \big\}\\[.2em]
 & = & \!\! \bigcup \big\{ \left] -\infty, s + t \right[ \; : \, s \in S, \, t \in T \big\}\\[.2em]
 & = & \!\! \bigcup \big\{ \left] -\infty, s \right[ \; : \, s \in S \big\} \, + \, \bigcup \big\{ \; \left] -\infty, t \right[ \; : \, t \in T \big\}\\[.2em]
& = & \!\! \alpha \, \gamma + \beta \, \gamma
\end{eqnarray*}

\eproof

\midvspace

\blemm
\label{lemm injection dyad reals}
Let $g : {\mathbb D}_+ \longrightarrow {\mathbb R}_+$, $g(d) = \; \left] -\infty, d \right[ \;$. Then $g$ is injective. For every $d, e \in {\mathbb D}_+$ we have
\benum
\item $g(\lfloor 0, 0 \rfloor) = \O$
\item $d < e \;\; \Longleftrightarrow \;\; g(d) < g(e)$
\item $g(d + e) = g(d) + g(e)$
\item $g(d \cdot e) = g(d) \cdot g(e)$
\eenum

\elemm

\bproof
Exercise.
\eproof

Regarding the injection from the positive dyadic rationals to the positive reals in Lemma and Definition~\ref{lemm injection dyad reals}, the same comments apply as regarding the injection from the natural numbers to the positive dyadic rationals defined in Lemma~\ref{lemm nat pos dyad embed}. That is, it preserves the orderings $<$ and~$\leq$ as well as the binary functions $+$ and~$\cdot$. This, again, justifies the usage of the same symbols and allows one to write mixed expressions of positive dyadic rationals and positive reals, but also of natural numbers and positive reals. In the latter case the notation of both injections is then suppressed. For instance, $\alpha + d$, $m \cdot \alpha$, and $\alpha \leq m$ are valid expressions, where $m \in \naturalnumbers$, $d \in {\mathbb D}_+$, and $\alpha \in {\mathbb R}_+$. Occasionally, we may even write $\alpha \in \naturalnumbers$, or $A \subset \naturalnumbers$ although actually $A \subset {\mathbb R}_+$.

\midvspace

\bprop
\label{prop pos real estim}
For every $\delta, \varepsilon \in {\mathbb R}_+$ and $r \in {\mathbb D}_+$, the inequality $\delta \leq r$ implies $\delta + \varepsilon \leq r + \varepsilon$ and $\delta \, \varepsilon \leq r \, \varepsilon$.
\eprop

\bproof
Let $G, H \subset {\mathbb D}_+$, $G, H \neq \O$, such that
\[
\delta \, = \, \bigcup \big\{ \left] -\infty, g \right[ \; : \, g \in G \big\} \, , \quad \quad
\varepsilon \, = \, \bigcup \big\{ \left] -\infty, h \right[ \; : \, h \in H \big\}
\]

If $\delta \leq r$, then $g \leq r$ for every $g \in G$. It follows that
\begin{eqnarray*}
\delta + \varepsilon \!\!\!
& = & \!\! \bigcup \big\{ \left] -\infty, g + h \right[ \; : \, g \in G, \, h \in H \big\}\\[.2em]
& \leq & \!\! \bigcup \big\{ \left] -\infty, r + h \right[ \; : \, h \in H \big\} \, = \; r + \varepsilon
\end{eqnarray*}

and
\begin{eqnarray*}
\delta \, \varepsilon \!\!\!
& = & \!\! \bigcup \big\{ \left] -\infty, g \, h \right[ \; : \, g \in G, \, h \in H \big\}\\[.2em]
& \leq & \!\! \bigcup \big\{ \left] -\infty, r \, h \right[ \; : \, h \in H \big\} \, = \; r \, \varepsilon
\end{eqnarray*}

\eproof

\blemm
\label{lemm inequalities reals}
Given $\alpha, \beta, \gamma \in {\mathbb R}_+$, the following implications hold:
\begin{eqnarray*}
\alpha < \beta \; & \Longrightarrow & \; \alpha + \gamma < \beta + \gamma\\
\left( \alpha < \beta \right) \, \wedge \, \left( 0 < \gamma \right) \; & \Longrightarrow & \; \alpha \, \gamma < \beta \, \gamma
\end{eqnarray*}

\elemm

\bproof
Let $D, E, F \subset {\mathbb D}_+$ such that $D, E, F \neq \O$ and

\begin{center}
\begin{tabular}{ll}
$\alpha \, = \, {\displaystyle \bigcup} \, \big\{ \left] -\infty, d \right[ \; : \, d \in D \big\}$\,, \quad \quad &
$\beta \, = \, {\displaystyle \bigcup} \, \big\{ \left] -\infty, e \right[ \; : \, e \in E \big\}$\,, \\[.8em]
$\gamma \, = \, {\displaystyle \bigcup} \, \big\{ \left] -\infty, f \right[ \; : \, f \in F \big\}$ &
\end{tabular}
\end{center}

Assume that $\alpha < \beta$. We may choose $a, b \in {\mathbb D}_+$ such that $\alpha < a < b < \beta$ by Proposition~\ref{lemm D0 dense}. Then we have $d < a$ for every $d \in D$. Further there exists $e_0 \in E$ such that $b < e_0$. 

Notice that the first implication is clear for $\gamma = 0$. Now assume that $0 < \gamma$. We may choose $m \in \naturalnumbers$ such that $a + \lfloor 1, m \rfloor < b$ by Proposition~\ref{prop dyadic diff}. Let $k$ be the maximum of~$\left\{ n \in \naturalnumbers \, : \, \lfloor n, m \rfloor < \gamma \right\}$.
\\

\hspace{0.05\textwidth}
\parbox{0.95\textwidth}
{[There exists $\lfloor p, v \rfloor \in {\mathbb D}_+$ such that $\gamma < \lfloor p, v \rfloor$. We have
\[
\lfloor p, v \rfloor = \lfloor p \, 2^m, v + m \rfloor \leq \lfloor p \, 2^m, m \rfloor
\]

Therefore the considered set is finite and non-empty since $0 < \gamma$. Thus it has a maximum by Lemma~\ref{lemm tot ord subset min max}.]}
\\

There is $c \in F$ such that $\lfloor k, m \rfloor < c$. We have
\begin{center}
\begin{tabular}{l}
$\alpha + \gamma \; \leq \; a + \gamma \; \leq \; a + \lfloor \sigma(k), m \rfloor \; = \; a + \lfloor k, m \rfloor + \lfloor 1, m \rfloor \; <  \; b + \lfloor k, m \rfloor$\\[.8em]
\quad $ < \; b + c \; \leq \; {\displaystyle \bigcup} \, \big\{ \left] -\infty, b + f \right[ \; : \, f \in F \big\}$\\[.8em]
\quad $ \leq \; {\displaystyle \bigcup} \, \big\{ \left] -\infty, e + f \right[ \; : \, e \in E, \, f \in F \big\} \, = \, \beta + \gamma$
\end{tabular}
\end{center}\

where the first and second inequality follow by Proposition~\ref{prop pos real estim}.

To prove the second implication, let $h \in {\mathbb D}_+$ such that $a + h = b$ by Proposition~\ref{prop dyadic diff}, and assume that $0 < \gamma$. We have $0 < h$, and thus $0 < h \, \gamma$. It follows that
\begin{center}
\begin{tabular}{l}
$\alpha \, \gamma \; \leq \; a \, \gamma \; < \; a \, \gamma + h \, \gamma \; = \; b \, \gamma \; = \; {\displaystyle \bigcup} \, \big\{ \left] -\infty, b \, f \right[ \; : \, f \in F \big\}$\\[.8em]
\quad $ \leq \; {\displaystyle \bigcup} \, \big\{ \left] -\infty, e \, f \right[ \; : \, e \in E, \, f \in F \big\} \, = \, \beta \, \gamma$
\end{tabular}
\end{center}\

The first inequality follows by Proposition~\ref{prop pos real estim}. The second inequality is a consequence of the inequality for sums.
\eproof

\blemm
Let $\alpha, \beta \in {\mathbb R}_+$. We have $\alpha + \beta = \left\{ a + b \, : \, a \in \alpha,\; b \in \beta \right\}$.
\elemm

\bproof
Let $D, E \subset {\mathbb D}_+$ with $D, E \neq \O$ such that
\[
\alpha \, = \, \bigcup \big\{ \left] -\infty, d \right[ \; : \, d \in D \big\} \, , \quad \quad
\beta \, = \, \bigcup \big\{ \left] -\infty, e \right[ \; : \, e \in E \big\}
\]

Now let $c \in \alpha + \beta$. There are $d \in D$ and $e \in E$ such that $c < d + e$. If $c < d$, then we have $c \in \alpha$. If $c < e$, then we have $c \in \beta$. If $c \geq d$ and $c \geq e$, then we define $g \in {\mathbb D}_+$ such that $c + g = d + e$ by Proposition~\ref{prop dyadic diff}. It follows that $0 < g$, $g \leq d \leq 2d$, and $g \leq e \leq 2e$. We further define $a \in {\mathbb D}_+$ such that $g + 2a = 2d$, as well as $b \in {\mathbb D}_+$ such that $g + 2b = 2e$. Thus we obtain $a \in \alpha$, $b \in \beta$, and $a + b = c$.

The converse is clear.
\eproof

\bprop
\label{prop pos real diff}
Let $\alpha, \beta \in {\mathbb R}_+$. If $\alpha < \beta$, then there is $\gamma \in {\mathbb R}_+$ such that $\alpha + \gamma = \beta$.
\eprop

\bproof
Assume the condition. We define $\gamma = \sup \left\{ \delta \in {\mathbb R}_+ \, : \, \alpha + \delta < \beta \right\}$. 
\\

\hspace{0.05\textwidth}
\parbox{0.95\textwidth}
{[The supremum is well-defined by Lemmas~\ref{lemm inequalities reals} and~\ref{lemm least upper bound pos reals} since $\alpha + \beta \geq \beta$.]}
\\

We may choose $A \subset {\mathbb D}_+$ with $A \neq \O$ and, for every $\delta \in {\mathbb R}_+$, $D_{\delta} \subset {\mathbb D}_+$ with $D_{\delta} \neq \O$ such that
\[
\alpha \, = \, \bigcup \big\{ \left] -\infty, a \right[ \; : \, a \in A \big\}, \quad \quad
\delta \, = \, \bigcup \big\{ \left] -\infty, d \right[ \; : \, d \in D_{\delta} \big\}
\]

It follows that
\[
\gamma \, = \, \bigcup \big\{ \left] -\infty, d \right[ \; : \, d \in D_{\delta}, \; \alpha + \delta < \beta \big\} \; = \; \bigcup \big\{ \left] -\infty, d \right[ \; : \, d \in D \big\}
\]

where
\[
D \, = \, \bigcup \big\{ D_{\delta} \, : \, \alpha + \delta < \beta \big\}
\]

Thus we obtain
\begin{eqnarray*}
\alpha + \gamma \!\!\! & = &  \!\!\! \bigcup \big\{ \left] -\infty, a + d \right[ \; : \, a \in A, \; d \in D \big\}\\[.2em]
 & = & \!\!\! \bigcup \big\{ \left] -\infty, a + d \right[ \; : \, a \in A, \; d \in D_{\delta}, \; \alpha + \delta < \beta \big\} \, \leq \; \beta
\end{eqnarray*}

Now assume that $\alpha + \gamma < \beta$. Then there are $d, e \in {\mathbb D}_+$ such that $\alpha + \gamma < d < e < \beta$ by Lemma~\ref{lemm D0 dense}. Further there is $f \in {\mathbb D}_+$ with $f > 0$ such that $d + f = e$ by Proposition~\ref{prop dyadic diff}. It follows that $\alpha + \gamma + f < \beta$, which is a contradiction.
\eproof

\bprop
\label{prop prod pos reals}
Let $\alpha, \beta \in {\mathbb R}_+ \!\!\setminus\! \left\{ 0 \right\}$. There is $\gamma \in {\mathbb R}_+$ such that $0 < \alpha \, \gamma < \beta$.
\eprop

\bproof
We may choose $a, b \in {\mathbb D}_+$ such that $\alpha < a$ and $0 < b < \beta$ by Lemma~\ref{lemm D0 dense}. There is $c \in {\mathbb D}_+$ such that $0 < a \, c < b$ by Proposition~\ref{prop multipl dense}. It follows that $0 < \alpha \, c < a \, c < \beta$.
\eproof

\bprop
\label{prop prod inequ}
Let $\alpha, \beta, \gamma, \delta \in {\mathbb R}_+$. If $\alpha < \beta$ and $\gamma < \delta$, then $\alpha \, \delta + \beta \, \gamma < \alpha \, \gamma + \beta \, \delta$.
\eprop

\bproof
This follows by the corresponding result for positive dyadic rationals, see Lemma~\ref{prop dyadic prod inequ}.
\\

\hspace{0.05\textwidth}
\parbox{0.95\textwidth}
{[We may choose $A, B, C, D \subset {\mathbb D}_+$ such that
\begin{center}
\begin{tabular}{ll}
$\alpha \, = \; {\displaystyle \bigcup} \, \big\{ \left] -\infty, a \right[ \; : \, a \in A \big\}$\,, \quad \quad &
$\beta \, = \; {\displaystyle \bigcup} \, \big\{ \left] -\infty, b \right[ \; : \, b \in B \big\}$\,, \\[.8em]
$\gamma \, = \; {\displaystyle \bigcup} \, \big\{ \left] -\infty, c \right[ \; : \, c \in C \big\}$\,, &
$\delta \, = \; {\displaystyle \bigcup} \, \big\{ \left] -\infty, d \right[ \; : \, d \in D \big\}$
\end{tabular}
\end{center}\

By definition of addition and multiplication we have
\begin{eqnarray*}
\alpha \, \delta + \beta \, \gamma \!\!\! & = & \!\!\! \bigcup \big\{ \left] -\infty, a \, d + b \, c \right[ \; : \, a \in A,\; b \in B,\; c \in C,\; d \in D \big\} \, ,\\[.5em]
\alpha \, \gamma + \beta \, \delta \!\!\! & = & \!\!\! \bigcup \big\{ \left] -\infty, a \, c + b \, d \right[ \; : \, a \in A,\; b \in B,\; c \in C,\; d \in D \big\}
\end{eqnarray*}

Now let $e \in \alpha \, \delta + \beta \, \gamma$. We may choose $a \in A$, $b_0 \in B$, $c \in C$, and $d_0 \in D$ such that $e < a \, d_0 + b_0 \, c$. Further, there is $b \in B$ such that $a < b$ and $b_0 \leq b$. Similarly, there is $d \in D$ such that $c < d$ and $d_0 \leq d$. It follows that $e < a \, d + b \, c$. We have $a \, d + b \, c < a \, c + b \, d$ by Lemma~\ref{prop dyadic prod inequ}. Hence $e \in \alpha \, \gamma + \beta \, \delta$.]}
\eproof

\midvspace

\bprop
\label{prop greater prod est}
Let $\alpha, \beta, \gamma \in {\mathbb R}_+$ with $\alpha \, \beta < \gamma$. There are $a, b \in {\mathbb D}_+$ such that $\alpha < a$, $\beta < b$, and $a \, b < \gamma$.
\eprop

\bproof
We may choose $c, d \in {\mathbb D}_+$ such that $\alpha \, \beta < c < d < \gamma$ by Lemma~\ref{lemm D0 dense}. There is $m \in \naturalnumbers$ with $0 < m$ such that $c + \lfloor 1, m \rfloor < d$ by Proposition~\ref{prop dyadic diff}. We may choose \mbox{$k \in \naturalnumbers$} such that $\left( \alpha + \beta\ \right) \cdot \lfloor 1, k \rfloor < \lfloor 1, m + 1 \rfloor$ by Proposition~\ref{prop multipl dense}. Let $n = \sup \left\{ m, k \right\}$. We may choose $a, b \in {\mathbb D}_+$ such that $\alpha < a < \alpha + \lfloor 1, n \rfloor$ and $\beta < b < \beta + \lfloor 1, n \rfloor$ by Lemma~\ref{lemm D0 dense}. Further, notice that $m + 1 \leq 2m$. 
\\

\hspace{0.05\textwidth}
\parbox{0.95\textwidth}
{[The claim is clear for $m = 1$. Assume it holds for some $m \in \naturalnumbers$ with $m \geq 1$. Then we have $\sigma(m) + 1 \leq 2m + 1 < 2 \, \sigma(m)$.]}
\\

It follows that
\begin{eqnarray*}
a \, b \!\!\! & < & \!\!\! \alpha \, \beta + \left( \alpha + \beta\ \right) \cdot \lfloor 1, n \rfloor + \lfloor 1, 2 n \rfloor\\[.2em]
 & < & \!\!\! \alpha \, \beta + \lfloor 1, m + 1 \rfloor + \lfloor 1, 2 m \rfloor\\[.2em]
 & \leq & \!\!\! \alpha \, \beta + \lfloor 1, m + 1 \rfloor + \lfloor 1, m + 1 \rfloor \, = \; \alpha \, \beta + \lfloor 1, m \rfloor
\end{eqnarray*}

\eproof

\btheo
\label{theo pos reals group}
The triple $\left({\mathbb R}_+ \!\!\setminus\! \left\{ 0 \right\},\, \cdot \, , 1 \right)$ is an Abelian group. Let $\alpha, \beta \in {\mathbb R}_+ \!\!\setminus\! \left\{ 0 \right\}$. The inverse of $\alpha$ with respect to multiplication is denoted by $\alpha^{-1}$, $\left( 1 / \alpha \right)$ or~${\displaystyle \left(\frac{1}{\alpha}\right)}$. We also write $\left( \alpha / \beta \right)$ or ${\displaystyle \left( \frac{\alpha}{\beta} \right)}$ for $\alpha \cdot (1 / \beta)$. In the absence of brackets the superscript~"$-1$" is evaluated before sums and products.
\etheo

\bproof
Multiplication is associative and commutative by Lemma and Definition~\ref{lede add mult pos reals}.

Moreover, we have $\alpha \cdot 1 = \alpha$ for every $\alpha \in {\mathbb R}_+$ by definition.

Now let $\alpha \in {\mathbb R}_+ \!\! \setminus \! \left\{ 0 \right\}$. We define $A = \left\{ \beta \in {\mathbb R}_+ \!\! \setminus\! \left\{ 0 \right\} \, : \, \beta \, \alpha < 1 \right\}$. It follows that $A \neq \O$ by Proposition~\ref{prop prod pos reals}. Moreover, $A$ has an upper bound by Proposition~\ref{prop multipl dense}. Hence the supremum of~$A$ is well-defined by Lemma~\ref{lemm least upper bound pos reals}, and we have \,\! $\sup A = \bigcup A$. We define $1/\alpha = \sup A$. We may choose $D \subset {\mathbb D}_+$ with $D \neq \O$ and, for every $\beta \in {\mathbb R}_+ \!\!\setminus\! \left\{ 0 \right\}$, $E_{\beta} \subset {\mathbb D}_+$ with $E_{\beta} \neq \O$ such that
\[
\alpha \, = \, \bigcup \big\{ \left] -\infty, d \right[ \; : \, d \in D \big\}, \quad \quad
\beta \, = \, \bigcup \big\{ \left] -\infty, e \right[ \; : \, e \in E_{\beta} \big\}
\]

It follows that
\[
1/\alpha \, = \; \bigcup \big\{ \left] -\infty, e \right[ \; : \, e \in E_{\beta},\; \beta \in {\mathbb R}_+ \!\! \setminus \! \left\{ 0 \right\},\; \beta \, \alpha < 1 \big\}
\]

and hence
\[
\alpha \cdot \left( 1 / \alpha \right) \, = \; \bigcup \big\{ \left] -\infty, d \, e \right[ \; : \, d \in D,\; e \in E_{\beta},\; \beta \in {\mathbb R}_+ \!\! \setminus \! \left\{ 0 \right\},\; \beta \, \alpha < 1 \big\} \; \leq \; 1
\]

Now assume that $\alpha \cdot \left( 1 / \alpha \right) < 1$. Then there is $f \in {\mathbb D}_+$ such that $1/\alpha < f$ and $\alpha \, f < 1$ by Proposition~\ref{prop greater prod est}, which is a contradiction.

Thus we have $\alpha \cdot \left( 1 / \alpha \right) = 1$.
\eproof

\bcoro
\label{coro inequality inv}
Let $\alpha, \beta \in {\mathbb R}_+$ with $0 < \alpha < \beta$. Then we have $\beta^{-1} < \alpha^{-1}$.
\ecoro

\bproof
The inequality $\alpha < \beta$ implies $1 < \beta \, \alpha^{-1}$. The claim follows.
\eproof

\brema
\label{rema inv prod}
Let $\alpha, \beta \in {\mathbb R}_+ \!\! \setminus \! \left\{ 0 \right\}$. We have $\left( \alpha \, \beta \right)^{-1} = \alpha^{-1} \, \beta^{-1}$.
\erema

\bcoro
\label{coro prod greater}
Let $\alpha, \beta \in {\mathbb R}_+ \!\! \setminus \! \left\{ 0 \right\}$. There is $\gamma \in {\mathbb R}_+$ such that $\beta < \alpha \, \gamma$.
\ecoro

\bproof
There is $\delta \in {\mathbb R}_+$ such that $\alpha^{-1} \, \delta < \beta^{-1}$ by Proposition~\ref{prop prod pos reals}. It follows that $\beta < \alpha \, \delta^{-1}$ by Corollary~\ref{coro inequality inv} and Remark~\ref{rema inv prod}.
\eproof

We continue by defining exponentiation on the positive reals where the exponent is a natural number.

\midvspace

\blede
\label{lede exp pos reals}
\index{Exponentiation!positive real numbers}
We define a function $h : {\mathbb R}_+ \!\times \naturalnumbers \longrightarrow {\mathbb R}_+$ recursively by
\benum
\item $h(\alpha, 0) = 1$
\item $h(\alpha, \sigma(m)) = h(\alpha, m) \cdot \alpha$ 
\eenum

for every $\alpha \in {\mathbb R}_+$ and every $m \in \naturalnumbers$. This function is called {\bf exponentiation on}~${\mathbb R}_+$. We also write $\alpha^m$ for~$h(\alpha, m)$ and call $\alpha$ the {\bf base} and $m$ the {\bf exponent} or {\bf power}. In the absence of brackets we define the following priorities:
\[
\alpha^{m+n} = \alpha^{(m+n)}\,, \quad \alpha^{m \cdot n} = \alpha^{(m \cdot n)}\,, \quad \alpha + \beta^m = \alpha + \left( \beta^m \right)\,, \quad \alpha \, \beta^m = \alpha \, \left( \beta^m \right)
\] 

We have, for every $\alpha, \beta \in {\mathbb R}_+$ and $m, n \in \naturalnumbers$,
\[
\alpha^{m+n} = \alpha^m \, \alpha^n, \quad \quad \left( \alpha^m \right)^n = \alpha^{m \cdot n}, \quad \quad (\alpha \, \beta)^m = \alpha^m \, \beta^m
\]

and the implications
\begin{eqnarray*}
\left( \alpha < \beta \right) \; \wedge \; \left( 0 < m \right) \; & \Longrightarrow & \; \alpha^m < \beta^m\\
\left( 0 < \alpha < 1 \right) \; \wedge \; \left( m < n \right) \; & \Longrightarrow & \; \alpha^n < \alpha^m\\
\left( 1 < \alpha \right) \; \wedge \; \left( m < n \right) \; & \Longrightarrow & \; \alpha^m < \alpha^n
\end{eqnarray*}

Furthermore we define $\alpha^{-m} = \left(\alpha^{-1}\right)^m$ for $\alpha \in {\mathbb R}_+ \!\! \setminus \! \left\{ 0 \right\}$ and $m \in \naturalnumbers$. 

Given $\alpha \in {\mathbb R}_+ \!\! \setminus \! \left\{ 0 \right\}$ and $m, n \in \naturalnumbers$, we have $\left(\alpha^m\right)^{-1} = \alpha^{-m}$. If $m \leq n$, then
$\alpha^p = \alpha^n \, \alpha^{-m}$ where $p \in \naturalnumbers$ such that $m + p = n$. If $m > n$, then $\alpha^{-p} = \alpha^n \, \alpha^{-m}$ where $p \in \naturalnumbers$ such that $n + p = m$.

\elede

\bproof
The existence and uniqueness of the function follow by Theorem~\ref{theo recursive def}.

The three equations follow by the Induction principle.

We now show the three implications, again by means of the Induction principle.

To see the first implication, assume that $\alpha < \beta$. The implication clearly holds for $m = 1$. Now assume that it holds for some $m \in \naturalnumbers$ with $0 < m$. Then we have $\alpha^{\sigma(m)} = \alpha^m \, \alpha < \beta^m \, \beta = \beta^{\sigma(m)}$.

To see the second implication, let $m \in \naturalnumbers$ and $0 < \alpha < 1$. We have $\alpha^{\sigma(m)} = \alpha^m \, \alpha < \alpha^m$.
\\

\hspace{0.05\textwidth}
\parbox{0.95\textwidth}
{[The second inequality follows by Lemma~\ref{lemm inequalities reals} since $0 < \alpha^m$, which in turn is proven by the Induction principle.]}
\\

Now assume that the implication holds for some $n \in \naturalnumbers$ with $n \geq \sigma(m)$. It follows that $\alpha^{\sigma(n)} = \alpha^n \, \alpha < \alpha^m \, \alpha < \alpha^m$.

To see the third implication, let $m \in \naturalnumbers$ and $1 < \alpha$. We have $\alpha^m < \alpha^m \, \alpha < \alpha^{\sigma(m)}$. Now assume that the implication holds for some $n \in \naturalnumbers$ where $n \geq \sigma(m)$. It follows that $\alpha^m < \alpha^m \, \alpha < \alpha^n \, \alpha < \alpha^{\sigma(n)}$.

The equation $\left(\alpha^m\right)^{-1} = \alpha^{-m}$ clearly holds for every $\alpha \in {\mathbb R}_+ \!\! \setminus \! \left\{ 0 \right\}$ and $m = 0$. Now assume that it holds for every $\alpha \in {\mathbb R}_+ \!\! \setminus \! \left\{ 0 \right\}$ and some $m \in \naturalnumbers$. We have
\[
\big(\alpha^{\sigma(m)}\big)^{-1} = \left(\alpha^m \, \alpha\right)^{-1} = \left(\alpha^m\right)^{-1} \, \alpha^{-1} = \left(\alpha^{-1}\right)^m \, \alpha^{-1} = \left(\alpha^{-1}\right)^{\sigma(m)}
\]

Finally let $\alpha \in {\mathbb R}_+ \!\! \setminus \! \left\{ 0 \right\}$ and $m, n \in \naturalnumbers$. If $m \leq n$, then we have
\[
\alpha^n \, \alpha^{-m} = \alpha^{m + p} \, \alpha^{-m} = \alpha^m \, \alpha^p \, \alpha^{-m} = \alpha^p
\]

If $m > n$, then we have
\[
\alpha^p \, \alpha^n \, \alpha^{-m} = \alpha^{p + n} \, \alpha^{-m} = 1
\]

Thus $\alpha^n \, \alpha^{-m}$ is the inverse of $\alpha^p$.
\eproof

\blemm
\label{lemm exp embed}
For every $\lfloor m, u \rfloor \in {\mathbb D}_+$ and $n \in \naturalnumbers$, we have $\left(g \left( \lfloor m, u \rfloor \right) \right)^n = g \left(\lfloor m^n, n \, u \rfloor \right)$ where $g$ is defined in Lemma and Definition~\ref{lede add mult pos reals}.
\elemm

\bproof
This follows by the Induction principle.
\eproof

\bcoro
\label{coro nat reals exp embed}
Let $g : \naturalnumbers \longrightarrow {\mathbb R}_+$, $g(m) = \, \left] -\infty, \lfloor m, 0 \rfloor \right[ \,$. Then $g$ is injective, and we have for every $m, n \in \naturalnumbers$:
\[
g(m^n) = g(m)^n
\]

\ecoro

\bproof
The map $g$ is injective as it is the composition of two injections by Lemma~\ref{lemm nat pos dyad embed} and Lemma and Definition~\ref{lede add mult pos reals}.

The equation follows by Lemma~\ref{lemm exp embed}.
\eproof

Corollary~\ref{coro nat reals exp embed} shows that exponentiation on~$\naturalnumbers$ as defined in Lemma and Definition~\ref{lede exponentiation natural numbers} and exponentiation on~${\mathbb R}_+$ as defined in Lemma and Definition~\ref{lede exp pos reals} are in agreement with the injection from~$\naturalnumbers$ to~${\mathbb R}_+$. Therefore we may also use mixed notation even when exponentiation occurs.

\midvspace

\bcoro
Using the notation of negative exponents in Lemma and Definition~\ref{lede exp pos reals}, we have $\lfloor m, u \rfloor = m \, 2^{-u}$ for $m, u \in \naturalnumbers$.
\ecoro

\bproof
We have $\lfloor 1, u \rfloor \cdot 2^u = \lfloor 2^u, u \rfloor = 1$. Therefore $\lfloor 1, u \rfloor$ is the inverse of~$2^u$, i.e.\ $\lfloor 1, u \rfloor = 2^{-u}$.
\eproof

\blemm
\label{lemm dyadic compl dense}
${\mathbb R}_+ \!\!\setminus\! {\mathbb D}_0$ is $<$-dense in~${\mathbb R}_+$.
\elemm

\bproof
We first show that $1/3 \notin {\mathbb D}_+$. Assume there are $m, u \in \naturalnumbers$ such that $3 \cdot m \, 2^{-u} = 1$. It follows that $3m = 2^u$. This is clearly false for every $m \in \naturalnumbers$ and $u = 0$. Assume it is false for every $m \in \naturalnumbers$ and some $u \in \naturalnumbers$. If $m$ is even, then there is $n \in \naturalnumbers$ such that $m = 2n$, and therefore $3m = 2^{u+1}$ implies $3n = 2^u$, which is a contradiction. If $m$ is odd, then there is $n \in \naturalnumbers$ such that $m = 2n+1$, and thus $3m = 2^{u+1}$ implies $3 \cdot 2n + 3 = 2^{u+1}$. The left hand side of the last equation is odd whereas the right hand side is even, which is again a contradiction.

Now let $\alpha, \beta \in {\mathbb R}_+$ with $0 < \alpha < \beta$. There are $a, b \in {\mathbb D}_+$ such that $\alpha < a < b < \beta$ by Lemma~\ref{lemm D0 dense}. We may choose $m, u \in \naturalnumbers$ such that $a = \lfloor m, u \rfloor$ and $\lfloor m + 1 , u \rfloor < b$. We define $\gamma = \left( m + 1/3 \right) \, 2^{-u}$. Assume that $\gamma = n \, 2^{-v}$ for some $n, v \in \naturalnumbers$. It follows that $m + 1/3 = n \, 2^u \, 2^{-v}$, which is a contradiction to the first part of the proof.
\eproof

\section{Real numbers}
\label{real numbers}

In this Section it remains to construct the full number systems, i.e.\ those containing positive and negative numbers. Since natural numbers and positive dyadic rationals can be identified with a subset of the positive reals as shown above, it is enough to construct the system of positive and negative real numbers, its orderings, as well as addition and multiplication on the reals.

\midvspace

\blede
\label{lede definition reals}
\index{Real numbers}
\index{Standard ordering in the sense of~$<$}
\index{Standard ordering in the sense of~$\leq$}
Let $P$ be the equivalence relation on ${\mathbb R}_+^2$ defined by
\[
\big((\alpha, \beta), (\gamma, \delta)\big) \in P \quad \Longleftrightarrow \quad \alpha + \delta = \gamma + \beta
\]

and ${\mathbb R} = {\mathbb R}_+ / P$. The equivalence classes are called {\bf real numbers}. For every $\alpha, \beta \in {\mathbb R}_+$, the equivalence class of $(\alpha, \beta)$ is denoted by~$\langle \alpha, \beta \rangle$. We define a total ordering in the sense of~"$<$" on ${\mathbb R}$ by
\[
\langle \alpha, \beta \rangle < \langle \gamma, \delta \rangle \quad \Longleftrightarrow \quad \alpha + \delta < \gamma + \beta
\]

It is called the {\bf standard ordering in the sense of}~"$<$". Moreover, we define $\leq$ to be the total ordering in the sense of~"$\leq$" on~${\mathbb R}$ obtained from the ordering~$<$ by the method of Lemma~\ref{lemm both orderings}. It is called the {\bf standard ordering in the sense of}~"$\leq$".
\elede

\bproof
Exercise.
\eproof

\bdefi
\label{defi real intervals}
\index{Closed interval}
\index{Interval!closed}
We adopt the convention that all notions related to orderings on~${\mathbb R}$, in particular intervals, refer to the standard ordering in the sense of~"$<$" as defined in Lemma and Definition~\ref{lede definition reals} unless otherwise specified.

Furthermore, for every $x,y \in {\mathbb R}$ with $x < y$ we define
\begin{center}
\begin{tabular}{ll}
$\left] -\infty, y \right] \, = \; \left] -\infty, y \right[ \; \cup \left\{ y \right\}$\,, \quad \quad &
$\left[ x, \infty \right[ \; = \; \left] x, \infty \right[ \; \cup \left\{ x \right\}\,,$\\[.8em]
$\left] x, y \right] \, = \; \left] x, y \right[ \; \cup \left\{ y \right\}$\,, \quad \quad &
$\left[ x, y \right[ \; = \; \left] x, y \right[ \; \cup \left\{ x \right\}$\,,\\[.8em]
$\left[ x, y \right] \, = \; \left] x, y \right[ \; \cup \left\{ x, y \right\}$ &
\end{tabular}
\end{center}

The set~$\left[ x, y \right]$ with $x, y \in {\mathbb R}$ is called {\bf closed interval}. Notice that it is a proper interval with respect to the ordering~$\leq$.

We further agree that the sets $\left[ 0, \infty \right[\;$ and $\;\left] 0, \infty \right[\;$ always refer to subsets of~${\mathbb R}$, or equivalently to the sets ${\mathbb R}_+$ and ${\mathbb R}_+ \!\!\setminus\! \left\{ 0 \right\}$ unless otherwise specified.
\edefi

This convention is in agreement with the ones adopted in the context of natural numbers, Definition~\ref{defi convention ordering natural numbers}, and positive dyadic rational numbers, Definition~\ref{defi convention ordering pos dyadic numbers}. We remark again that, apart from the definition of intervals, it is mostly irrelevant whether the ordering~$<$ or the ordering~$\leq$ on~${\mathbb R}$ is considered, cf.\ Lemmas~\ref{lemm max min inv}, \ref{lemm bound sup inv}, and~\ref{lemm mon invar}.

\midvspace

\begin{lede}
\label{lede add mult reals}
\index{Addition!real numbers}
\index{Sum!real numbers}
\index{Multiplication!real numbers}
\index{Product!real numbers}
\index{Exponentiation!real numbers}
We define two binary functions $+$ (called {\bf addition}) and $\cdot$ (called {\bf multiplication}) on~$\mathbb R$ by
\begin{eqnarray*}
\langle \alpha, \beta \rangle + \langle \gamma, \delta \rangle \!\!\! & = & \!\!\! \langle \alpha + \gamma,\; \beta + \delta \rangle\\
\langle \alpha, \beta \rangle \cdot \langle \gamma, \delta \rangle \!\!\! & = & \!\!\! \langle \alpha \, \gamma + \beta \, \delta,\; \alpha \, \delta + \beta \, \gamma \rangle
\end{eqnarray*}

$\langle \alpha, \beta \rangle + \langle \gamma, \delta \rangle$ and $\langle \alpha, \beta \rangle \cdot \langle \gamma, \delta \rangle$ are called the {\bf sum} and the {\bf product of} $\langle \alpha, \beta \rangle$ {\bf and} $\langle \gamma, \delta \rangle$, respectively. For every $x, y \in {\mathbb R}$ we also write $x \, y$ for $x \cdot y$. We further define that in the absence of brackets products are evaluated before sums. Both addition and multiplication are commutative and associative, and the distributive law
\[
(x + y) \cdot z \, = \, (x \cdot z) + (y \cdot z)
\]

as well as the implications
\begin{eqnarray*}
x < y \; & \Longrightarrow & \; x + z < y + z\\
x < y \;\; \wedge \;\; \langle 0, 0 \rangle < z \; & \Longrightarrow & \; x \, z < y \, z\\
x < y \;\; \wedge \;\; z < \langle 0, 0 \rangle \; & \Longrightarrow & \; y \, z < x \, z
\end{eqnarray*}

hold for $x, y, z \in {\mathbb R}$. For every $x \in {\mathbb R}$ with $x > \langle 0, 0 \rangle$ there is $\alpha \in {\mathbb R}_+$ such that $x = \langle \alpha, 0 \rangle$. Furthermore, for every $x \in {\mathbb R}$ with $x < \langle 0, 0 \rangle$ there is $\alpha \in {\mathbb R}_+$ such that $x = \langle 0, \alpha \rangle$.

Further, let $g : {\mathbb R}_+ \longrightarrow {\mathbb R}$, $g(\alpha) = \langle \alpha, 0 \rangle$, and $h : {\mathbb R}_+ \longrightarrow {\mathbb R}$, $h(\beta) = \langle 0, \beta \rangle$. Then we have
\[
g \left[ {\mathbb R}_+ \right] = \left\{ x \in {\mathbb R} \, : \, x \geq \langle 0, 0 \rangle \right\}\,, \quad \quad h \left[ {\mathbb R}_+ \right] = \left\{ x \in {\mathbb R} \, : \, x \leq \langle 0, 0 \rangle \right\}
\]

The functions $g$ and $h$ are injective. We have 
\benum
\item \label{lede add mult reals 1} $\alpha < \beta \;\; \Longleftrightarrow \;\; g(\alpha) < g(\beta) \;\; \Longleftrightarrow \;\; h(\alpha) > h(\beta)$
\item \label{lede add mult reals 2} $g(\alpha + \beta) = g(\alpha) + g(\beta)$
\item \label{lede add mult reals 3} $h(\alpha + \beta) = h(\alpha) + h(\beta)$
\item \label{lede add mult reals 4} $g(\alpha \, \beta) = g(\alpha) \, g(\beta)$
\eenum

Furthermore, we define the {\bf exponentiation} on the positive subset by
\[
f : g \left[ {\mathbb R}_+ \right]  \times \naturalnumbers \longrightarrow g \left[ {\mathbb R}_+ \right]\,, \quad f ( \langle \alpha, 0 \rangle , m ) = \langle \alpha^m , 0 \rangle
\]

For every $m \in \naturalnumbers$ and $x \in {\mathbb R}$ with $x \geq \langle 0, 0 \rangle$, we also write $x^m$ for~$f(x,m)$. We define the same rules regarding the order of evaluation as for~${\mathbb R}_+$. We have
\[
g \left( \alpha^m \right) = (g(\alpha))^m
\]

for every $m \in \naturalnumbers$ and $\alpha \in {\mathbb R}_+$.
\end{lede}
\vspace{.0in}

\bproof
The proofs that addition and multiplication are well-defined, and that they are commutative and associative, as well as the proof of the distributive law are left as exercise.

To see the three implications, let $\langle \alpha, \beta \rangle,\, \langle \gamma, \delta \rangle,\, \langle \chi, \psi \rangle \in  \mathbb R$ with $\langle \alpha, \beta \rangle < \langle \gamma, \delta \rangle$. Hence we have $\alpha + \delta < \beta + \gamma$.

To show the first implication, notice that $\alpha + \delta + \chi + \psi < \beta + \gamma + \chi + \psi$. It follows that $\langle \alpha, \beta \rangle + \langle \chi, \psi \rangle = \langle \alpha + \chi, \; \beta + \psi \rangle <  \langle \gamma + \chi, \; \delta + \psi \rangle = \langle \gamma, \delta \rangle + \langle \chi, \psi \rangle$.

To show the second implication, assume that $\langle \chi, \psi \rangle > \langle 0, 0 \rangle$. This implies $\psi < \chi$. By Proposition~\ref{prop prod inequ} we obtain
\[
(\alpha + \delta) \, \chi + (\beta + \gamma) \, \psi \, < \, (\alpha + \delta) \, \psi + (\beta + \gamma) \, \chi
\]

and thus
\[
\alpha \, \chi + \beta \, \psi + \gamma \, \psi + \delta \, \chi \, < \, \alpha \, \psi + \beta \, \chi + \gamma \, \chi + \delta \, \psi
\]

It follows that $\langle \alpha, \beta \rangle \cdot \langle \chi, \psi \rangle = \langle \alpha \, \chi + \beta \, \psi, \; \alpha \, \psi + \beta \, \chi \rangle < \langle \gamma \, \chi + \delta \, \psi, \; \gamma \, \psi + \delta \, \chi \rangle = \langle \gamma, \delta \rangle \cdot \langle \chi, \psi \rangle$.

To see the third implication, assume that $\langle \chi, \psi \rangle < \langle 0, 0 \rangle$. This implies $\psi > \chi$. In this case we obtain
\[
\alpha \, \chi + \beta \, \psi + \gamma \, \psi + \delta \, \chi \, > \, \alpha \, \psi + \beta \, \chi + \gamma \, \chi + \delta \, \psi
\]

by Proposition~\ref{prop prod inequ}. This implies $\langle \alpha, \beta \rangle \cdot \langle \chi, \psi \rangle > \langle \gamma, \delta \rangle \cdot \langle \chi, \psi \rangle$.

Now let $\alpha, \beta \in {\mathbb R}_+$. If $\langle \alpha, \beta \rangle > \langle 0, 0 \rangle$, then we have $\beta < \alpha$. Thus there exists $\gamma \in {\mathbb R}_+$ such that $\beta + \gamma = \alpha$ by Proposition~\ref{prop pos real diff}. Therefore we have $\langle \alpha, \beta \rangle = \langle \gamma, 0 \rangle$.

If $\langle \alpha, \beta \rangle < \langle 0, 0 \rangle$, then there is $\gamma \in {\mathbb R}_+$ such that $\alpha + \gamma = \beta$ by Proposition~\ref{prop pos real diff}. Thus $\langle \alpha, \beta \rangle = \langle 0, \gamma \rangle$.

Finally the results (\ref{lede add mult reals 1}) to (\ref{lede add mult reals 4}) clearly follow by definition.
\eproof

Here as in the previous cases the fact that the injection from the positive reals to the reals preserves the orderings and binary functions explains the usage of the same symbols and allows us to deliberately mix the different kinds of numbers in expressions, such as $x + m$, $x \cdot \lfloor m, u \rfloor$, or $x \cdot \alpha$, where $m, u \in \naturalnumbers$, $\alpha \in {\mathbb R}_+$, and $x \in {\mathbb R}$.

The next Proposition, on which the subsequent Lemma is based, is almost obvious though its derivation is a bit lengthy.

\midvspace

\bprop
\label{prop reals pos reals}
Let $g$ and $h$ be defined as in Lemma and Definition~\ref{lede add mult reals}, $B_+ \subset {\mathbb R}_+$ with $B_+ \neq \O$, and ${\mathbb R}_- = \left\{ x \in {\mathbb R} \, : \, x \leq 0 \right\}$. The following statements hold:
\benum
\item \label{prop reals pos reals 1} If $B_+$ has a minimum (maximum), say $\alpha$, then $g(\alpha)$ is a minimum (maximum) of~$g \left[ B_+ \right]$ and $h(\alpha)$ is a maximum (minimum) of~$h \left[ B_+ \right]$.
\item \label{prop reals pos reals 2} Let $L_+$ and $U_+$ be the sets of all lower and upper bounds of $B_+$\,, respectively, $U_g$ the set of all upper bounds of~$g \left[ B_+ \right]$, and $U_h$ the set of all upper bounds of $h \left[ B_+ \right]$, i.e.\
\begin{eqnarray*}
L_+ \!\!\! & = & \!\!\! \big\{ \alpha \in {\mathbb R}_+ : \, \forall \beta \in B_+ \!\setminus\! \left\{ \alpha \right\} \;\; \alpha < \beta \big\}\\
U_+ \!\!\! & = & \!\!\! \big\{ \alpha \in {\mathbb R}_+ : \, \forall \beta \in B_+ \!\setminus\! \left\{ \alpha \right\} \;\; \beta < \alpha \big\}\\
U_g \!\!\! & = & \!\!\! \big\{ x \in {\mathbb R} \, : \, \forall y \in g \left[ B_+ \right] \setminus\! \left\{ x \right\} \;\; y < x \big\}\\
U_h \!\!\! & = & \!\!\! \big\{ x \in {\mathbb R} \, : \, \forall y \in h \left[ B_+ \right] \setminus\! \left\{ x \right\} \;\; y < x \big\}
\end{eqnarray*}

Then we have $g \left[ U_+ \right] = U_g$ and $h \left[ L_+ \right] = U_h \cap {\mathbb R}_-$\,.
\item \label{prop reals pos reals 3} If $B_+$ has a supremum, then $g \left[ B_+ \right]$ has a supremum and we have $\sup g \left[ B_+ \right] = g \left(\sup B_+ \right)$.
\item \label{prop reals pos reals 4} $B_+$ has an infimum, $h \left[ B_+ \right]$ has a supremum, and we have \mbox{$\sup h \left[ B_+ \right] = h \left( \inf B_+ \right)$}.
\eenum

\eprop

\bproof
(\ref{prop reals pos reals 1}) and (\ref{prop reals pos reals 2}) are consequences of Lemma and Definition~\ref{lede add mult reals}~(\ref{lede add mult reals 1}).
\\

\hspace{0.05\textwidth}
\parbox{0.95\textwidth}
{[If $\alpha$ is the minimum of~$B_+$\,, then we have $\alpha \in B_+$ and $\alpha < \beta$ for every $\beta \in B_+ \!\!\setminus\! \left\{ \alpha \right\}$. It follows that $g(\alpha) \in g \left[ B_+ \right]$ and $h(\alpha) \in h \left[ B_+ \right]$. Moreover, we have $g(\alpha) < g(\beta)$ and $h(\beta) < h(\alpha)$ for every $\beta \in B_+ \!\setminus\! \left\{ \alpha \right\}$. Therefore we have $g(\alpha) < x$ for every $x \in g \left[ B_+ \!\setminus\! \left\{ \alpha \right\} \right] = g \left[ B_+ \right] \setminus\! \left\{ g(\alpha) \right\}$. Additionally, we have $x < h(\alpha)$ for every $x \in h \left[ B_+ \!\setminus\! \left\{ \alpha \right\} \right] = h \left[ B_+ \right] \setminus\! \left\{ h(\alpha) \right\}$. The claim in~(\ref{prop reals pos reals 1}) that is stated in brackets is shown similarly.

Further, we obtain~(\ref{prop reals pos reals 2}) as follows:
\begin{eqnarray*}
g \left[ U_+ \right] \!\!\! & = & \!\!\! g \left[ \big\{ \alpha \in {\mathbb R}_+ : \, \forall \beta \in B_+ \!\setminus\! \left\{ \alpha \right\} \;\; \beta < \alpha \big\} \right]\\[.2em]
 & = & \!\!\! \big\{ g(\alpha) \, : \, \alpha \in {\mathbb R}_+\,,\, \forall \beta \in B_+ \!\setminus\! \left\{ \alpha \right\} \;\; \beta < \alpha \big\}\\[.2em]
 & = &  \!\!\! \big\{ g(\alpha) \, : \, \alpha \in {\mathbb R}_+\,,\, \forall \beta \in B_+ \!\setminus\! \left\{ \alpha \right\} \;\; g(\beta) < g(\alpha) \big\}\\[.2em]
 & = &  \!\!\! \big\{ g(\alpha) \, : \, \alpha \in {\mathbb R}_+\,,\, \forall y \in g \left[ B_+ \!\setminus\! \left\{ \alpha \right\} \right] \;\; y < g(\alpha) \big\}\\[.2em]
 & = &  \!\!\! \big\{ g(\alpha) \, : \, \alpha \in {\mathbb R}_+\,,\, \forall y \in g \left[ B_+ \right] \setminus\! \left\{ g(\alpha) \right\} \;\; y < g(\alpha) \big\}\\[.2em]
 & = &  \!\!\! \big\{ x \in {\mathbb R} \, : \, x \geq 0,\, \forall y \in g \left[ B_+ \right] \setminus\! \left\{ x \right\} \;\; y < x \big\}\\[.2em]
 & = &  \!\!\! \big\{ x \in {\mathbb R} \, : \, \forall y \in g \left[ B_+ \right] \setminus\! \left\{ x \right\} \;\; y < x \big\} \, = \; U_g
\end{eqnarray*}
\begin{eqnarray*}
h \left[ L_+ \right] \!\!\! & = & \!\!\! h \left[ \big\{ \alpha \in {\mathbb R}_+ : \, \forall \beta \in B_+ \!\setminus\! \left\{ \alpha \right\} \;\; \alpha < \beta \big\} \right]\\[.2em]
 & = &  \!\!\! \big\{ h(\alpha) \, : \, \alpha \in {\mathbb R}_+\,,\, \forall \beta \in B_+ \!\setminus\! \left\{ \alpha \right\} \;\; \alpha < \beta \big\}\\[.2em]
 & = &  \!\!\! \big\{ h(\alpha) \, : \, \alpha \in {\mathbb R}_+\,,\, \forall \beta \in B_+ \!\setminus\! \left\{ \alpha \right\} \;\; h(\beta) < h(\alpha) \big\}\\[.2em]
 & = &  \!\!\! \big\{ h(\alpha) \, : \, \alpha \in {\mathbb R}_+\,,\, \forall y \in h \left[ B_+ \!\setminus\! \left\{ \alpha \right\} \right] \;\; y < h(\alpha) \big\}\\[.2em]
 & = &  \!\!\! \big\{ h(\alpha) \, : \, \alpha \in {\mathbb R}_+\,,\, \forall y \in h \left[ B_+ \right] \setminus\! \left\{ h(\alpha) \right\} \;\; y < h(\alpha) \big\}\\[.2em]
 & = &  \!\!\! \big\{ x \in {\mathbb R} \, : \, \forall y \in h \left[ B_+ \right] \setminus\! \left\{ x \right\} \;\; y < x \big\} \cap \, {\mathbb R}_- \, = \; U_h \cap {\mathbb R}_-
\end{eqnarray*}]}

Now (\ref{prop reals pos reals 3}) and (\ref{prop reals pos reals 4}) follow by (\ref{prop reals pos reals 1}) and (\ref{prop reals pos reals 2}).
\\

\hspace{0.05\textwidth}
\parbox{0.95\textwidth}
{[If the condition of (\ref{prop reals pos reals 3}) is satisfied, then $\sup B_+$ is the minimum of~$U_+$ as defined in~(\ref{prop reals pos reals 2}). Hence, by~(\ref{prop reals pos reals 1}), $g \left( \sup B_+ \right)$ is the minimum of~$g \left[ U_+ \right]$, which is the supremum of~$g \left[ B_+ \right]$ by~(\ref{prop reals pos reals 2}).

The infimum of~$B_+$ exists by Lemma~\ref{lemm least upper bound pos reals} and Theorem~\ref{theo least upper greatest lower}. Moreover, $\inf B_+$ is the maximum of~$L_+$ as defined in~(\ref{prop reals pos reals 2}). Hence, by~(\ref{prop reals pos reals 1}), $h \left( \inf B_+ \right)$ is the minimum of~$h \left[ L_+ \right]$. This in turn is the minimum of~$U_h$\,, which is the supremum of~$h \left[ B_+ \right]$ by~(\ref{prop reals pos reals 2}).]}
\eproof

\blemm
\label{lemm least upper bound reals}
The ordered space $({\mathbb R},<)$ has the least upper bound property.
\elemm

\bproof
Let $A \subset {\mathbb R}$ such that $A$ has an upper bound and $A \neq \O$. Further let $U$ be the set of all upper bounds of~$A$.

First assume there exists $y \in A$ with $y \geq 0$. We define $B = \left\{ y \in A \, : \, y \geq 0 \right\}$ and $B_+ = g^{-1} \left[ B \right]$ where $g$ is defined as in Lemma and Definition~\ref{lede add mult reals}. Then $B_+ \neq \O$ and $g \left[ B_+ \right] = B$. Moreover, $U$ is the set of all upper bounds of~$B$. Thus also $B_+$ has an upper bound by Proposition~\ref{prop reals pos reals}~(\ref{prop reals pos reals 2}). Therefore $\sup B_+$ exists by Lemma~\ref{lemm least upper bound pos reals}, and $g(\sup B_+) = \sup B = \sup A$ by Proposition~\ref{prop reals pos reals}~(\ref{prop reals pos reals 3}).

Now assume that $y < 0$ for every $y \in A$. We define $A_+ = h^{-1} \left[ A \right]$ where $h$ is defined as in Lemma and Definition~\ref{lede add mult reals}. Then we have $A_+ \neq \O$. By Proposition~\ref{prop reals pos reals}~(\ref{prop reals pos reals 3}), the infimum of~$A_+$ and the supremum of~$h \left[ A_+ \right]$ exist, and we have $h(\inf A_+) = \sup h \left[ A_+ \right] = \sup A$.
\eproof

\blemm
\label{lemm op increasing}
Given $x \in \mathbb R$, the function $f : {\mathbb R} \longrightarrow {\mathbb R}$, $f(y) = y + x$, is strictly increasing. If $x > 0$, then the function $g : {\mathbb R} \longrightarrow {\mathbb R}$, $g(y) = y \, x$, is strictly increasing. If $x < 0$, $g$ is strictly decreasing. For every $m \in \naturalnumbers$, $m \geq 1$, the function $h_m : {\mathbb R}_+ \longrightarrow {\mathbb R}_+$\,, $h_m(y) = y^m$, is strictly increasing. For every $x \in \mathbb R$ and $m \in \naturalnumbers$, $m \geq 1$, the functions $f$, $g$, and $h_m$ are unbounded.
\elemm

\bproof
The fact that $f$ is strictly increasing, and the fact that $g$ is strictly increasing or decreasing under the respective conditions follows by Lemma and Definition~\ref{lede add mult reals}. The fact that $h_m$ is strictly increasing for every $m \in \naturalnumbers$, $m \geq 1$, follows by Lemma and Definition~\ref{lede exp pos reals}.

To see that the functions are unbounded, let $\gamma, \delta, \psi, \chi \in {\mathbb R}_+$ such that $x = \langle \gamma, \delta \rangle$.

For the case of $f$, we define $\alpha = \psi + \delta + 1$ and $\beta = \chi + \gamma$. Then we have $\langle \alpha, \beta \rangle > \langle \psi, \chi \rangle + \langle \delta, \gamma \rangle$. It follows that $\langle \alpha, \beta \rangle + \langle \gamma, \delta \rangle > \langle \psi, \chi \rangle$. 

For the case of~$g$, we assume that $\psi > \chi$. If $\gamma > \delta$, then there is $\varepsilon \in {\mathbb R}_+$ such that $\gamma = \delta + \varepsilon$ by Proposition~\ref{prop pos real diff}. We may choose $\zeta \in {\mathbb R}_+$ such that $\varepsilon \, \zeta > \psi$ by Corollary~\ref{coro prod greater}. We define $\alpha = \gamma + \zeta$ and $\beta = \gamma$. It follows that
\begin{eqnarray*}
 && \!\! \alpha \, \gamma + \beta \, \delta + \chi\\
& = & \!\! \gamma \, \gamma + \zeta \, \gamma + \gamma \, \delta + \chi\\
& = & \!\! \gamma \, \gamma + \zeta \, \delta + \zeta \, \varepsilon + \gamma \, \delta + \chi\\
& > & \!\! \gamma \, \gamma + \zeta \, \delta + \psi + \gamma \, \delta + \chi\\
& \geq & \!\! \gamma \, \gamma + \zeta \, \delta + \gamma \, \delta + \psi\\
& = & \!\! \alpha \, \delta + \beta \, \gamma + \psi
\end{eqnarray*}

Hence $\langle \alpha, \beta \rangle \cdot \langle \gamma, \delta \rangle = \langle \alpha \, \gamma + \beta \, \delta, \; \alpha \, \delta + \beta \, \gamma \rangle > \langle \psi, \chi \rangle$. If $\gamma < \delta$, then there is  $\varepsilon \in {\mathbb R}_+$ such that $\delta = \gamma + \varepsilon$. Then we may again choose $\zeta \in {\mathbb R}_+$ such that $\varepsilon \, \zeta > \psi$. We define $\alpha = \gamma$ and $\beta = \gamma + \zeta$. It follows that
\begin{eqnarray*}
 && \!\! \alpha \, \gamma + \beta \, \delta + \chi\\
& = & \!\! \gamma \, \gamma + \gamma \, \delta + \zeta \, \delta + \chi\\
& = & \!\! \gamma \, \gamma + \gamma \, \delta + \zeta \, \gamma + \zeta \, \varepsilon + \chi\\
& > & \!\! \gamma \, \gamma + \gamma \, \delta + \zeta \, \gamma + \psi\\
& = & \!\! \alpha \, \delta + \beta \, \gamma + \psi
\end{eqnarray*}

also in this case.

Let $m \in \naturalnumbers$ with $m \geq 1$. To see that $h_m$ is unbounded, let $\beta \in {\mathbb R}_+$. We may choose $p \in \naturalnumbers$ such that $p > 1$ and $p > \beta$. It follows that $p^m \geq p > \beta$ by Lemma and Definition~\ref{lede exponentiation natural numbers}.
\eproof

\blede
The triple $({\mathbb R},+,0)$ is an Abelian group. Let $x, y \in {\mathbb R}$ and $\alpha, \beta \in {\mathbb R}_+$ such that $x = \langle \alpha, \beta \rangle$. The inverse of~$x$ with respect to addition is given by $\langle \beta, \alpha \rangle$ and denoted by~$-x$. We also write $y - x$ for $y + (-x)$. In the absence of brackets we define the following priorities:
\[
- x + y = (-x) + y, \quad -x \, y = -(x \, y), \quad - x^m = - \left( x^m \right)
\]

\elede

\bproof
The addition is associative and commutative by Lemma and Definition \ref{lede add mult reals}. The other assertions are clear.
\eproof

\brema
Let $\alpha, \beta \in {\mathbb R}_+$ and $x, y \in {\mathbb R}$. We have
\benum
\item $\langle \alpha, \beta \rangle = \alpha - \beta$
\item If $\alpha > \beta$, then $(\alpha - \beta) \in {\mathbb R}_+$.
\item $0 - x = -x$
\item $(-1) \cdot x = -x$
\item $x < y \;\; \Longrightarrow \;\; -y < -x$
\item $\left( x > 0 \right) \; \wedge \; \left( y > 0 \right) \;\; \Longrightarrow \;\; x \, y > 0$
\item $\left( x > 0 \right) \; \wedge \; \left( y < 0 \right) \;\; \Longrightarrow \;\; x \, y < 0$
\eenum

\erema

\blede
The triple $({\mathbb R} \!\setminus\! \{ 0 \},\cdot,1)$ is an Abelian group. Let $x, y \in {\mathbb R} \!\setminus\! \{ 0 \}$. The inverse of $x$ with respect to multiplication is denoted by $x^{-1}$, $\left( 1 / x \right)$, or~$\left( \frac{\textstyle 1}{\textstyle x} \right)$. We also write $\left( x / y \right)$ or $\left( \frac{\textstyle x}{\textstyle y} \right)$ for $x \cdot (1 / y)$. In the absence of brackets the superscript~"$-1$" is evaluated before sums and products, and we define $-x^{-1} = - \left(x^{-1} \right)$. We have $(-x)^{-1} = -\left( x^{-1}\right)$.
\elede

\bproof
The multiplication is associative and commutative by Lemma and Definition~\ref{lede add mult reals}. For every $\alpha, \beta \in {\mathbb R}_+$ we have
\[
\langle \alpha, \beta \rangle \cdot \langle 1, 0 \rangle \, = \, \langle \alpha, \beta \rangle
\]

For $x \in {\mathbb R} \!\setminus\! \{ 0 \}$ we define
\[
1 / x = \left\{ \begin{array}{ll}
\left\langle {\displaystyle \frac{1}{\alpha - \beta}}, \, 0 \right\rangle & \; \mbox{if \, $\alpha > \beta$}\\[1.5em]
\left\langle 0,\, {\displaystyle \frac{1}{\beta - \alpha}} \right\rangle & \; \mbox{if \, $\beta > \alpha$}
\end{array} \right.
\]

where $\alpha, \beta \in {\mathbb R}_+$ such that $x = \langle \alpha, \beta \rangle$. The inverses on the right hand side are defined according to Theorem~\ref{theo pos reals group}. Notice that this definition is independent of the specific choice of~$\alpha$ and~$\beta$. Then $x \cdot (1/x) = 1$, that is $(1 / x)$ is the inverse of~$x$. The last claim is clear.
\eproof

\bdefi
\index{Dyadic rational numbers}
\index{Numbers!dyadic rational}
The members of the set
\[
{\mathbb D} = \big\{ x \in {\mathbb R} \, : \, x \in {\mathbb D}_+ \, \vee \, -x \in {\mathbb D}_+ \big\}
\]

are called {\bf dyadic rational numbers}.
\edefi

\brema
\label{rema dyadic countable}
${\mathbb D}$ is countable by Corollary~\ref{coro dyadic countable} and Lemma~\ref{lemm count times count}.
\erema

\blemm
\label{lemm D dense}
We have
\benum
\item \label{lemm D dense 1} ${\mathbb D}$ is $<$-dense in~${\mathbb R}$.
\item \label{lemm D dense 2} ${\mathbb R} \!\setminus\! {\mathbb D}$ is $<$-dense in~${\mathbb R}$.
\item \label{lemm D dense 3} ${\mathbb R}$ is $<$-dense.
\eenum
\elemm

\bproof
(\ref{lemm D dense 1}) follows by Lemma~\ref{lemm D0 dense}.

(\ref{lemm D dense 2}) follows by Lemma~\ref{lemm dyadic compl dense}.

(\ref{lemm D dense 3}) is a consequence of~(\ref{lemm D dense 1}).
\eproof

\bdefi
\index{Absolute value}
We define the function $b : {\mathbb R} \longrightarrow {\mathbb R}_+$ by
\[
b(x) = \left\{\begin{array}{ll}
x & \; \mbox{if \, $x \geq 0$}\\[.5em]
-x & \; \mbox{if \, $x < 0$}
\end{array}\right.
\]

We also write $|x|$ for $b(x)$. $|x|$ is called {\bf absolute value of}~$x$.
\edefi

\brema
\label{rema prop abs val}
The function $b$ is clearly surjective. Moreover, we have
\[
|x|^2 = x^2, \quad \quad |x + y| \leq |x| + |y|\,, \quad \quad |x \, y| = |x| \, |y|
\]

for every $x, y \in {\mathbb R}$.
\erema

The following result is applied in the proof of Lemma~\ref{lemm dist cont}.

\midvspace

\blemm
\label{lemm inf const}
Let $X$ be a set, $f : X \longrightarrow {\mathbb R}$ a map, and $c \in {\mathbb R}$. Then
\[
\inf_{x \in X} \big( f(x) + c \big) \, = \, \inf_{x \in X} f(x) \, + \, c
\]

\elemm

\bproof
Let $a$ be a lower bound of $\left\{ f(x) + c \, : \, x \in X \right\}$, i.e.\ we have $a \leq f(x) + c$ for every $x \in X$. It follows that $a - c \leq f(x)$ for every $x \in X$, i.e.\ $(a - c)$ is a lower bound of $\left\{ f(x) \, : \, x \in X \right\}$. Thus we have $a - c \leq \inf_{x \in X} f(x)$, and hence
\[
\inf_{x \in X} \big( f(x) + c \big) \, \leq \, \inf_{x \in X} f(x) \, + \, c
\]

Applying this result to $-c$ instead of $c$, and to the function
\[
g : X \longrightarrow {\mathbb R},\quad  g(x) = f(x) + c
\]

we obtain the reverse inequality.
\eproof

We conclude this Section with some examples of orderings and functions involving the real numbers.

\midvspace

\bexam
We recall Example~\ref{exam product pre-orderings}: Let $(X_i,R_i)$ ($i \in I$) be pre-ordered spaces, where $I$ is an index set, and $X = \bigtimes_{\!\! i \in I}\, X_i$. Then ${\mc R} = \left\{ p_i^{-1} \left[ R_i \right] \, : \, i \in I \right\}$ is a system of pre-orderings on~$X$. 

Now, if $(X_i,R_i) = (\mathbb{R},<)$ ($i \in I$), then the members of~$\mc R$ are orderings in the sense of~"$<$". However, they are not total orderings unless $I$ is a singleton. Clearly $\mc R$ is independent.
\eexam

\bexam
\label{exam rn ordering}
Let $(X_i,R_i)$ ($i \in I$) be pre-ordered spaces, where $I$ is an index set, and $X = \bigtimes_{\!\! i \in I}\, X_i$. Then $S = \bigcap \left\{ p_i^{-1} \left[ R_i \right] \, : \, i \in I \right\}$ is a pre-ordering on~$X$ (cf.\ Example~\ref{exam product section ordering}). Now let $n \in \naturalnumbers$, $n > 0$, and $I = \sigma(n) \!\setminus\! \left\{ 0 \right\}$. If $(X_k,R_k) = (\mathbb{R},<)$ ($k \in \naturalnumbers$, $1 \leq k \leq n$), then $S$ is an ordering in the sense of~"$<$". However, it is not a total ordering unless $n = 1$. For $x, y \in \mathbb{R}^n$ we have $x < y$ iff $x_k < y_k$ ($1 \leq k \leq n$). For the same~$I$, if $(X_k,R_k) = (\mathbb{R},\leq)$ ($1 \leq k \leq n$), then $S$ is an ordering in the sense of~"$\leq$" since the ordering $\leq$ on $\mathbb{R}$ is antisymmetric and $\left\{ p_i \, : \, i \in I \right\}$ distinguishes points. However the ordering $\leq$ on $\mathbb{R}^n$ is not a total ordering unless $n = 1$. For $x, y \in \mathbb{R}^n$ we have $x \leq y$ iff $x_k \leq y_k$ ($1 \leq k \leq n$). Both $<$ and $\leq$ on $\mathbb{R}^n$ have full range and full domain. They are not connective unless $n = 1$. Since $\mathbb{D}$ is dense in~$\mathbb{R}$, $\mathbb{D}^n$ is dense in $\mathbb{R}^n$ with respect to both orderings.
\eexam

\bexam
The pair $(\mathbb{R},\leq)$, where $\leq$ denotes the standard ordering in the sense of~"$\leq$", is a pre-ordered space. The function $f : \mathbb{R} \longrightarrow \mathbb{R}$ that maps every real number $x$ to the smallest integer greater or equal than~$x$ is $\leq$-increasing and projective. Similarly, the function $g : \mathbb{R} \longrightarrow \mathbb{R}$ that maps every real number $x$ to the smallest {\it even} integer greater or equal than~$x$ is $\leq$-increasing and projective. Thus also the composition $g \circ f$ is $\leq$-increasing.
\eexam


\backmatter


\cleardoublepage
\phantomsection
\addcontentsline{toc}{chapter}{Index}
\printindex

\end{document}